\documentclass[reqno,12pt]{amsart}
 \usepackage{amsmath}
 \usepackage{amsthm}
 \usepackage{amssymb}
 \usepackage{latexsym,longtable}
 \usepackage{graphicx}
 \usepackage{multicol}
 \usepackage{mathrsfs}
 \usepackage{young}
 \usepackage[vcentermath,enableskew,stdtext]{youngtab}
 \usepackage[left=1.1in,right=1.1in,top=1.06in,bottom=1.1in]{geometry}
 \usepackage[all]{xy}
    \SelectTips{cm}{10}     
    \everyxy={<2.5em,0em>:} 
 \usepackage{fancyhdr}      
 \usepackage{multicol}
 \usepackage{multirow}
 \usepackage{enumitem}
 \usepackage{bm}
 \usepackage{stmaryrd}
 \usepackage{etex}

\usepackage{tabmac}
\usepackage{ytableau}
\usepackage{tikz}
\usepackage{tkz-euclide}
\usetikzlibrary{backgrounds}
\usetikzlibrary{arrows}

 \usepackage{float}
 \usepackage{array}
 \usepackage{colortbl}
 \usepackage{mathtools}
 \restylefloat{figure}
\numberwithin{equation}{section}
 \usetikzlibrary{shapes}
 \usetikzlibrary{arrows}

 \usetikzlibrary{snakes}
\usepackage{url}

\usepackage{amsfonts,epsfig,color}
\makeatletter
\newcommand*{\centerfloat}{%
  \parindent \z@
  \leftskip \z@ \@plus 1fil \@minus \textwidth
  \rightskip\leftskip
  \parfillskip \z@skip}
\makeatother

\theoremstyle{plain}
\newtheorem{theorem}{Theorem}[section]
\newtheorem{lemma}[theorem]{Lemma}
\newtheorem{corollary}[theorem]{Corollary}
\newtheorem{proposition}[theorem]{Proposition}
\newtheorem{conjecture}[theorem]{Conjecture}

\theoremstyle{definition}
\newtheorem{definition}[theorem]{Definition}

\newtheorem{remark}[theorem]{Remark}
\newtheorem{example}[theorem]{Example}

\newtheorem{question}[theorem]{Question}

\newcommand{\ignore}[1]{}

\renewcommand{\a}{\ensuremath{\mathbf{a}}}
\newcommand{\aut}{\text{\rm Aut}\,}
\newcommand{\A}{{\ensuremath{\mathcal{A}}}}

\renewcommand{\b}{\ensuremath{\mathbf{b}}}
\renewcommand{\c}{\ensuremath{\mathbf{c}}}
\renewcommand{\d}{\ensuremath{\mathbf{d}}}

\newcommand{\gl}{\ensuremath{\mathfrak{gl}}}

\newcommand{\PP}{\ensuremath{\mathbb{P}}}

\newcommand{\QQ}{\ensuremath{\mathbb{Q}}}

\newcommand{\sgn}{\text{\rm sgn}}

\newcommand{\U}{\ensuremath{\mathcal{U}}}

\newcommand{\ZZ}{\ensuremath{\mathbb{Z}}}

\newcommand{\be}{\begin{equation}}
\newcommand{\ee}{\end{equation}}

\renewcommand{\SS}{\ensuremath{\mathcal{S}}}
\newcommand{\tsr}{\ensuremath{\otimes}}

\newcommand{\creading}{\text{\rm colword}}

\newcommand{\sh}{\text{\rm sh}}

\def\Tiny{\fontsize{6pt}{6pt}\selectfont}

\newcommand{\rev}{\ensuremath{{\text{\rm rev}}}}

\newcommand{\Iplacp}[1]{\ensuremath{{{I_\text{\rm plac}^{#1}}}}}
\newcommand{\Icommp}[1]{\ensuremath{{{I_{\text{\rm C}}^{#1} }}}}
\newcommand{\Iplac}{\ensuremath{{{I_\text{\rm plac}}}}}
\newcommand{\Ipol}{\ensuremath{I_\text{\rm pol}}}  

\newcommand{\inv}{\ensuremath{{\text{\rm inv}}}}
\newcommand{\Des}{\ensuremath{{\text{\rm Des}}}}

\newcommand{\Wo}[1]{\ensuremath{{\text{\rm W$_{#1}$}}}}

\newcommand{\e}{\mathsf}
\newcommand{\diagread}{\text{\rm diagread}}

\newcommand{\SYT}{\text{\rm SYT}}

\newcommand{\SSYT}{\text{\rm SSYT}}

\newcommand{\KT}{\text{\rm KT}}
\newcommand{\LT}{\text{\rm LT}}

\newcommand{\ver}{\text{\rm Vert}}

\newcommand{\spa}{\hspace{.3mm}}
\newcommand{\dotto}{\dashrightarrow}  
\DeclareMathOperator{\dote}{\ensuremath{\text{\spa-\spa-\spa-}}}
\newcommand{\leftto}{\leftarrow} 
\newcommand{\inc}{\text{\rm inc}}
\newcommand{\asc}{\text{\rm asc}}
\newcommand{\Asc}{\text{\rm Asc}}
\newcommand{\defeq}{\overset{\text{{\em def}}}{=}}
\newcommand{\bx}[2]{{\boldsymbol {#1}[#2]}}
\newcommand{\threeone}{\ensuremath{(\mathbf{3}+\mathbf{1})}}
\newcommand{\twotwo}{\ensuremath{(\mathbf{2}+\mathbf{2})}}
\newcommand{\twoone}{\ensuremath{(\mathbf{2}+\mathbf{1})}}

\usetikzlibrary{positioning}
\usetikzlibrary{arrows.meta}
\tikzstyle{mystealthbig} = [-{Latex[length=2.1mm, width=1.53mm]}]
\tikzstyle{mystealth} = [-{Latex[length=1.8mm, width=1.34mm]}]
\tikzstyle{mystealthsmall} = [-latex]

\tikzstyle{vertex}=[circle, fill, inner sep=0pt, minimum size=4pt, outer sep = .3pt]
\newcommand{\vertex}{\node[vertex]}

\newcommand{\Jnot}[1]{\ensuremath{{\hspace{1pt}\text{\rm:}~\raisebox{1pt}{$#1$}~\text{\rm:}\hspace{1pt}}}}

\title{Noncommutative Schur functions for posets}

\author{J. Blasiak}
   \author{H. Eriksson}
   \author{P. Pylyavskyy}
   \author{I. Siegl}

   \address[Blasiak]
   {Dept.\ of Mathematics\\
    Drexel University\\
    Philadelphia, PA}
   \email{jblasiak@gmail.com}

   \address[Eriksson]
   {Dept.\ of Mathematics\\
    Drexel University\\
    Philadelphia, PA}
   \email{hae34@drexel.edu}

   \address[Pylyavskyy]
    {Dept.\ of Mathematics\\
     University of Minnesota\\
     Minneapolis, MN}
   \email{pylyavskyy@gmail.com}

   \address[Siegl]
   {Dept.\ of Mathematics\\
    University of Washington\\
    Seattle, WA}
   \email{isaiahsiegl@gmail.com}

   \thanks{Authors were supported by NSF Grants DMS-1855784 and DMS-2154282 (J.~B.)
    and DMS-1949896 (P.~P.).}

\begin{document}

\begin{abstract}
The machinery of noncommutative Schur
functions is a general approach to
Schur positivity of symmetric functions initiated by Fomin-Greene \cite{FG}.
Hwang \cite{Hwang} recently adapted this theory to posets to give a new approach to the
Stanley-Stembridge conjecture.
We further develop this theory, incorporating
the $P$-Knuth equivalences introduced by Kim and the third author \cite{KPchromatic},
and prove that the symmetric function of any  $P$-Knuth equivalence graph is Schur positive.
This settles the main conjecture in \cite{KPchromatic} and refines results of Gasharov  \cite{Gasharov}, Shareshian-Wachs \cite{SWchromatic}, and Hwang \cite{Hwang}
 on the Schur positivity of chromatic symmetric functions.
\end{abstract}%

\maketitle

\section{Introduction}
\label{s intro}

The chromatic symmetric function  of a  graph $G$
is the sum
$
X_G(\mathbf{x})  = \sum_\kappa \prod_{v \in \ver(G)} x_{\kappa(v)}
$
where $\kappa$ ranges over all proper colorings  $G$ with colors taken from the set of positive integers.
It can be regarded as a generalization of the chromatic polynomial $\chi_G(k)$ of  $G$ since the specialization of $X_G$ at  $x_1=\cdots = x_k = 1$,  $x_{k+1}= x_{k+2}= \cdots = 0$ is $\chi_G(k)$.

Studies of chromatic symmetric functions have centered around the Stanley-Stembridge conjecture {\cite{Stanleychromatic, StanleyStembridge}}:
\begin{conjecture}
\label{cj SS intro}
If $G = \inc(P)$ is the incomparability graph of a \threeone-free poset  $P$,
then $X_G(\mathbf{x})$ is a positive sum of elementary symmetric functions  $e_\lambda(\mathbf{x})$ (or  $e$-positive for short).
\end{conjecture}

This conjecture originated in studies of immanants of Jacobi-Trudi matrices
of Goulden-Jackson \cite{GouldenJackson}, Greene \cite{GreeneImmanant}, and Stembridge  \cite{StembridgeImmanant}, and Haiman's work \cite{Himmanant} uncovered deep connections to Kazhdan-Lusztig theory.

This area has become quite active in recent years.
Shareshian-Wachs \cite{SWchromatic}
introduced a $t$-version  $X_{\inc(P)}(\mathbf{x},t)$ of chromatic symmetric functions (see \S\ref{ss t chromatic})
for any natural unit interval order $P$, and
conjectured that they encode characters of a symmetric group action on the cohomology of regular semisimple Hessenberg varieties.
This conjecture was proven several years later by Brosnan and Chow \cite{bc18} and independently by Guay-Paquet \cite{g16}.
Alexandersson-Panova discovered connections to LLT polynomials \cite{APchromatic} and Clearman-Hyatt-Shelton-Skandera \cite{CHSS} further developed connections to Kazhdan-Lusztig theory. 
Shareshian-Wachs \cite{SWchromatic} also conjectured a $t$-refinement of the Stanley-Stembridge conjecture:

\begin{conjecture}
\label{cj SW intro}
The $t$-chromatic symmetric functions $X_{\inc(P)}(\mathbf{x},t)$ are
$e$-positive.
\end{conjecture}

Partial progress towards Conjectures \ref{cj SS intro} and \ref{cj SW intro}
includes results of
Stembridge \cite{StembridgeImmanant}, Wolfgang \cite{Wolfgangthesis}, Clearman-Hyatt-Shelton-Skandera \cite{CHSS}, and Hwang \cite{Hwang}
establishing positivity of the coefficient of  $e_\lambda$ for several classes of shapes $\lambda$,
and work by Gebhard-Sagan \cite{gebhard2001chromatic},  Dahlberg-van Willigenburg \cite{dahlberg2018lollipop},
Harada-Precup \cite{harada2019cohomology},
Cho-Huh \cite{cho2019positivity},
Dahlberg \cite{Dahlberg}, and
Cho-Hong \cite{ChoHong} proving $e$-positivity of  $X_G$ for various classes of graphs  $G$.


Towards Conjecture \ref{cj SS intro}, a theorem of Haiman \cite{Himmanant} established the weaker statement that
$X_{\inc(P)}(\mathbf{x})$ is Schur positive.
Gasharov \cite{Gasharov} gave an explicit formula for the Schur expansion of  $X_{\inc(P)}(\mathbf{x})$ in terms of an adaptation of tableaux known as  $P$-tableaux (Definition \ref{def:poset SSYT}).

Many have sought to adapt combinatorics such as the RSK-correspondence, jeu-de-taquin, and crystals to posets to explain and strengthen Gasharov's result.  Sundquist-Wagner-West \cite{SWW}, Chow \cite{ChowDescents},
and Kim and the third author \cite{KPchromatic} develop versions of the Robinson-Schensted correspondence for certain unit interval orders, which refines Gasharov's result for these posets.  In a similar vein, Ehrhard's recent work \cite{Ehrhard} constructs a crystal on  $P$-arrays whose components are not isomorphic to usual $\gl_n$-crystals but
do have Schur positive  generating functions.

Building off of work of \cite{KPchromatic}, we define a version of a $P$-Knuth equivalence graph
for any \threeone-free poset  $P$ whose vertices are words in the alphabet  $P$ (Definition \ref{def-switchboard}).  
To any $P$-Knuth equivalence graph  $\Gamma$, we associate
\begin{align}
\label{eq:F Gamma def intro}
F_\Gamma(\mathbf{x})=\displaystyle\sum_{\e{w} \in \ver(\Gamma)}
Q_{\Des_P(\e{w})}(\mathbf{x}),
\end{align}
where $\Des_P(\e{w})= \{i  : \e{w}_i >_P \e{w}_{i+1} \}$
and $Q_{\Des_P(\e{w})}(\mathbf{x})$ is Gessel's fundamental quasisymmetric function \cite{GesselPPartition}.

We resolve the main conjecture in \cite{KPchromatic} and generalize it from the setting of natural unit interval orders to \threeone-free posets:

\begin{theorem}
\label{t P KE schur pos intro}
For any  $P$-Knuth equivalence graph  $\Gamma$,
$F_\Gamma(\mathbf{x})$ is a Schur positive symmetric function.
The coefficient of $s_\lambda(\mathbf{x})$ in  $F_\Gamma(\mathbf{x})$ is
the number of  $P$-tableaux $T$ such that
the column reading word of $T$ appears as a vertex in~$\Gamma$.
\end{theorem}

This further refines the  refinements of Gasharov's result established by Shareshian-Wachs \cite[Theorem 6.3]{SWchromatic} and Hwang \cite[Theorem 4.19]{Hwang}.

The starting point for our proof of Theorem \ref{t P KE schur pos intro} is the theory of \emph{noncommutative Schur functions}, a general technique for studying Schur positivity.
Drawing on ideas of Lascoux and Sch\"utzenberger
\cite{LS, Sch}, Fomin and Greene \cite{FG} developed this theory 
to give positive formulae for the Schur expansions of a large class of symmetric functions
that includes the Stanley symmetric functions and stable Grothendieck polynomials.
More recently, Fomin,  Liu, and the first author \cite{BLamLLT,BF,BLKronecker} further developed this theory, incorporating ideas of
Assaf \cite{SamiOct13, SamiForum} and Lam \cite{LamRibbon},
and used it to prove Schur positivity of LLT polynomials indexed by a 3-tuple of skew shapes and to give a new proof of a formula for Kronecker coefficients when one of the shapes is a hook.

The technique
adapts various notions from the theory of symmetric functions such as Schur functions and
Cauchy identities to versions in which these objects belong to an algebra $\U$ of noncommuting variables
 $u_1,u_2,\dots$.
It relates the Schur positivity of ordinary symmetric functions to the problem of finding a  monomial positive expression for a noncommutative version of a Schur function  $\mathfrak{J}_\lambda \in \U$ modulo a carefully chosen ideal.

Hwang \cite{Hwang} recently adapted this theory to posets to give a new approach to the
Stanley-Stembridge conjecture.  We further develop this theory
to prove Theorem \ref{t P KE schur pos intro}
and give additional insights into Hwang's approach.
More specifically, we

\begin{itemize}[leftmargin=.86cm]
\item 
show that $P$-analogs of elementary symmetric functions commute in a  $P$-version  of the plactic algebra      $\U_P/\Iplacp{P}$  when  $P$ is a \threeone-free poset (Theorem \ref{th:IplacP elem commute}).
\item reduce Theorem \ref{t P KE schur pos intro} to the problem of showing that a noncommutative  $P$-version of  a Schur function
$\mathfrak{J}^P_\lambda$ has a  monomial positive expression mod  $\Iplacp{P}$ (\S\ref{ss P descents}--\ref{ss P Schur etc}), and
then prove this by an inductive argument on a flagged version of the $\mathfrak{J}^P_\lambda$ (\S\ref{sec:proof-of-IplacP-positivity}).
\item study a family of graphs introduced by Hwang \cite{Hwang} called  $H$-graphs obtained by adding additional edges to $P$-Knuth equivalence graphs; we conjecture that their symmetric functions as in
    \eqref{eq:F Gamma def intro} are  $h$-positive
    for any  \threeone-free poset  $P$, slightly strengthening a conjecture of \cite{Hwang} (\S\ref{ss e pos}).
\item give, for  $\lambda$ a hook or two-column shape,  a formula for the coefficient of  $h_\lambda$ in  $\omega\spa X_{\inc(P)}(\mathbf{x},t)$ and more generally for the symmetric function of any $H$-graph,
    giving variations on results of \cite{Hwang}, \cite{CHSS}, \cite{Wolfgangthesis}
     (\S\ref{ss e pos}).
   \item
introduce an algebra called the \emph{arrow algebra},
  which encompasses all the  $\U_P$ simultaneously in a tidy way
and suggests further strengthenings of Conjectures~\ref{cj SS intro}~and~\ref{cj SW intro} as well as variations
on Stanley's  $P$-partitions conjecture
(\S\ref{sec:poset algebra}).
\end{itemize}

\subsection{Comparison with results of Hwang}
In the final stages of writing this paper, we learned of the recent article \cite{Hwang}
which
has significant overlap with this work.
Essentially all of \S\ref{s sec 2} and much of \S\ref{ss e pos}
appears in \cite{Hwang}, so this material can now be regarded as expository
(however see the beginnings of these (sub)sections and Remark \ref{r how generalize} (a) for an account of some slight differences in our viewpoints).
Our Theorem \ref{t P KE schur pos intro} is new and can be regarded as a strengthening of
\cite[Theorem 4.19]{Hwang}.  
Our Theorem \ref{th:IplacP elem commute}, which shows that $P$-versions of elementary symmetric functions commute
modulo  $\Iplacp{P}$, is a strengthening of \cite[Theorem 4.2]{Hwang} which establishes the same for the larger ideal  $I_H^P$.

Much of our paper takes place in the setting of \threeone-free posets whereas Hwang mostly works in
setting of natural unit interval orders.  Though Guay-Paquet \cite{GPchromatic} reduced the Stanley-Stembridge conjecture to proving $e$-positivity for natural unit interval orders,
the additional generality considered here may still be significant since (1) it is not clear that the reduction
extends to the refinements of Schur and  $e$-positivity of  $X_{\inc(P)}$ studied in our papers,
and (2)
it is desirable to have positive formulas which are uniform over all \threeone-free posets, and such formulas are more readily obtained from our approach than via the Guay-Paquet reduction.

\section{Noncommutative  Schur and  monomial functions for posets}
\label{s sec 2}

In this section, we present Hwang's adaptation \cite{Hwang} of the basic setup of \cite{FG,BF} to posets.
We will then apply this general tool in \S\ref{s P Knuth and H graphs} to study the Schur and  $e$-positivity of (refinements of) chromatic symmetric functions.

Our presentation closely follows the material in \cite[(2.1) to (2.15)]{BF}.
Our results are stated in a more general setting than the corresponding ones in \cite{Hwang}, but
it is easily seen that the proofs therein apply to the present setting as well.

Throughout this paper, a poset will always mean a finite poset.

\subsection{Elementary and complete homogeneous functions for posets}
\label{ss poset elementary symmetric functions}

Let  $P$ be a poset of size $|P| = N$.
For  $a,b \in P$, we write
\begin{align}
\label{e atob def1}
& \text{$a \to_P b$  \ \ \ \spa \spa if   $a$ is less than  $b$ in  $P$,} \\
\label{e atob def2}
& \text{$a \dote_P b$ \ \  \ \spa \spa if  $a$ and  $b$ are incomparable or equal elements of  $P$}, \\
\label{e atob def3}
& \text{$a \dotto_P b$ \ \  if  $a \dote_P b$ or $a \to_P b$.}
\end{align}
We also write  $b \leftto_P a$ for  $a \to_P b$ and omit the subscript  $P$ when it is clear.

Let $\U_P= \ZZ\langle u_a : a \in P \rangle$ be the free associative ring
 with generators $\{ u_a : a \in P\}$, which we regard as noncommuting variables.

For $S\subset P$ and $k\in\mathbb{Z}$, we define
the \emph{noncommutative $P$-elementary function} $e_k^P(\mathbf{u}_S)   \in  \U_P$ by
\begin{align}\label{e ek def}
e_k^P(\mathbf{u}_S)=\sum_{\substack{a_1 \leftto_P a_2 \leftto_P \cdots \leftto_P a_k \\ a_1,\dots,a_k \in S}}u_{a_1}u_{a_2}\cdots u_{a_k}\,;
\end{align}
by convention,
$e_0(\mathbf{u}_S)=1$ (even when  $S = \varnothing$) and $e_k^P(\mathbf{u}_S) = 0$ for $k<0$ or $k>|S|$.
For $S= P$, we use the notation $e_k^P(\mathbf{u})$.

The \emph{noncommutative $P$-complete homogeneous functions} are defined by
\[
h^P_\ell(\mathbf{u})=\sum_{ a_1 \spa \dotto_P \spa a_2 \spa \dotto_P \spa \cdots \spa \dotto_P \spa a_\ell}u_{a_1}u_{a_2}\cdots u_{a_\ell};
\]
by convention, $h^P_0(\mathbf{u})=1$
and $h^P_{\ell}(\mathbf{u}) = 0$ for $\ell<0$.

\begin{definition} 
Let  $\Icommp{P}$ denote the (two-sided) ideal of  $\U_P$ generated by
\begin{equation}
\label{eq:e's-commute}
e^P_k(\mathbf{u}_S)\,e^P_\ell(\mathbf{u}_S) - e^P_\ell(\mathbf{u}_S)\,e^P_k(\mathbf{u}_S)
\end{equation}
over all $k,\ell$ and $S \subset P$.
\end{definition}

Also define $\Ipol^P$ to be the ideal of  $\U_P$ generated by
\begin{align}
  u_a u_b - u_b u_a \ \ \ \ \text{ for all }  a,b \in P,
\end{align}
so that  $\U_P/\Ipol^P$ is just the polynomial ring in the variables  $\{u_a : a \in P\}$.
Since the  $e_k^P(\mathbf{u})$'s commute in  $\U_P/\Ipol^P$, we have  $\Ipol^P \supset \Icommp{P}$.

The most natural setting for our basic setup takes place in  $\U_P/\Icommp{P}$.
However all of these results will carry over to further quotients  $\U_P/I$ for  $I \supset \Icommp{P}$,
 $\Ipol^P$ being the simplest example.
Later in \S\ref{s P Knuth and H graphs} we will study three other ideals sandwiched between  $\Ipol^P$ and  $\Icommp{P}$, which will
be the setting where we study various positivity questions such as the Stanley-Stembridge conjecture.

\begin{proposition}
\label{cor:h's-commute}
Let $I\subset\U_P$ be an ideal containing~$\Icommp{P}$.
Then the elements $h^P_1(\mathbf{u}), h^P_2(\mathbf{u}),\dots$
commute pairwise modulo~$I$.
\end{proposition}

\begin{proof}
For any~$m>0$, the following identity holds in  $\U_P$\,:
\begin{align}\label{e E(-x)H(x)}
h^P_m(\mathbf{u})-h^P_{m-1}(\mathbf{u})\,e^P_1(\mathbf{u})+h^P_{m-2}(\mathbf{u})\,e^P_2(\mathbf{u})-\cdots
+(-1)^m e^P_m(\mathbf{u})=0.
\end{align}
To prove this, note that any monomial $u_{a_1} \cdots u_{a_m}$ such that  $a_1 \dotto a_2 \dotto \cdots \dotto a_{s} \leftto a_{s+1} \leftto \cdots \leftto a_m$ for some  $s$  appears in exactly two of the summands in \eqref{e E(-x)H(x)}, namely  $(-1)^{m-s+1}h^P_{s-1}(\mathbf{u})\,e^P_{m-s+1}(\mathbf{u})$ and  $(-1)^{m-s}h^P_{s}(\mathbf{u})\,e^P_{m-s}(\mathbf{u})$; so its contribution to the left side of \eqref{e E(-x)H(x)} is 0.
Further, any monomial contributing to the left side of \eqref{e E(-x)H(x)} has this form.
Thus $h^P_1(\mathbf{u}), h^P_2(\mathbf{u}), \dots$ can be recursively expressed  in terms
of $e^P_1(\mathbf{u}), e^P_2(\mathbf{u}), \dots$. 
\end{proof}

\subsection{$P$-descents and the noncommutative $P$-Cauchy product}
\label{ss P descents}

Let $\U_P^*$ be the free $\mathbb{Z}$-module with basis the set of words
in the alphabet  $P$.
Let $\langle\cdot ,\cdot  \rangle \colon \U_P \times \U_P^* \to \ZZ$ denote the
natural pairing
in which the basis of noncommutative monomials is dual to the basis of words, i.e.,
if $\mathbf{u}_{\e{w}}=u_{\e{w}_1}\cdots u_{\e{w}_n}\in\U_P$
is a monomial corresponding to the word
$\e{w} = \e{w}_1 \cdots \e{w}_n\in\U_P^*$,
and $\e{v}\in\U_P^*$ is another word, then
$\langle \mathbf{u}_{\e{w}},\e{v}\rangle =\delta_{\e{vw}}$. 

For an ideal  $I\subset\U_P$, we denote by $I^\perp$ the orthogonal complement
\[
I^\perp = \big\{ \gamma \in \U_P^* : \langle z, \gamma \rangle = 0 \text{ for all }z
\in I \big\}.
\]
Note that since $f\equiv g \bmod I$ implies
$\langle f,\gamma \rangle=\langle g,\gamma \rangle$
for all $\gamma\in I^\perp$,
any element of $\U_P/I$ has a well-defined pairing with any
element of~$I^\perp$.

The \emph{content} of a word  $\e{w} \in \U_P^*$ is the multiset of letters occurring in  $\e{w}$.
For example, the content of  $\e{w} = \e{431411}$ is the multiset $\{1,1,1,3,4,4\}$.
For any multiset  $\beta$ of elements of  $P$, define $\Wo{\beta} \in \U_P^*$ by
\begin{align}
\label{eq: Wo beta}
\Wo{\beta} = \text{ the sum of all words of content  $\beta$}.
\end{align}
In particular,  $\Wo{P}$ is the sum of all  $N!$ permutations of  $P$.  The elements
$\Wo{\beta}$ form a  $\ZZ$-basis for $(\Ipol^P)^\perp \subset \Icommp{P}^\perp$ and are
the simplest and most important examples of elements in $(\Icommp{P})^\perp$.

As touched on in the introduction, to certain elements of $\U_P^*$, we associate 
elements of the ring of symmetric functions  $\Lambda(\mathbf{x}) \subset \ZZ[[x_1,x_2,\dots]]$ in
\emph{commuting} variables $\mathbf{x}=(x_1,x_2,\dots)$.

\begin{definition}
For a word $\e{w} = \e{w}_1 \cdots \e{w}_n \in \U_P^*$,
its \emph{$P$-descent set} is
\begin{align}
\Des_P(\e{w})= \{i \in \{1,\dots,n-1\} : \e{w}_i \leftto_P \e{w}_{i+1} \}.
\end{align}
For any $\gamma=\displaystyle\sum_{\e{w}} \gamma_{\e{w}} \,\e{w}\in\U_P^*$,
define the formal
power series $F_\gamma(\mathbf{x})\in \mathbb{Z}[[x_1,x_2,\dots]]$ by
\begin{equation}
\label{eq:Fgamma-via-Q}
F_\gamma(\mathbf{x})=\displaystyle\sum_{\e{w}} \gamma_{\e{w}}\,
Q_{\Des_P(\e{w})}(\mathbf{x}),
\end{equation}
where
\[
Q_{\Des_P(\e{w})}(\mathbf{x}) = \sum_{\substack{1 \le i_1 \le \,
    \cdots \, \le i_n\\j \in \Des_P(\e{w}) \implies i_j < i_{j+1} }}
\!\! x_{i_1}\cdots x_{i_n}
\]
is the \emph{fundamental quasisymmetric function} \cite{GesselPPartition} associated to  $\Des_P(\e{w})$.
\end{definition}

The $F_\gamma(\mathbf{x})$ can also be defined  in
terms of a variant of the Cauchy product.

\begin{definition}
\label{def:cauchy-product}
Define the \emph{ noncommutative $P$-Cauchy product} $\Omega(\mathbf{x},\mathbf{u})\in\U_P[[x_1,x_2,\dots]]$~by
\begin{equation}
\label{eq:cauchy-product}
\Omega(\mathbf{x}, \mathbf{u})
=\Big( \sum_{\ell = 0}^\infty x_1^\ell h^P_\ell(\mathbf{u}) \Big)\Big( \sum_{\ell = 0}^\infty x_2^\ell h^P_\ell(\mathbf{u}) \Big) \Big( \sum_{\ell = 0}^\infty x_3^\ell h^P_\ell(\mathbf{u}) \Big)\cdots.
\end{equation}
\end{definition}

Collecting the terms in~\eqref{eq:cauchy-product} involving each noncommutative monomial~$\mathbf{u}_{\e{w}}$, we have $\Omega(\mathbf{x}, \mathbf{u}) = \sum_{\e{w}}
Q_{\Des_P(\e{w})}(\mathbf{x}) \mathbf{u}_{\e{w}}$, and hence
\begin{equation}
\label{eq:F-gamma}
F_\gamma(\mathbf{x})=\big\langle \Omega(\mathbf{x}, \mathbf{u}) ,\gamma\big\rangle,
\end{equation}
where  $\langle\cdot ,\cdot  \rangle$ is
the bilinear pairing between $\U_P[[x_1,x_2,\dots]]$ and~$\U_P^*$ in which $\langle \sum_{\alpha} f_\alpha(\mathbf{u}) \mathbf{x}^\alpha , \e{v}\rangle = \sum_\alpha \mathbf{x}^\alpha \langle f_\alpha(\mathbf{u}) , \, \e{v} \rangle$ for
any  $\sum_\alpha f_\alpha(\mathbf{u}) \mathbf{x}^\alpha  \in \U_P[[x_1,x_2,\dots]]$.

\begin{proposition}
\label{Icomm-symmetric}
If $\gamma\in (\Icommp{P})^\perp$, then $F_\gamma(\mathbf{x})$
is symmetric in $x_1, x_2,\dots$.
\end{proposition}
\begin{proof}
Follows from \eqref{eq:F-gamma} and Proposition \ref{cor:h's-commute}, by
the same proof as in \cite[Prop. 2.6]{BF}.
\end{proof}

For  $P$ a total order, the  $F_\gamma$ include many important classes of symmetric functions such as LLT polynomials and stable Grothendieck polynomials (see \cite{BF}). For this paper,
the prime examples of the $F_\gamma$ are the chromatic symmetric functions:

\begin{theorem}[{\cite[Corollary 2]{ChowDescents}}]
\label{th:chromatic as F}
For any poset  $P$, the symmetric function $F_{\Wo{P}}(\mathbf{x})$ associated to the sum  $\Wo{P}$ of permutations of  $P$ is equal to $\omega$ of the chromatic symmetric function of the incomparability graph of  $P$:
\begin{align}
\label{eq:chromatic as F 1}
F_{\Wo{P}}(\mathbf{x}) = \omega \, X_{\inc(P)}(\mathbf{x}),
\end{align}
where $\omega \colon \Lambda(\mathbf{x}) \rightarrow \Lambda(\mathbf{x}) $
is the algebra involution
determined by
$\omega \spa h_\lambda(\mathbf{x}) = e_\lambda(\mathbf{x})$.
\end{theorem}

Further, we will see in \S\ref{s P Knuth and H graphs}
that several refinements of chromatic symmetric functions, including the
$t$-chromatic symmetric functions, are examples of $F_\gamma$.

\subsection{Noncommutative $P$-Schur and  $P$-monomial functions and positivity}
\label{ss P Schur etc}

Extending the machinery of \cite{FG, BF},
we formulate the first step in an approach to (refinements of) the Stanley-Stembridge conjecture
by providing a
tool to study the expansion of the symmetric functions  $F_\gamma(\mathbf{x})$
in various bases.

\begin{definition}
For the three bases  $e_\lambda$,  $s_\lambda$,   $m_\lambda$ of symmetric functions, we  define noncommutative $P$-versions $\mathfrak{e}^P_\lambda(\mathbf{u})$,  $\mathfrak{J}^P_\lambda(\mathbf{u})$,   $\mathfrak{m}^P_\lambda(\mathbf{u})$ as certain elements of  $\U_P$.
Let $\lambda = (\lambda_1, \dots, \lambda_\ell)$ be a partition. The first is defined simply as
\begin{align}
\mathfrak{e}^P_\lambda(\mathbf{u}) = e_{\lambda_1}(\mathbf{u}) \cdots e_{\lambda_\ell}(\mathbf{u}).
\end{align}
Next, define the \emph{noncommutative $P$-Schur function} by the following $P$-version of the classical Kostka-Naegelsbach determinantal formula:
\begin{equation}
\label{eq:P J def}
\mathfrak{J}^P_\lambda(\mathbf{u}) = \sum_{\pi\in \SS_{k}}
\sgn(\pi) \, e^P_{\lambda'_1+\pi(1)-1}(\mathbf{u})\,
e^P_{\lambda'_2+\pi(2)-2}(\mathbf{u}) \cdots
e^P_{\lambda'_{k}+\pi(k)-k}(\mathbf{u}),
\end{equation}
where $\lambda'$ is the transpose of the partition $\lambda$ and  $k = \lambda_1$.
Define the \emph{noncommutative $P$-monomial function} by
\begin{equation}
\label{eq:P m def}
\mathfrak{m}^P_\lambda(\mathbf{u}) = \sum_{\mu} (K^{-1})_{\lambda \mu}  \,  \mathfrak{J}^P_\mu(\mathbf{u}),
\end{equation}
where  $K^{-1}$ is the inverse Kostka matrix.
\end{definition}

Let $\Lambda(\mathbf{y})$ denote the ordinary ring of symmetric polynomials in
an infinite alphabet of commuting variables $\mathbf{y}=(y_1,y_2,\dots)$.
This is a polynomial ring with algebraically independent generators
$e_1(\mathbf{y}),e_2(\mathbf{y}), \dots$.
Define the ring homomorphism
\begin{align}
\label{eq:psi def}
\psi \colon \Lambda(\mathbf{y})\to \U_P/\Icommp{P}  \quad \text{ by } \ \
 e_k(\mathbf{y})\mapsto e^P_k(\mathbf{u}).
\end{align}
The elements $e_\lambda(\mathbf{y})$,  $s_\lambda(\mathbf{y})$,   $m_\lambda(\mathbf{y})$ of $\Lambda(\mathbf{y})$ are mapped by  $\psi$ to the images of  $\mathfrak{e}^P_\lambda(\mathbf{u})$,  $\mathfrak{J}^P_\lambda(\mathbf{u})$, and  $\mathfrak{m}^P_\lambda(\mathbf{u})$ in  $\U_P/\Icommp{P}$, respectively.

\begin{remark}
We could get away with defining $\mathfrak{e}^P_\lambda(\mathbf{u}),$ $\mathfrak{J}^P_\lambda(\mathbf{u})$, $\mathfrak{m}^P_\lambda(\mathbf{u})$ as the elements  $\psi(e_\lambda),$ $\psi(s_\lambda),$
$\psi(m_\lambda)$ of
$\U_P/\Icommp{P}$, 
but it is preferable for proofs later on to realize these as
actual elements of  $\U_P$.
\end{remark}

\begin{proposition}[Noncommutative $P$-Cauchy formulas {\cite{Hwang}}]
\label{pr:P Cauchy}
The following identities hold modulo $\Icommp{P}[[\mathbf{x}]]$
\begin{equation}
\label{eq:P Cauchy}
\Omega(\mathbf{x}, \mathbf{u})
\, \equiv \,
\sum_{\lambda} s_\lambda(\mathbf x)
\mathfrak J^P_{\lambda}(\mathbf{u})
\, \equiv \,
\sum_{\lambda} h_\lambda(\mathbf x)
\mathfrak m^P_{\lambda}(\mathbf{u})
\, \equiv \,
\sum_{\lambda} \omega \spa m_\lambda(\mathbf x)
e^P_{\lambda}(\mathbf{u}).
\end{equation}
\end{proposition}
\begin{proof}
Apply $\psi$ to the classical Cauchy identity  $\sum_{\lambda} m_\lambda(\mathbf x)h_\lambda(\mathbf{y})
= \sum_\lambda s_\lambda(\mathbf{x}) s_\lambda(\mathbf{y})$
to obtain the first congruence.
The others follow similarly from well-know variations of the Cauchy identity (see, e.g. \cite[I Ch. 4]{Macdonald95}).
\end{proof}

\begin{corollary}[\cite{Hwang}]
\label{Icomm-Schur}
For any $\gamma\in (\Icommp{P})^\perp$, we have
\begin{alignat}{1}
\label{eq:F-via-schurs}
F_\gamma(\mathbf{x}) 
& = \sum_{\lambda} s_{\lambda}(\mathbf x) \big\langle
\mathfrak J^P_{\lambda}(\mathbf{u}),  \gamma \big\rangle,  \\
\label{eq:F-via-m}
\omega\spa F_\gamma(\mathbf{x}) & = \sum_{\lambda} e_\lambda(\mathbf x) \big\langle
\mathfrak m^P_{\lambda}(\mathbf{u}),  \gamma \big\rangle, \\
\label{eq:F-via-e}
\omega\spa F_\gamma(\mathbf{x}) & = \sum_{\lambda} m_\lambda(\mathbf x) \big\langle
 e^P_{\lambda}(\mathbf{u}),  \gamma \big\rangle.
\end{alignat}
\end{corollary}
\begin{proof}
Apply  $\langle \cdot , \gamma\rangle$ to \eqref{eq:P Cauchy} and combine with
\eqref{eq:F-gamma}.
\end{proof}

This gives a route to showing that the functions $\omega \spa F_\gamma$
are Schur positive or  $e$-positive.
An element $f\in\U_P$ is \emph{$\mathbf{u}$-monomial positive modulo an ideal~$I$}
if $f$ can be written as an element of $I$ plus a nonnegative integer combination of
noncommutative monomials  $\mathbf{u}_{\e{w}}$.
Let  $(\U_P^*)_{\ge 0}$ denote the set of nonnegative integer combinations of words in  $\U_P^*$.

\begin{corollary}
\label{c positivity}
Let $I\subset\U$ be an ideal containing~$\Icommp{P}$.

(i) If $\mathfrak{J}^P_\lambda(\mathbf{u})$ is $\mathbf{u}$-monomial
positive modulo~$I$, then for any
$\gamma \in (\U_P^*)_{\ge 0} \cap I^\perp$, the coefficient of  $s_\lambda(\mathbf{x})$
in the Schur expansion of $F_\gamma(\mathbf{x})$ is
the nonnegative integer
$\big\langle \mathfrak J^P_{\lambda}(\mathbf{u}),  \gamma \big\rangle$.

(ii) If  $\mathfrak{m}^P_\lambda(\mathbf{u})$ is $\mathbf{u}$-monomial
positive modulo~$I$, then for any $\gamma \in (\U_P^*)_{\ge 0} \cap I^\perp$, the coefficient of  $e_\lambda(\mathbf{x})$ in the $e$-expansion of $\omega F_\gamma(\mathbf{x})$ is the nonnegative integer
 $\big\langle \mathfrak m^P_{\lambda}(\mathbf{u}),  \gamma \big\rangle$.
\end{corollary}

\section{$P$-Knuth equivalence and refinements of the Stanley-Stembridge conjecture}

We now turn to the deeper aspects of the theory---the selection of an ideal  $I$ with which to
apply Corollary \ref{c positivity} and establishing  $\mathbf{u}$-monomial positivity of
$\mathfrak{J}^P_\lambda(\mathbf{u})$ or  $\mathfrak{m}^P_\lambda(\mathbf{u})$
 mod  $I$.
Selecting $I$ is delicate:
choose it too small and the desired  $\mathbf{u}$-monomial positivity
will fail; but otherwise, the smaller the better
because (1) a smaller  $I$ means a larger  $I^\perp$, thus proving positivity for a larger class, and (2) the smaller  the $I$, the fewer the possible $\mathbf{u}$-monomial positive formulas for  $\mathfrak{J}^P_\lambda(\mathbf{u})$ or  $\mathfrak{m}^P_\lambda(\mathbf{u})$ mod  $I$, which often makes it easier to single out a nice canonical one
(for more on this,
see \S\ref{ss e pos} as well as Section 7 and the discussion after Remark 2.11 in \cite{BF}).

We focus on three ideals  which seem to be good settings to study chromatic symmetric functions.
For the first,  $\Iplacp{P}$, we prove  $\mathbf{u}$-monomial positivity of  $\mathfrak{J}_\lambda^P(\mathbf{u})$ mod $\Iplacp{P}$ to establish our main theorem (\S\S\ref{ss P Knuth}--\ref{ss Schur pos P Knuth}).
The second,  $I_T^P$, is the natural setting with which to reformulate results and conjectures
of Shareshian-Wachs in our language (\S\ref{ss t chromatic}).
The last,  $I_H^P$,
is suited to studying 
the Stanley-Stembridge conjecture  (\S\ref{ss e pos}).

\label{s P Knuth and H graphs}
\subsection{The  $P$-plactic algebra and  $P$-Knuth equivalence graphs}
\label{ss P Knuth}
Building off of work in \cite{KPchromatic}, we define a version of the plactic algebra and Knuth equivalence graphs
for any \threeone-free poset  $P$.  (A poset is \threeone-free if it contains no induced subposet isomorphic to
the disjoint union of a 3 element chain and a 1 element chain.)
Recall our notation $a \to_P b$,  $a \dote_P b$,  $a \dotto_P b$ for poset relations from \eqref{e atob def1}--\eqref{e atob def3}.

\begin{theorem} 
\label{th:IplacP elem commute}
Let $P$ be a (finite) \threeone-free  poset. Let $\Iplacp{P}$ be the ideal of~$\U_P$ generated by the following elements:
\begin{alignat}{3}
&u_b u_a u_c - u_b u_c u_a  &  \qquad(a \to_P b \dotto_P c \ \text{ and } \ a \to_P c), \label{poset rel knuth etc bca}\\
&u_c u_a u_b - u_a u_c u_b &  \qquad(a \dotto_P b \to_P c \ \text{ and } \ a \to_P c),
\label{poset rel knuth etc cab}\\
& u_c u_a u_b -u_b u_c u_a  &  \qquad(a \dote_P b \dote_P c \ \text{ and } \  a \to_P c).\, \, \spa \label{poset rel bca cab}
\end{alignat}
The noncommutative elementary functions $e_k^P(\mathbf{u}_S)$ and
$e_\ell^P(\mathbf{u}_S)$ commute with each other modulo~$\Iplacp{P}$, for any $k,\ell$ and
any~$S \subset P$
\begin{equation}
\label{eq:e's-commute}
e_k^P(\mathbf{u}_S)\,e_\ell^P(\mathbf{u}_S)\equiv e_\ell^P(\mathbf{u}_S)\,e_k^P(\mathbf{u}_S) \ \bmod \, \Iplacp{P}.
\end{equation}
\end{theorem}

The proof will be given in Section \ref{sec:elem commute}. Note that this theorem is equivalent to the statement that $\Iplacp{P} \supset \Icommp{P}$.
When  $P$ is a total order, $\U_P/\Iplacp{P}$ is the plactic algebra of Lascoux and Sch\"utzenberger~\cite{LS}.


The elements of  $(\U_P^*)_{\ge 0} \cap (\Iplacp{P})^\perp$ have a natural description in terms of  graphs.

\begin{definition}[\emph{$P$-Knuth equivalence graphs}]
\label{def-switchboard}
Let  $P$ be a \threeone-free poset.
Let $\e{w}=\e{w}_1\cdots \e{w}_n$ and $\e{w'}=\e{w}_1'\cdots
  \e{w}_n'$
be two words of the same length~$n$ in the
alphabet $P$.
We say that $\e{w}$ and $\e{w'}$  are related by a
\emph{$P$-Knuth transformation} in position~$i$ if
$\e{w}_j=\e{w}_j'$ for $j\notin\{i-1,i,i+1\}$,
and the pair of subwords
$\{\e{w_{i-1}w_iw_{i+1}},
\e{w}'_{i-1}\e{w}'_i\e{w}'_{i+1}\}$
has one of the following forms:
\pagebreak[2]
\begin{itemize}
\item[(1)]
$\{\e{bac}, \e{bca}\}$, with $a \to_P b \dotto_P c \ \text{ and } \ a \to_P c$;
\item[(2)]
$\{\e{cab}, \e{acb}\}$, with $a \dotto_P b \to_P c \ \text{ and } \ a \to_P c$;
\item[(3)]
$\{\e{cab}, \e{bca}\}$, with $a \dote_P b \dote_P c \ \ \text{ and } \  a \to_P c$.
\end{itemize}
Let  $\Gamma^n$ denote the edge-labeled graph on the
set of words of length~$n$ whose  $i$-edges are all possible  $P$-Knuth transformations in position  $i$
on words of length  $n$.
A \emph{$P$-Knuth equivalence graph} is a
disjoint union of connected components of  $\Gamma^n$.

A  $P$-Knuth equivalence graph  is \emph{standard} if its vertices are permutations of  $P$.
\end{definition}

See Figures \ref{f P switchboard 1234}, \ref{f 2p2}, \ref{f P switchboard 11234}, and \ref{fig:H disconnected} for examples.  These figures and many of our other examples use the following class of \threeone-free posets:
\begin{definition}
\label{def:jump poset}
For a positive integer $k$, let $\PP_{k}$ be the poset on $\{1,2,\dots, N\}$ with relations $a \to_P c$ whenever $c-a \ge k$.
When these posets are used in examples, we take $N$ to be any number larger than those appearing in the example.
\end{definition}

\begin{remark}
The $P$-Knuth transformations can be illustrated more explicitly  by organizing them
according to the induced subposet on the two or three letters involved:
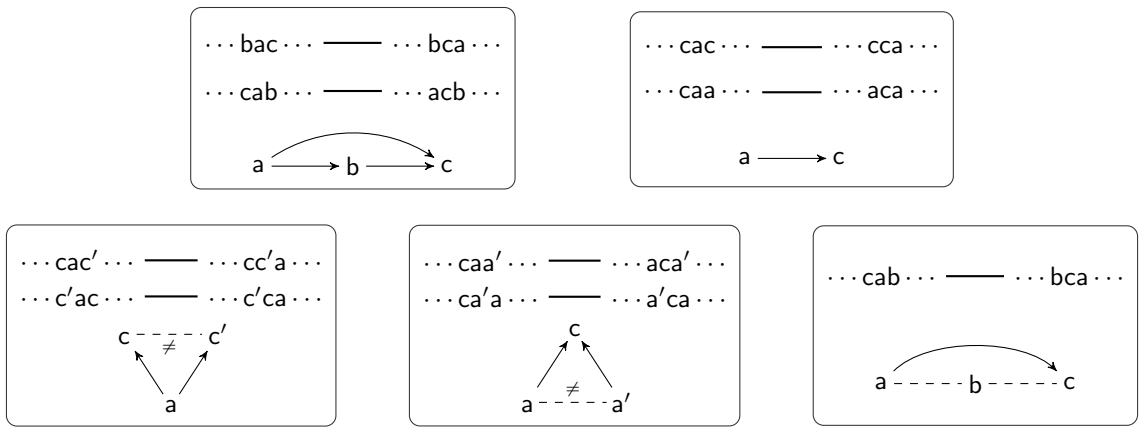
\begin{center}
\begin{tikzpicture}[xscale = 2.5,yscale = 2, framed, background rectangle/.style={draw=black!80, rounded corners}]
\tikzstyle{vertex}=[inner sep=0pt, outer sep=4pt]
\tikzstyle{framedvertex}=[inner sep=3pt, outer sep=4pt, draw=gray]
\tikzstyle{edge} = [draw, thick, -,black]
\tikzstyle{LabelStyleH} = [text=black, anchor=south]
\tikzstyle{LabelStyleHn} = [text=black, anchor=north]
\tikzstyle{LabelStyleV} = [text=black, anchor=east]
\tikzstyle{LabelStyleVw} = [text=black, anchor=west]
\node[vertex] (v1) at (0,0){\footnotesize $\e{\cdots cab\cdots} $};
\node[vertex] (v2) at (1,0){\footnotesize $\e{\cdots acb\cdots} $};
\node[vertex] (v3) at (0,0.32){\footnotesize $\e{\cdots bac\cdots} $};
\node[vertex] (v4) at (1,0.32){\footnotesize $\e{\cdots bca\cdots} $};
\draw[edge] (v1) to node[LabelStyleH]{} (v2);
\draw[edge] (v3) to node[LabelStyleH]{} (v4);
\tikzstyle{vertex}=[inner sep=0pt, outer sep= 2.5pt]
\tikzstyle{aedge} = [draw, ->,>=stealth', black]
\tikzstyle{aedgecurve} = [draw, ->,>=stealth', black, bend left=40]
\tikzstyle{edge} = [draw, thick, -,black]
\tikzstyle{dashededge} = [draw, -,  dashed, black]
\node[vertex] (va) at (0,-0.5){\footnotesize$\e{a}$};
\node[vertex] (vb) at (0.5,-0.5){\footnotesize$\e{b}$};
\node[vertex] (vc) at (1,-0.5){\footnotesize$\e{c}$};
\draw[aedge] (va) to (vb);
\draw[aedge] (vb) to (vc);
\draw[aedgecurve] (va) to (vc);
\end{tikzpicture}
\quad \quad \quad
\raisebox{.3mm}{
\begin{tikzpicture}[xscale = 2.5,yscale = 2, framed, background rectangle/.style={draw=black!80, rounded corners}]
\tikzstyle{vertex}=[inner sep=0pt, outer sep=4pt]
\tikzstyle{framedvertex}=[inner sep=3pt, outer sep=4pt, draw=gray]
\tikzstyle{aedge} = [draw, thin, ->,black]
\tikzstyle{edge} = [draw, thick, -,black]
\tikzstyle{doubleedge} = [draw, thick, double distance=1pt, -,black]
\tikzstyle{hiddenedge} = [draw=none, thick, double distance=1pt, -,black]
\tikzstyle{dashededge} = [draw, very thick, dashed, black]
\tikzstyle{LabelStyleH} = [text=black, anchor=south]
\tikzstyle{LabelStyleHn} = [text=black, anchor=north]
\tikzstyle{LabelStyleV} = [text=black, anchor=east]
\tikzstyle{LabelStyleVw} = [text=black, anchor=west]
\node[vertex] (v1) at (0,0){\footnotesize $\e{\cdots caa\cdots} $};
\node[vertex] (v2) at (1,0){\footnotesize $\e{\cdots aca\cdots} $};
\node[vertex] (v3) at (0,0.3){\footnotesize $\e{\cdots cac\cdots} $};
\node[vertex] (v4) at (1,0.3){\footnotesize $\e{\cdots cca\cdots} $};
\draw[edge] (v1) to node[LabelStyleH]{} (v2);
\draw[edge] (v3) to node[LabelStyleH]{} (v4);
\tikzstyle{vertex}=[inner sep=0pt, outer sep= 2.5pt]
\tikzstyle{aedge} = [draw, ->,>=stealth', black]
\tikzstyle{aedgecurve} = [draw, ->,>=stealth', black, bend left=40]
\tikzstyle{edge} = [draw, thick, -,black]
\tikzstyle{dashededge} = [draw, -,  dashed, black]
\node[vertex] (va) at (0.25,-0.44){\footnotesize$\e{a}$};
\node[vertex] (vb) at (0.75,-0.44){\footnotesize$\e{c}$};
\node[vertex] (g) at (0.75,-0.5){\footnotesize$\phantom{c} $};
\node[vertex] (gg) at (0.75,.33){\footnotesize$\phantom{a} $};
\draw[aedge] (va) to (vb);
\end{tikzpicture}}
\end{center}
\vspace{1mm}
\begin{center}
\begin{tikzpicture}[xscale = 2.5,yscale = 1.2, framed, background rectangle/.style={draw=black!80, rounded corners}]
\tikzstyle{vertex}=[inner sep=0pt, outer sep=4pt]
\tikzstyle{framedvertex}=[inner sep=3pt, outer sep=4pt, draw=gray]
\tikzstyle{aedge} = [draw, thin, ->,black]
\tikzstyle{edge} = [draw, thick, -,black]
\tikzstyle{doubleedge} = [draw, thick, double distance=1pt, -,black]
\tikzstyle{hiddenedge} = [draw=none, thick, double distance=1pt, -,black]
\tikzstyle{dashededge} = [draw, very thick, dashed, black]
\tikzstyle{LabelStyleH} = [text=black, anchor=south]
\tikzstyle{LabelStyleHn} = [text=black, anchor=north]
\tikzstyle{LabelStyleV} = [text=black, anchor=east]
\tikzstyle{LabelStyleVw} = [text=black, anchor=west]
\node[vertex] (v1) at (0,0.1){\footnotesize $\e{\cdots c'ac\cdots} $};
\node[vertex] (v2) at (1,0.1){\footnotesize $\e{\cdots c'ca\cdots} $};
\node[vertex] (v3) at (0,0.5){\footnotesize $\e{\cdots cac'\cdots} $};
\node[vertex] (v4) at (1,0.5){\footnotesize $\e{\cdots cc'a\cdots} $};
\draw[edge] (v1) to node[LabelStyleH]{} (v2);
\draw[edge] (v3) to node[LabelStyleH]{} (v4);
%
\tikzstyle{vertex}=[inner sep=0pt, outer sep= 2.5pt]
\tikzstyle{aedge} = [draw, ->,>=stealth', black]
\tikzstyle{aedgecurve} = [draw, ->,>=stealth', black, bend left=40]
\tikzstyle{edge} = [draw, thick, -,black]
\tikzstyle{dashededge} = [draw, -,  dashed, black]
\node[vertex] (va) at (0.25,-0.32){\footnotesize$\e{c\vphantom{'}}$};
\node[vertex] (vb) at (0.75,-0.32){\footnotesize$\e{c'}$};
\node[vertex] (vc) at (0.5,-1.121){\footnotesize$\e{a}$};
\draw[dashededge] (va) -- (vb) node [midway, below, inner sep = 2pt] {\Tiny$\ne$};
\draw[aedge] (vc) to (vb);
\draw[aedge] (vc) to (va);
\end{tikzpicture}
\quad \quad
\begin{tikzpicture}[xscale = 2.5,yscale = 1.2, framed, background rectangle/.style={draw=black!80, rounded corners}]
\tikzstyle{vertex}=[inner sep=0pt, outer sep=4pt]
\tikzstyle{framedvertex}=[inner sep=3pt, outer sep=4pt, draw=gray]
\tikzstyle{aedge} = [draw, thin, ->,black]
\tikzstyle{edge} = [draw, thick, -,black]
\tikzstyle{doubleedge} = [draw, thick, double distance=1pt, -,black]
\tikzstyle{hiddenedge} = [draw=none, thick, double distance=1pt, -,black]
\tikzstyle{dashededge} = [draw, very thick, dashed, black]
\tikzstyle{LabelStyleH} = [text=black, anchor=south]
\tikzstyle{LabelStyleHn} = [text=black, anchor=north]
\tikzstyle{LabelStyleV} = [text=black, anchor=east]
\tikzstyle{LabelStyleVw} = [text=black, anchor=west]

\node[vertex] (v1) at (0,.1){\footnotesize $\e{\cdots ca'a\cdots} $};
\node[vertex] (v2) at (1,.1){\footnotesize $\e{\cdots a'ca\cdots} $};
\node[vertex] (v3) at (0,0.5){\footnotesize $\e{\cdots caa'\cdots} $};
\node[vertex] (v4) at (1,0.5){\footnotesize $\e{\cdots aca'\cdots} $};
\draw[edge] (v1) to node[LabelStyleH]{} (v2);
\draw[edge] (v3) to node[LabelStyleH]{} (v4);

%
\tikzstyle{vertex}=[inner sep=0pt, outer sep= 1.9pt]
\tikzstyle{aedge} = [draw, ->,>=stealth', black]
\tikzstyle{aedgecurve} = [draw, ->,>=stealth', black, bend left=40]
\tikzstyle{edge} = [draw, thick, -,black]
\tikzstyle{dashededge} = [draw, -,  dashed, black]
\node[vertex] (va) at (0.25,-1.071){\footnotesize$\e{a\vphantom{'}}$};
\node[vertex] (vb) at (0.75,-1.071){\footnotesize$\e{a'}$};
\node[vertex] (vc) at (0.5,-0.27){\footnotesize$\e{c}$};
\draw[dashededge] (va) -- (vb) node [midway, above, inner sep = 2pt] {\Tiny$\ne$};
\draw[aedge] (vb) to (vc);
\draw[aedge] (va) to (vc);
\end{tikzpicture}
\quad \ \
\raisebox{.1mm}{
\begin{tikzpicture}[xscale = 2.5,yscale = 1.57, framed, background rectangle/.style={draw=black!80, rounded corners}]
\tikzstyle{vertex}=[inner sep=0pt, outer sep=4pt]
\tikzstyle{framedvertex}=[inner sep=3pt, outer sep=4pt, draw=gray]
\tikzstyle{edge} = [draw, thick, -,black]
\tikzstyle{LabelStyleH} = [text=black, anchor=south]
\tikzstyle{LabelStyleHn} = [text=black, anchor=north]
\tikzstyle{LabelStyleV} = [text=black, anchor=east]
\tikzstyle{LabelStyleVw} = [text=black, anchor=west]
\node[vertex] (v00) at (0,0.3){\footnotesize $ \ \ $};
\node[vertex] (v1) at (0,0.0){\footnotesize $\e{\cdots cab\cdots} $};
\node[vertex] (v2) at (1,0.0){\footnotesize $\e{\cdots bca\cdots} $};
\draw[edge] (v1) to node[LabelStyleH]{} (v2);
\tikzstyle{vertex}=[inner sep=0pt, outer sep= 2.5pt]
\tikzstyle{aedge} = [draw, ->,>=stealth', black]
\tikzstyle{aedgecurve} = [draw, ->,>=stealth', black, bend left=58]
\tikzstyle{edge} = [draw, thick, -,black]
\tikzstyle{dashededge} = [draw, -,  dashed, black]
\node[vertex] (va) at (0,-0.91){\footnotesize$\e{a}$};
\node[vertex] (vb) at (0.5,-0.91){\footnotesize$\e{b}$};
\node[vertex] (vc) at (1,-0.91){\footnotesize$\e{c}$};
\draw[dashededge] (va) to (vb);
\draw[dashededge] (vb) to (vc);
\draw[aedgecurve] (va) to (vc);
\node[vertex] (g) at (0.5,-1.1){\footnotesize$\phantom{a} $};
\node[vertex] (gg) at (0.5,0.14){\footnotesize$\phantom{a} $};
\end{tikzpicture}}
\end{center}

Here,   $\begin{tikzpicture}[xscale = 2.5,yscale = 1.4]
\tikzstyle{vertex}=[inner sep=0pt, outer sep= 2.5pt]
\tikzstyle{aedge} = [draw, ->,>=stealth', black]
\tikzstyle{aedgecurve} = [draw, ->,>=stealth', black, bend left=40]
\tikzstyle{edge} = [draw, thick, -,black]
\tikzstyle{dashededge} = [draw, -,  dashed, black]
\node[vertex] (va) at (0.25,-1.071){\footnotesize$\e{a\vphantom{'}}$};
\node[vertex] (vb) at (0.75,-1.071){\footnotesize$\e{a'}$};
\draw[dashededge] (va) -- (vb) node [midway, above, inner sep = 2pt] {\Tiny$\ne$};
\end{tikzpicture}$ indicates that  $\e{a}$ and  $\e{a'}$ are incomparable and unequal elements of  $P$.

Note that 
a word $\e{w}$ is related to exactly one other word  $\e{w'}$  by a
$P$-Knuth transformation in position~$i$ if $|\Des_P(\e{w}) \cap \{i-1,i\}| = 1$, and is related to no other words by such a transformation if $|\Des_P(\e{w}) \cap \{i-1,i\}| \ne 1$.
\end{remark}

\begin{remark}
(i) Our definition of  $P$-Knuth equivalence graphs extends that in \cite[Definition 4.6]{KPchromatic}---from
natural unit interval orders to \threeone-free posets and from standard $P$-Knuth equivalence graphs to
ones on words of any content.

(ii) When  $P$ is a total order,
$P$-Knuth equivalence graphs are the same as the \emph{plactic switchboards} in \cite[\S5]{BF}, and, after forgetting vertex labels and with slightly different conventions,
they are the \emph{standard dual equivalence graphs} of~\cite{SamiOct13}.
See \cite[\S5]{BF} for further discussion.
\end{remark}

The  $0,1$-vectors in  $(\Iplacp{P})^\perp$ which are supported on words of a fixed length are exactly the vectors $\sum_{\e{w} \in  \ver(\Gamma)} \e{w}$ for  $P$-Knuth equivalence graphs  $\Gamma$ (by a $0,1$-vector in  $\U_P^*$ we mean a linear combination of words with coefficients 0 and 1).

\begin{definition}
\label{def:FGamma}
For a  $P$-Knuth equivalence graph~$\Gamma$, let $\gamma=\sum_{\e{w} \in  \ver(\Gamma)} \e{w}$, and define
the \emph{symmetric function of  $\Gamma$} by
\begin{align}
F_\Gamma(\mathbf{x})=F_\gamma(\mathbf{x}) =
\sum_{\e{w} \in \ver(\Gamma)}
Q_{\Des_P(\e{w})}(\mathbf{x}).
\end{align}
\end{definition}

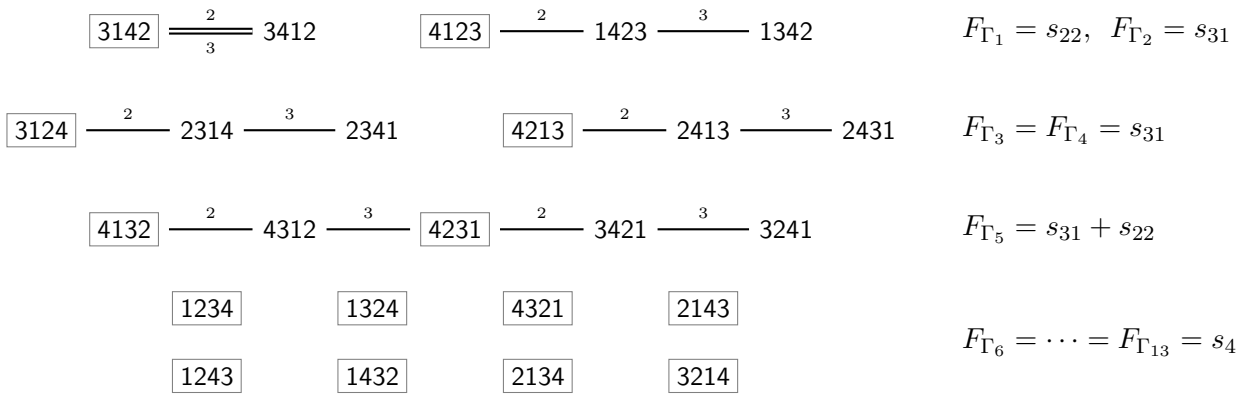
\begin{figure}
\centerfloat
\begin{tikzpicture}[xscale = 2.2,yscale = 1.5]
\tikzstyle{vertex}=[inner sep=0pt, outer sep=4pt]
\tikzstyle{framedvertex}=[inner sep=3pt, outer sep=4pt, draw=gray]
\tikzstyle{aedge} = [draw, thin, ->,black]
\tikzstyle{edge} = [draw, thick, -,black]
\tikzstyle{doubleedge} = [draw, thick, double distance=1pt, -,black]
\tikzstyle{hiddenedge} = [draw=none, thick, double distance=1pt, -,black]
\tikzstyle{dashededge} = [draw, very thick, dashed, black]
\tikzstyle{LabelStyleH} = [text=black, anchor=south]
\tikzstyle{LabelStyleHn} = [text=black, anchor=north]
\tikzstyle{LabelStyleV} = [text=black, anchor=east]

\begin{scope}[yshift=0]
\node[framedvertex] (v1) at (-2,0){\footnotesize$\e{3142} $};
\node[vertex] (v2) at (-1,0){\footnotesize$\e{3412} $};
\draw[doubleedge] (v1) to node[LabelStyleH]{\Tiny$2 $} (v2);
\draw[hiddenedge] (v1) to node[LabelStyleHn]{\Tiny$3 $} (v2);
\node[vertex] (w1) at (2,0){\footnotesize$\e{1342} $};
\node[vertex] (w2) at (1,0){\footnotesize$\e{1423} $};
\node[framedvertex] (w3) at (0,0){\footnotesize$\e{4123} $};
\draw[edge] (w1) to node[LabelStyleH]{\Tiny$3 $} (w2);
\draw[edge] (w2) to node[LabelStyleH]{\Tiny$2 $} (w3);
\node[vertex,anchor = west] (vv) at (3,0){\small$F_{\Gamma_1} = s_{22}$, \ $F_{\Gamma_2} = s_{31}$};
\end{scope}

\begin{scope}[yshift=-1*25]
\node[vertex] (v1) at (-.5,0){\footnotesize$\e{2341} $};
\node[framedvertex] (v2) at (-2.5,0){\footnotesize$\e{3124} $};
\node[vertex] (v3) at (-1.5,0){\footnotesize$\e{2314} $};
\draw[edge] (v1) to node[LabelStyleH]{\Tiny$3 $} (v3);
\draw[edge] (v2) to node[LabelStyleH]{\Tiny$2 $} (v3);
\node[vertex] (w1) at (2.5,0){\footnotesize$\e{2431} $};
\node[vertex] (w2) at (1.5,0){\footnotesize$\e{2413} $};
\node[framedvertex] (w3) at (0.5,0){\footnotesize$\e{4213} $};
\draw[edge] (w1) to node[LabelStyleH]{\Tiny$3 $} (w2);
\draw[edge] (w2) to node[LabelStyleH]{\Tiny$2 $} (w3);
\node[vertex,anchor = west] (vv) at (3,0){\small$F_{\Gamma_3} = F_{\Gamma_4} = s_{31}$};
\end{scope}

\begin{scope}[yshift=-2*25]
\node[vertex] (v1) at (2,0){\footnotesize$\e{3241} $};
\node[vertex] (v2) at (1,0){\footnotesize$\e{3421} $};
\node[framedvertex] (v3) at (0,0){\footnotesize$\e{4231} $};
\node[vertex] (v4) at (-1,0){\footnotesize$\e{4312} $};
\node[framedvertex] (v5) at (-2,0){\footnotesize$\e{4132} $};
\draw[edge] (v1) to node[LabelStyleH]{\Tiny$3 $} (v2);
\draw[edge] (v2) to node[LabelStyleH]{\Tiny$2 $} (v3);
\draw[edge] (v3) to node[LabelStyleH]{\Tiny$3 $} (v4);
\draw[edge] (v4) to node[LabelStyleH]{\Tiny$2 $} (v5);
\node[vertex,anchor = west] (vv) at (3,0){\small$F_{\Gamma_5} = s_{31}+s_{22}$};
\end{scope}

\begin{scope}[yshift=-3*25-12, yscale = .6]
\node[framedvertex] (v1) at (-1.5,1){\footnotesize$\e{1234} $};
\node[framedvertex] (v2) at (-.5,1){\footnotesize$\e{1324} $};
\node[framedvertex] (v3) at (.5,1){\footnotesize$\e{4321} $};
\node[framedvertex] (v4) at (1.5,1){\footnotesize$\e{2143} $};

\node[framedvertex] (w1) at (-1.5,0){\footnotesize$\e{1243} $};
\node[framedvertex] (w2) at (-.5,0){\footnotesize$\e{1432} $};
\node[framedvertex] (w3) at (.5,0){\footnotesize$\e{2134} $};
\node[framedvertex] (w4) at (1.5,0){\footnotesize$\e{3214} $};
\node[vertex,anchor = west] (vv) at (3,0.5){\small$F_{\Gamma_6}= \cdots = F_{\Gamma_{13}} = s_{4}$};
\end{scope}

\end{tikzpicture}
\caption{\label{f P switchboard 1234}
For  $P = \PP_2$, the $P$-Knuth equivalence graph on the words  $\e{w}$ of content
$\{1,2,3,4\}$.
On the right are the symmetric functions  $F_{\Gamma_i}(\mathbf{x})$ of
each of its components $\Gamma_1, \dots, \Gamma_{13}$,
which by Theorem~\ref{c P KE schur pos} can be read off from the
\boxed{\text{vertices}}
which are column reading words of $P$-tableaux.
}
\end{figure}

\begin{figure}
\begin{center}
\begin{tikzpicture}[xscale = 2.2,yscale = 2.3]
\tikzstyle{vertex}=[inner sep=0pt, outer sep=4pt]
\tikzstyle{framedvertex}=[inner sep=3pt, outer sep=4pt, draw=gray]
\tikzstyle{aedge} = [draw, thin, ->,black]
\tikzstyle{edge} = [draw, thick, -,black]
\tikzstyle{doubleedge} = [draw, thick, double distance=1pt, -,black]
\tikzstyle{hiddenedge} = [draw=none, thick, double distance=1pt, -,black]
\tikzstyle{dashededge} = [draw, very thick, dashed, black]
\tikzstyle{LabelStyleH} = [text=black, anchor=south]
\tikzstyle{LabelStyleHn} = [text=black, anchor=north]
\tikzstyle{LabelStyleV} = [text=black, anchor=east]
\node[framedvertex] (v1) at (1,1){\footnotesize$\e{3142} $};
\node[vertex] (v2) at (2,1){\footnotesize$\e{4312} $};
\node[vertex] (v3) at (1,0.5){\footnotesize$\e{3421} $};
\node[framedvertex] (v4) at (2,0.5){\footnotesize$\e{4231} $};

\draw[edge] (v1) to node[LabelStyleH]{\Tiny$2 $} (v2);
\draw[edge] (v1) to node[LabelStyleV]{\Tiny$3 $} (v3);
\draw[edge] (v4) to node[LabelStyleV]{\Tiny$3 $} (v2);
\draw[edge] (v4) to node[LabelStyleH]{\Tiny$2
 $} (v3);

\node[vertex] (v) at (-0.75,0.75){\small$  \  $};
\node[vertex,anchor = west] (vv) at (4.4,0.75){\small$F_{\Gamma_1} = 2 \spa s_{22}$};
\end{tikzpicture}
\end{center}
\vspace{2.8mm}

\begin{center}
\begin{tikzpicture}[xscale = 2.2,yscale = 1.5]
\tikzstyle{vertex}=[inner sep=0pt, outer sep=4pt]
\tikzstyle{framedvertex}=[inner sep=3pt, outer sep=4pt, draw=gray]
\tikzstyle{aedge} = [draw, thin, ->,black]
\tikzstyle{edge} = [draw, thick, -,black]
\tikzstyle{doubleedge} = [draw, thick, double distance=1pt, -,black]
\tikzstyle{hiddenedge} = [draw=none, thick, double distance=1pt, -,black]
\tikzstyle{dashededge} = [draw, very thick, dashed, black]
\tikzstyle{LabelStyleH} = [text=black, anchor=south]
\tikzstyle{LabelStyleHn} = [text=black, anchor=north]
\tikzstyle{LabelStyleV} = [text=black, anchor=east]
\node[framedvertex] (v1) at (1,1){\footnotesize$\e{4213} $};
\node[vertex] (v2) at (2,1){\footnotesize$\e{1423} $};
\node[vertex] (v3) at (3,1){\footnotesize$\e{1342} $};

\node[framedvertex] (v4) at (4,1){\footnotesize$\e{3124} $};
\node[vertex] (v5) at (5,1){\footnotesize$\e{2314} $};
\node[vertex] (v6) at (6,1){\footnotesize$\e{2431} $};

\draw[edge] (v1) to node[LabelStyleH]{\Tiny$2 $} (v2);
\draw[edge] (v2) to node[LabelStyleH]{\Tiny$3 $} (v3);

\draw[edge] (v4) to node[LabelStyleH]{\Tiny$2 $} (v5);
\draw[edge] (v5) to node[LabelStyleH]{\Tiny$3 $} (v6);

\node[vertex,anchor = west] (vv) at (6.4,1){\small$F_{\Gamma_2} = F_{\Gamma_3} = s_{31}$};
\end{tikzpicture}
\end{center}
\vspace{2.8mm}

\begin{center}
\begin{tikzpicture}[xscale = 2.2,yscale = 1.5]
\tikzstyle{vertex}=[inner sep=0pt, outer sep=4pt]
\tikzstyle{framedvertex}=[inner sep=3pt, outer sep=4pt, draw=gray]
\tikzstyle{aedge} = [draw, thin, ->,black]
\tikzstyle{edge} = [draw, thick, -,black]
\tikzstyle{doubleedge} = [draw, thick, double distance=1pt, -,black]
\tikzstyle{hiddenedge} = [draw=none, thick, double distance=1pt, -,black]
\tikzstyle{dashededge} = [draw, very thick, dashed, black]
\tikzstyle{LabelStyleH} = [text=black, anchor=south]
\tikzstyle{LabelStyleHn} = [text=black, anchor=north]
\tikzstyle{LabelStyleV} = [text=black, anchor=east]
\node[framedvertex] (v1) at (1,1){\footnotesize$\e{1234} $};
\node[framedvertex] (v2) at (2,1){\footnotesize$\e{1243} $};
\node[framedvertex] (v3) at (3,1){\footnotesize$\e{1324} $};
\node[framedvertex] (v4) at (4,1){\footnotesize$\e{1432} $};
\node[framedvertex] (v5) at (5,1){\footnotesize$\e{2134} $};
\node[framedvertex] (v6) at (1,0.5){\footnotesize$\e{2143} $};
\node[framedvertex] (v7) at (2,0.5){\footnotesize$\e{2341} $};
\node[framedvertex] (v8) at (3,0.5){\footnotesize$\e{2413} $};
\node[framedvertex] (v9) at (4,0.5){\footnotesize$\e{3214} $};
\node[framedvertex] (v6) at (5,0.5){\footnotesize$\e{3241} $};
\node[framedvertex] (v6) at (1,0){\footnotesize$\e{3412} $};
\node[framedvertex] (v6) at (2,0){\footnotesize$\e{4123} $};
\node[framedvertex] (v6) at (3,0){\footnotesize$\e{4132} $};
\node[framedvertex] (v6) at (4,0){\footnotesize$\e{4321} $};

\node[vertex] (v) at (0.43,0.75){\small$  \  $};
\node[vertex,anchor = west] (vv) at (4.93,0){\small$F_{\Gamma_4}= F_{\Gamma_5}= \ldots = F_{\Gamma_{17}} = s_{5}$};
\end{tikzpicture}
\end{center}
\caption{\label{f 2p2}
For  $P$ the poset $\twotwo$ with incomparable chains $3\leftto_P 1$ and $4\leftto_P 2$, the $P$-Knuth equivalence graph on the set of permutations of  $P$.
The symmetric functions of its components can be read off from the
framed vertices as in Figure \ref{f P switchboard 1234}.
}
\end{figure}
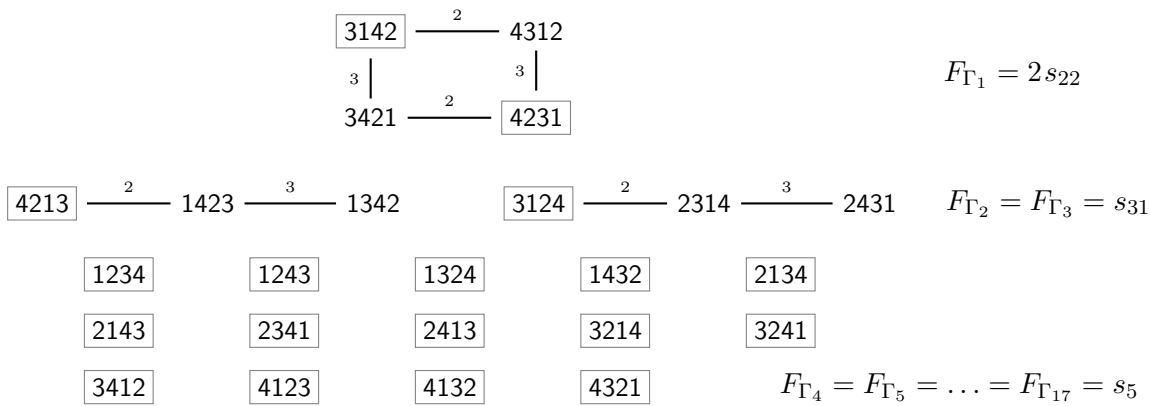

\begin{figure}
\begin{center}
\begin{tikzpicture}[xscale = 2.2,yscale = 2.3]
\tikzstyle{vertex}=[inner sep=0pt, outer sep=4pt]
\tikzstyle{framedvertex}=[inner sep=3pt, outer sep=4pt, draw=gray]
\tikzstyle{aedge} = [draw, thin, ->,black]
\tikzstyle{edge} = [draw, thick, -,black]
\tikzstyle{doubleedge} = [draw, thick, double distance=1pt, -,black]
\tikzstyle{hiddenedge} = [draw=none, thick, double distance=1pt, -,black]
\tikzstyle{dashededge} = [draw, very thick, dashed, black]
\tikzstyle{LabelStyleH} = [text=black, anchor=south]
\tikzstyle{LabelStyleHn} = [text=black, anchor=north]
\tikzstyle{LabelStyleV} = [text=black, anchor=east]
\node[framedvertex] (v1) at (1,1){\footnotesize$\e{41132} $};
\node[vertex] (v2) at (2,1){\footnotesize$\e{14132} $};
\node[vertex] (v3) at (4,1){\footnotesize$\e{14231} $};
\node[vertex] (v4) at (3,1){\footnotesize$\e{14312} $};
\node[framedvertex] (v5) at (6,1){\footnotesize$\e{41312} $};
\node[vertex] (v6) at (3,0.5){\footnotesize$\e{13241} $};
\node[vertex] (v7) at (7,1){\footnotesize$\e{43112} $};
\node[vertex] (v8) at (4,0.5){\footnotesize$\e{13421} $};
\node[vertex] (v9) at (5,1){\footnotesize$\e{41231} $};
\draw[edge] (v1) to node[LabelStyleH]{\Tiny$2 $} (v2);
\draw[edge] (v2) to node[LabelStyleH]{\Tiny$3 $} (v4);
\draw[edge] (v3) to node[LabelStyleH]{\Tiny$4 $} (v4);
\draw[edge] (v3) to node[LabelStyleV]{\Tiny$3 $} (v8);
\draw[edge] (v3) to node[LabelStyleH]{\Tiny$2 $} (v9);
\draw[doubleedge] (v5) to node[LabelStyleH]{\Tiny$2 $} (v7);
\draw[hiddenedge] (v5) to node[LabelStyleHn]{\Tiny$3 $} (v7);
\draw[edge] (v5) to node[LabelStyleH]{\Tiny$4 $} (v9);
\draw[edge] (v6) to node[LabelStyleH]{\Tiny$4 $} (v8);
\node[vertex,anchor = west] (vv) at (6.9,0.65){\small$F_{\Gamma_1} = s_{41}+s_{32}$};
\end{tikzpicture}
\end{center}
\vspace{2.8mm}

\begin{center}
\begin{tikzpicture}[xscale = 2.2,yscale = 1.5]
\tikzstyle{vertex}=[inner sep=0pt, outer sep=4pt]
\tikzstyle{framedvertex}=[inner sep=3pt, outer sep=4pt, draw=gray]
\tikzstyle{aedge} = [draw, thin, ->,black]
\tikzstyle{edge} = [draw, thick, -,black]
\tikzstyle{doubleedge} = [draw, thick, double distance=1pt, -,black]
\tikzstyle{hiddenedge} = [draw=none, thick, double distance=1pt, -,black]
\tikzstyle{dashededge} = [draw, very thick, dashed, black]
\tikzstyle{LabelStyleH} = [text=black, anchor=south]
\tikzstyle{LabelStyleHn} = [text=black, anchor=north]
\tikzstyle{LabelStyleV} = [text=black, anchor=east]
\node[vertex] (v1) at (3,1){\footnotesize$\e{23141} $};
\node[vertex] (v2) at (0,1){\footnotesize$\e{34121} $};
\node[framedvertex] (v3) at (1,1){\footnotesize$\e{31421} $};
\node[vertex] (v4) at (4,1){\footnotesize$\e{23411} $};
\node[vertex] (v5) at (2,1){\footnotesize$\e{31241} $};
\draw[doubleedge] (v1) to node[LabelStyleH]{\Tiny$3 $} (v4);
\draw[hiddenedge] (v1) to node[LabelStyleHn]{\Tiny$4 $} (v4);
\draw[edge] (v1) to node[LabelStyleH]{\Tiny$2 $} (v5);
\draw[doubleedge] (v2) to node[LabelStyleH]{\Tiny$2 $} (v3);
\draw[hiddenedge] (v2) to node[LabelStyleHn]{\Tiny$3 $} (v3);
\draw[edge] (v3) to node[LabelStyleH]{\Tiny$4 $} (v5);
\node[vertex] (v) at (-.9,0.75){\small$  \  $};
\node[vertex,anchor = west] (vv) at (4.9,1){\small$F_{\Gamma_2} = s_{32}$};
\end{tikzpicture}
\end{center}\vspace{1mm}

\begin{center}%
\begin{tikzpicture}[xscale = 2.2,yscale = 1.5]
\tikzstyle{vertex}=[inner sep=0pt, outer sep=4pt]
\tikzstyle{framedvertex}=[inner sep=3pt, outer sep=4pt, draw=gray]
\tikzstyle{aedge} = [draw, thin, ->,black]
\tikzstyle{edge} = [draw, thick, -,black]
\tikzstyle{doubleedge} = [draw, thick, double distance=1pt, -,black]
\tikzstyle{hiddenedge} = [draw=none, thick, double distance=1pt, -,black]
\tikzstyle{dashededge} = [draw, very thick, dashed, black]
\tikzstyle{LabelStyleH} = [text=black, anchor=south]
\tikzstyle{LabelStyleHn} = [text=black, anchor=north]
\tikzstyle{LabelStyleV} = [text=black, anchor=east]
\node[framedvertex] (v1) at (0,1){\footnotesize$\e{31214} $};
\node[vertex] (v2) at (1,1){\footnotesize$\e{23114} $};
\node[vertex] (v3) at (3,1){\footnotesize$\e{21341} $};
\node[vertex] (v4) at (2,1){\footnotesize$\e{21314} $};
\draw[edge] (v1) to node[LabelStyleH]{\Tiny$2 $} (v2);
\draw[edge] (v2) to node[LabelStyleH]{\Tiny$3 $} (v4);
\draw[edge] (v3) to node[LabelStyleH]{\Tiny$4 $} (v4);
\node[vertex] (v) at (-0.9,0.75){\small$  \  $};
\node[vertex,anchor = west] (vv) at (4.9,1){\small$F_{\Gamma_3} = s_{41}$};
\end{tikzpicture}
\end{center}\vspace{1.8mm}

\begin{center}%
\begin{tikzpicture}[xscale = 2.2,yscale = 1.5]
\tikzstyle{vertex}=[inner sep=0pt, outer sep=4pt]
\tikzstyle{framedvertex}=[inner sep=3pt, outer sep=4pt, draw=gray]
\tikzstyle{aedge} = [draw, thin, ->,black]
\tikzstyle{edge} = [draw, thick, -,black]
\tikzstyle{doubleedge} = [draw, thick, double distance=1pt, -,black]
\tikzstyle{hiddenedge} = [draw=none, thick, double distance=1pt, -,black]
\tikzstyle{dashededge} = [draw, very thick, dashed, black]
\tikzstyle{LabelStyleH} = [text=black, anchor=south]
\tikzstyle{LabelStyleHn} = [text=black, anchor=north]
\tikzstyle{LabelStyleV} = [text=black, anchor=east]
\node[vertex] (v1) at (2,1){\footnotesize$\e{12413} $};
\node[framedvertex] (v2) at (0,1){\footnotesize$\e{41213} $};
\node[vertex] (v3) at (3,1){\footnotesize$\e{12431} $};
\node[vertex] (v4) at (1,1){\footnotesize$\e{14213} $};
\draw[edge] (v1) to node[LabelStyleH]{\Tiny$4 $} (v3);
\draw[edge] (v1) to node[LabelStyleH]{\Tiny$3 $} (v4);
\draw[edge] (v2) to node[LabelStyleH]{\Tiny$2 $} (v4);
\node[vertex] (v) at (-0.9,0.75){\small$  \  $};
\node[vertex,anchor = west] (vv) at (4.9,1){\small$F_{\Gamma_4} = s_{41}$};
\end{tikzpicture}
\end{center}\vspace{-.3mm}

\begin{center}%
\begin{tikzpicture}[xscale = 2.2,yscale = 1.5]
\tikzstyle{vertex}=[inner sep=0pt, outer sep=4pt]
\tikzstyle{framedvertex}=[inner sep=3pt, outer sep=4pt, draw=gray]
\tikzstyle{aedge} = [draw, thin, ->,black]
\tikzstyle{edge} = [draw, thick, -,black]
\tikzstyle{doubleedge} = [draw, thick, double distance=1pt, -,black]
\tikzstyle{hiddenedge} = [draw=none, thick, double distance=1pt, -,black]
\tikzstyle{dashededge} = [draw, very thick, dashed, black]
\tikzstyle{LabelStyleH} = [text=black, anchor=south]
\tikzstyle{LabelStyleHn} = [text=black, anchor=north]
\tikzstyle{LabelStyleV} = [text=black, anchor=east]
\node[framedvertex] (v1) at (1,1){\footnotesize$\e{21134} $};
\node[framedvertex] (v2) at (2,1){\footnotesize$\e{12143} $};
\node[framedvertex] (v3) at (3,1){\footnotesize$\e{11432} $};
\node[framedvertex] (v4) at (4,1){\footnotesize$\e{13214} $};
\node[vertex] (v) at (-0.2,0.75){\small$  \  $};
\node[vertex,anchor = west] (vv) at (4.9,1){\small$F_{\Gamma_5}= F_{\Gamma_6}= F_{\Gamma_7} = F_{\Gamma_8} = s_{5}$};
\end{tikzpicture}
\end{center}
\caption{\label{f P switchboard 11234}
For  $P = \PP_{2}$, the $P$-Knuth equivalence graph on the set of words  $\e{w}$ of content
$\{1,1,2,3,4\}$ and with  $\inv_P(\e{w}) = 2$, where  $\inv_P$ is as defined in \eqref{e inv def}.
The symmetric functions of its components can be read off from the
framed vertices as in Figure \ref{f P switchboard 1234}.
}
\end{figure}
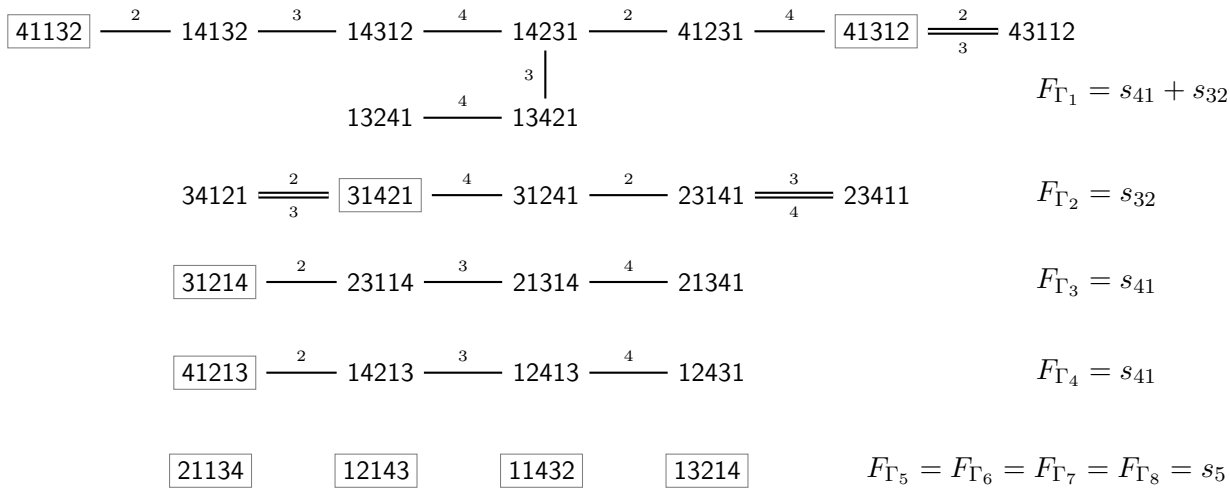

\subsection{Schur positivity of $P$-Knuth equivalence graphs}
\label{ss Schur pos P Knuth}

The \emph{Young diagram} of a partition $\lambda= (\lambda_1 \ge \cdots \ge
\lambda_k > 0)$ is the subset of boxes
$\{(r,c) \in \ZZ_{\geq 1} \times \ZZ_{\geq 1} :  c \le \lambda_r\}$
in the plane, drawn in English (matrix-style) notation so that row (resp. column) indices
increase from north to south (resp. west to east).
We write $\ell(\lambda )$ for the number of nonzero parts of  $\lambda$.

We say that a word  $\e{w} = \e{w}_1 \cdots \e{w}_n \in \U_P^*$ is  \emph{$P$-strictly decreasing} (resp.  \emph{$P$-weakly increasing}) if  $\e{w}_1 \leftto_P \e{w}_2 \leftto_P \cdots \leftto_P \e{w}_n$
(resp. $\e{w}_1 \dotto_P \e{w}_2 \dotto_P \cdots \dotto_P \e{w}_n$).

We use the shorthand  $[N]$ for the set  $\{1,\dots, N\}$.

\begin{definition}
\label{def:poset SSYT}
For a poset $P$, a \emph{$P$-tableau of shape $\lambda$} is a filling $T \colon \lambda \to P$ of the Young diagram of  $\lambda$ such that
\begin{itemize}
\item Each column of  $T$, read bottom to top, is  $P$-strictly decreasing, i.e., for all  $c \in [\lambda_1]$,  $T(\lambda'_c,c) \leftto_P   \cdots \leftto_P T(2,c) \leftto_P T(1,c)$.
\item Each row of  $T$, read left to right is $P$-weakly increasing, i.e., for all  $r\in [\ell(\lambda)]$,  $T(r,1) \dotto_P T(r,2)  \dotto_P \cdots \dotto_P  T(r,\lambda_r)$.
\end{itemize}
Let $\SSYT_P(\lambda)$ denote the set of  $P$-tableaux of shape $\lambda$.

The \emph{content} of a  $P$-tableau  $T$ is the multiset  $\beta$ of elements of $P$ in which the multiplicity of
$a$ in  $\beta$ is the number of occurrences of $a$ in  $T$.
A  \emph{standard $P$-tableau} is  a  $P$-tableau of content  $\beta = P$.
\end{definition}

\begin{remark}
Our  $P$-tableaux are the transpose of those defined in \cite{Gasharov}.
Our standard $P$-tableaux are the same as the  $P$-tableaux in \cite{KPchromatic} and are the transpose of the
$P$-tableaux in \cite{SWchromatic}.
\end{remark}

The \emph{column reading word} $\creading(T)\in\U^*_P$
of a tableau $T \in \SSYT_P(\lambda)$
is the word obtained by concatenating the columns of~$T$ (reading
each column bottom to top), starting with the leftmost column.
For example, for  $P=\PP_{2}$, the  $P$-tableau
\[
T=  \fontsize{8pt}{6pt}\selectfont \tableau{1&2&1&1 \\3&4&3\\5&6&7}
\]
of shape $\lambda=(433)$ has column reading word $\creading(T)=\e{5316427311}$
and content  $\{ 1,1,1, 2,3,3,4,5,6,7\}$.

\begin{theorem}
\label{th:IplacP-positivity}
For a \threeone-free poset $P$, the noncommutative $P$-Schur functions have the following  $\mathbf{u}$-monomial positive expression modulo $\Iplacp{P}$:
\[\mathfrak{J}_\lambda^P(\mathbf{u}) \equiv \sum_{T \in \SSYT_P(\lambda)} \mathbf{u}_{\creading(T)}.\]
\end{theorem}
This will be proved in Section \ref{sec:proof-of-IplacP-positivity}.
Combining Theorem \ref{th:IplacP-positivity} and Corollary \ref{c positivity} (i) establishes our main result, Theorem \ref{t P KE schur pos intro}, which we restate here:

\begin{theorem}
\label{c P KE schur pos}
For any  $P$-Knuth equivalence graph  $\Gamma$,
the symmetric function $F_\Gamma(\mathbf{x})$ is Schur positive.
The coefficient of $s_\lambda(\mathbf{x})$ in  $F_\Gamma(\mathbf{x})$ is
the number of  $P$-tableaux $T\in \SSYT_P(\lambda)$ such that
the column reading word $\creading(T)$ appears as a vertex in~$\Gamma$.
\end{theorem}

This settles \cite[Conjecture 4.14]{KPchromatic} and generalizes it in two directions---from natural unit interval orders $P$ to \threeone-free posets  $P$ and from standard $P$-Knuth equivalence graphs to arbitrary ones.
It can be regarded as a refinement of results of Gasharov and Shareshian-Wachs on the Schur expansion of  ($t$-)chromatic symmetric functions, as will be explained below in \S\ref{ss t chromatic}.


\subsection{Relation to  $t$-chromatic symmetric functions}
\label{ss t chromatic}
There is a natural and clean way to obtain $t$-chromatic symmetric functions from the  $P$-plactic algebra,
and we use this to show how to recover two results of Shareshian-Wachs from our approach.
This can be regarded as a summary and expansion of similar material in \cite[\S\S4.3--4.4]{KPchromatic}.

\begin{definition}
\label{def:interval order}
A \emph{natural unit interval order} is a poset  $P$ (with relation denoted  $a \to_P b$ as usual)
equipped with a total order $<$ on the set $P$,  satisfying
\begin{itemize}
\item[(i)]  $a \to_P b$ implies  $a< b$;
\item[(ii)] If  \spa
\begin{tikzpicture}[xscale = 1.8,yscale = 1.5]
\tikzstyle{vertex}=[inner sep=0pt, outer sep= 2.5pt]
\tikzstyle{aedge} = [draw, ->,>=stealth', black]
\tikzstyle{aedgecurve} = [draw, ->,>=stealth', black, bend left=40]
\tikzstyle{edge} = [draw, thick, -,black]
\tikzstyle{dashededge} = [draw, -,  dashed, black]
\node[vertex] (va) at (0,-0.8){\small${a}$};
\node[vertex] (vb) at (0.5,-0.8){\small${b}$};
\node[vertex] (vc) at (1,-0.8){\small${c}$};
\draw[dashededge] (va) to (vb);
\draw[dashededge] (vb) to (vc);
\draw[aedgecurve] (va) to (vc);
\end{tikzpicture}
\spa is an induced subposet of  $P$, then $a < b < c$.
\end{itemize}
\end{definition}

The posets  $\PP_{k}$ from Definition \ref{def:jump poset} are natural unit interval orders with total order the usual one on  $[N]$.

\begin{definition}[{\cite[Definition 1.2]{SWchromatic}}]
For a natural unit interval order  $P$, the  \emph{$t$-chromatic symmetric function} of  $\inc(P)$ is
\begin{align}
\label{eq:def chi t}
X_{\inc(P)}(\mathbf{x},t) = \sum_\kappa t^{\asc_P(\kappa)} \prod_{a \in P} x_{\kappa(a)}
\end{align}
where the sum is over all proper colorings $\kappa$ of  $\inc(P)$ and
\begin{align}
\asc_P(\kappa) \defeq \big| \big\{ a, b \in P : a< b,\, a \dote_P b, \, \text{ and } \kappa(a) < \kappa(b) \big\} \big|.
\end{align}
\end{definition}

It is natural from our perspective to also consider a multicolored version of ($t$-)chromatic symmetric functions.

\begin{definition}
Let  $P$ be a poset.  Let $\beta$ be a multiset of elements of  $P$ and  $\mathbf{n} = (n_a)_{a \in P}$ the corresponding vector of multiplicities.
Let  $P[\beta]$ denote the poset obtained from  $P$ by replacing each element  $a\in P$ with  $n_a$ copies  $a^{(1)}, \dots,  a^{(n_a)}$ of  $a$ and relations determined by
$a^{(i)} \to_{P[\beta]} b^{(j)}$ if and only if  $a \to_P b$ (hence the copies $a^{(i)}$ and  $a^{(j)}$ are incomparable for $i \ne j$).
If  $P$ is a natural unit interval order, we define $P[\beta]$ to be the natural unit interval order
whose underlying poset is as just defined and with total order  $<$
determined by (1) $a^{(1)} < \cdots < a^{(n_a)}$, and (2) $a^{(i)} < b^{(j)}$ whenever
$a < b$ in  $P$.

The \emph{$\beta$-multicolored chromatic symmetric function} of  $\inc(P)$ is
\begin{align}
\label{eq:def chi beta}
X_{\inc(P)}^{\beta}(\mathbf{x}) \, \defeq \, \sum_\kappa \prod_{b \in P[\beta]} x_{\kappa(b)},
\end{align}
where the sum is over all proper colorings  $\kappa$ of  $\inc(P[\beta])$ such that
$\kappa(a^{(1)}) > \cdots > \kappa(a^{(n_a)})$ for each  $a \in P$.
If  $P$ is a natural unit interval order,
the \emph{$\beta$-multicolored  $t$-chromatic symmetric function} of  $\inc(P)$ is
\begin{align}
\label{eq:t beta chi}
X_{\inc(P)}^{\beta}(\mathbf{x},t) \, \defeq \, \sum_\kappa t^{\asc_{P[\beta]}(\kappa)} \prod_{b \in P[\beta]} x_{\kappa(b)},
\end{align}
where the sum is just as in \eqref{eq:def chi beta}.
Note that with  $\beta = P$, we recover the usual \break ($t$-)chromatic symmetric functions.
\end{definition}

\begin{remark}
The  multicolored chromatic symmetric functions were introduced by Gasharov \cite{Gasharov}.
A coloring  $\kappa$ as in \eqref{eq:def chi beta} is essentially the same as
a proper  $\mathbf{n}$-multicoloring of  $\inc(P)$ defined in \cite{Gasharov}; the latter is by definition
a map  $\tilde{\kappa}$ from the vertices of $\inc(P)$ to the set of subsets of  $\ZZ_{\ge 1}$ such that  $|\tilde{\kappa}(a)| = n_a$ for all  $a\in P$ and  $\tilde{\kappa}(a) \cap \tilde{\kappa}(b) = \varnothing$ for all  $a \dote_P b$, and the correspondence  $\kappa \mapsto \tilde{\kappa}$ is given by setting
$\tilde{\kappa}(a) = \{\kappa(a^{(1)}), \dots, \kappa(a^{(n_a) } ) \}$.
\end{remark}

\begin{remark}
\label{r beta}
Let  $P$ be a poset and  $\beta$ and  $\mathbf{n} = (n_a)_{a \in P}$ be as above.
It is easily seen that
\begin{align}
X_{\inc(P)}^{\beta}(\mathbf{x}) =  \frac{1}{\prod_{a\in P} n_a!}\spa X_{\inc(P[\beta])}(\mathbf{x})
\end{align}
since each coloring of  $\inc(P[\beta])$ as in \eqref{eq:def chi beta} corresponds naturally to  $\prod_{a\in P} n_a!$  many colorings of $\inc(P[\beta])$.
For a natural unit interval order $P$,
it follows similarly using \cite[Example 2.4]{SWchromatic} that
\begin{align}
X_{\inc(P)}^{\beta}(\mathbf{x},t) =  \frac{1}{\prod_{a\in P} [n_a]_t!}\spa X_{\inc(P[\beta])}(\mathbf{x}, t),
\end{align}
where  $[d]_t = 1+t + \cdots +t^{d-1}$ and  $[n]_t! = \prod_{d = 1}^n [d]_t$.
Therefore the multicolored variants are not really any more general.  However, we mention
them since they arise naturally in our approach and they form important examples of the elements
$\omega F_\gamma (\mathbf{x})$.  
\end{remark}

Let $\QQ(t)$ denote the field of
rational functions in one variable~$t$.
The basic setup of \S\ref{s sec 2} extends
to the setting where objects 
are defined over the ground ring $\QQ(t)$ rather than  $\ZZ$.

\begin{definition}
For a natural unit interval order  $P$,
let $I_T^P$ denote the ideal of  $\QQ(t) \tsr_{\ZZ} \U_P$ generated by
\begin{alignat}{4}
&u_c u_a - u_a u_c & &  \qquad  (a \to_P c), \label{poset rel e pos v0} \\
&u_b u_a - t \spa u_a u_b & &  \qquad (a < b \ \text{ and } \  a \dote_P b). \label{poset rel t}
\end{alignat}
\end{definition}

\begin{proposition}
\label{pr:Itcomm and plac}
For a natural unit interval order  $P$, $\QQ(t) \tsr \Iplacp{P} \subset I_T^P$.
\end{proposition}
\begin{proof}
The elements  \eqref{poset rel knuth etc bca} and \eqref{poset rel knuth etc cab} are obtained by multiplying \eqref{poset rel e pos v0} on the left and right by  $u_b$.
To see that \eqref{poset rel bca cab} lies in  $I_T^P$, by Definition \ref{def:interval order} (ii),
$a \dote_P b \dote_P c$ and $a \to_P c$ implies  $a < b < c$; hence
  $u_c u_a u_b- u_b u_c u_a  = u_c (u_a u_b - t^{-1}u_bu_a) - (u_b u_c-t^{-1}u_c u_b)u_a \in I_T^P$.
\end{proof}

Let  $P$ be a natural unit interval order.
For a word  $\e{w} = \e{w}_1 \cdots \e{w}_n \in \U^*_P$, define
\begin{align}
\label{e inv def}
\inv_P(\e{w}) = \big| \big\{ (\e{w}_i,\e{w}_j) \, : \, i< j,\, \e{w}_i > \e{w}_j, \text{ and } \e{w}_i \dote_P \e{w}_j \big\} \big|.
\end{align}
For a multiset  $\beta$ of elements of  $P$,
define the following element of  $\QQ(t) \tsr \U^*_P$
\begin{align}
\Wo{\beta}(t)  = \sum_{\text{words $\e{w}$ of content  $\beta$}} \!\! t^{\inv_P(\e{w})} \spa \e{w}.
\end{align}
Note that the  $t=1$ specialization of $\Wo{\beta}(t)$ is the same as $\Wo{\beta}$ defined in \eqref{eq: Wo beta}.

\begin{proposition}
\label{pr:Wo beta t}
Maintain the notation above.
The element  
$\Wo{\beta}(t)$ lies in $(I_T^P)^\perp$. 
Moreover, the set of $\Wo{\beta}(t)$ over all multisets $\beta$ forms a  $\QQ(t)$-basis for
$(I_T^P)^\perp \subset \QQ(t) \tsr \U^*_P$.
\end{proposition}
\begin{proof}
Checking $\Wo{\beta}(t) \in (I_T^P)^\perp$ amounts to showing that for
words $\e{w}=\e{w}_1 \cdots \e{w}_n$ and  $\e{w}' = \e{w}_1 \cdots\e{w}_{i-1} \e{w}_{i+1} \e{w}_i \e{w}_{i+2} \cdots  \e{w}_n$,
the difference $\inv_P(\e{w})-\inv_P(\e{w}')$ is 1 if  $\e{w}_i > \e{w}_{i+1}$ and  $\e{w}_i \dote_P \e{w}_{i+1}$, and is 0 if  $\e{w}_i \leftto_P \e{w}_{i+1}$. This is clear.
 The $\Wo{\beta}(t)$ are clearly linearly independent, while to see that they span, note that for any words  $\e{v}$, $\e{v}'$ of the same content there is some element $\e{v}-t^k\e{v}' \in I_T^P$; hence
 any $\sum_{\e{w}} \gamma_{\e{w}} \,\e{w}\in (I_T^P)^\perp$
 is determined by selecting one  $\gamma_{\e{w}}$ from each bin $\{\gamma_{\e{w}} : \text{content}(\e{w}) = \beta\}$, for  each multiset  $\beta$.
\end{proof}

Recall Chow's result (Theorem \ref{th:chromatic as F}) which expresses a chromatic symmetric function in terms of quasisymmetric functions.
Shareshian-Wachs \cite[Theorem 3.1]{SWchromatic} generalized this to address the $t$-chromatic symmetric functions.
Our method gives a short proof of these results and their multicolored generalizations,
which demonstrates the utility of the  noncommutative $P$-Cauchy formulas.

\begin{theorem}
For any poset  $P$ and multiset  $\beta$ of elements of  $P$,
\begin{align}
\label{eq:F W beta}
F_{\Wo{\beta}}(\mathbf{x}) = \omega X_{\inc(P)}^{\beta}(\mathbf{x}).
\end{align}
If  $P$ is a natural unit interval order, then
\begin{align}
\label{eq:F W beta t}
F_{\Wo{\beta}(t)}(\mathbf{x}) = \omega X_{\inc(P)}^{\beta}(\mathbf{x},t).
\end{align}
\end{theorem}

\begin{proof}
We first prove \eqref{eq:F W beta}.  Plugging in  $\gamma = \Wo{\beta}$ into \eqref{eq:F-via-e} yields
\begin{equation}
\label{eq:P Cauchy m e 2}
\omega \spa F_{\Wo{\beta}}(\mathbf{x}) =
\sum_{\lambda} m_\lambda(\mathbf x) \langle e^P_{\lambda}(\mathbf{u}), \Wo{\beta} \rangle.
\end{equation}
The integer $\langle e^P_{\lambda}(\mathbf{u}), \Wo{\beta} \rangle$ is the number of ordered partitions of $\beta$ into chains of lengths  $\lambda_1, \dots, \lambda_\ell$ (where  $\ell$ is the length of  $\lambda$), that is, the number of  $\ell$-tuples $(C_1, \dots, C_\ell)$
such that $C_i$ is a chain of  $P$ of size $\lambda_i$ for all  $i$ and the multiset union of the  $C_i$ is  $\beta$.  This is exactly the coefficient of  $m_\lambda(\mathbf{x})$ in $X_{\inc(P)}^{\beta}(\mathbf{x})$.

Formula \eqref{eq:F W beta t} follows similarly: the integer
$\langle e^P_{\lambda}(\mathbf{u}), \Wo{\beta}(t) \rangle$ is now the sum  $\sum_{\mathbf{C}} t^{\inv_P(\e{w}(\mathbf{C}))}$
over ordered partitions  $\mathbf{C} = (C_1,\dots, C_\ell)$ of $\beta$ as above, where  $\e{w}(\mathbf{C})$ is the word obtained by concatenating the  $P$-strictly decreasing words corresponding to the chains $C_1$,  $C_2, \dots, C_\ell$.
On the other hand, noting that the coloring  $\kappa$ of  $\inc(P[\beta])$ in which  $C_i$ receives color  $\ell-i+1$ satisfies  $\asc_{P[\beta]}(\kappa) = \inv_P(\e{w}(\mathbf{C}))$, this integer is exactly the coefficient of  $x_1^{\lambda_\ell} \cdots x_\ell^{\lambda_1}$
 in $X_{\inc(P)}^{\beta}(\mathbf{x},t)$, and this is the same as
 the coefficient of $m_\lambda(\mathbf{x})$ in $X_{\inc(P)}^{\beta}(\mathbf{x},t)$.
\end{proof}

Gasharov \cite{Gasharov} gave a formula for the Schur expansion of  
$X_{\inc(P)}^\beta(\mathbf{x})$ for a \threeone-free poset  $P$ in terms of  $P$-tableaux
and Shareshian-Wachs \cite[Theorem 6.3]{SWchromatic} gave a similar formula for the  $t$-chromatic symmetric functions $X_{\inc(P)}(\mathbf{x}, t)$. 
We now show how to obtain these results from our
Theorem \ref{th:IplacP-positivity} or  \ref{c P KE schur pos}.

\begin{corollary}
\label{c chi beta t Schur pos}
For a \threeone-free poset  $P$
and any multiset $\beta$ of elements of  $P$,
\begin{align}
\label{eq:c chi beta Schur pos}
\omega\spa X_{\inc(P)}^{\beta}(\mathbf{x}) = \sum_{T\in \SSYT_P^\beta } s_{\sh(T)}(\mathbf{x}),
\end{align}
where  $\SSYT_P^\beta$ denotes the set of  $P$-tableaux of content  $\beta$.
If  $P$ is a natural unit interval order, then
\begin{align}
\label{eq:c chi beta t Schur pos}
\omega\spa  X_{\inc(P)}^{\beta}(\mathbf{x};t) = \sum_{T\in \SSYT_P^\beta } t^{\inv_P(\creading(T))}s_{\sh(T)}(\mathbf{x}).
\end{align}
\end{corollary}
\begin{proof}
By Propositions \ref{pr:Itcomm and plac} and \ref{pr:Wo beta t}, $\Wo{\beta}(t) \in (\QQ(t)\tsr \Iplacp{P})^\perp$.
Then
apply Theorem \ref{th:IplacP-positivity} and \eqref{eq:F-via-schurs} with  $\gamma = \Wo{\beta}(t)$ and combine with \eqref{eq:F W beta t}.  This yields \eqref{eq:c chi beta t Schur pos}.
Formula \eqref{eq:c chi beta Schur pos} can be proved in a similar way.
\end{proof}

\begin{remark}
For the usual  $t$-chromatic symmetric function ($\beta = P$),
\eqref{eq:c chi beta t Schur pos} is the same as \cite[Theorem 6.3]{SWchromatic}, and for general  $\beta$, it can be deduced easily from this case using Remark \ref{r beta}.
\end{remark}

\begin{remark}
Corollary \ref{c chi beta t Schur pos} can also be deduced
from Theorem \ref{c P KE schur pos} by observing that the
set of all words  $\e{w}$ of content  $\beta$ and with $\inv_P(\e{w}) = d$ is a  $P$-Knuth equivalence graph.
\end{remark}

\subsection{On the  $e$-positivity of chromatic symmetric functions}
\label{ss e pos}

Following up on the general discussion at the beginning of \S\ref{s P Knuth and H graphs},
towards the Stanley-Stembridge conjecture, we seek as an ideal  $I$ as small as possible such
that the $\omega\spa F_\gamma(\mathbf{x})$ are $e$-positive for any  $\gamma \in (\U_P^*)_{\ge 0} \cap I^\perp$.
Since this is not true for  $\Iplacp{P}$
(see Figures \ref{f P switchboard 1234}, \ref{f 2p2}, and \ref{f P switchboard 11234})
we consider the following larger ideal introduced by Hwang \cite{Hwang}
as a setting to study the Stanley-Stembridge conjecture
(we work in the slightly more general setting than Hwang where  $P$ is a \threeone-free poset rather than a natural unit interval order).

\begin{definition}
For a  \threeone-free poset  $P$,
let $I_H^P$ denote the ideal of  $\U_P$ generated by
\begin{alignat}{4}
&u_c u_a - u_a u_c & & \quad  (a \to_P c), \label{poset rel e pos} \\
&u_c u_a u_b - u_b u_c u_a  & & \quad(a \dote_P b \dote_P c \ \text{ and } \  a \to_P c). \label{poset rel bca cab 2}
\end{alignat}
\end{definition}

\begin{proposition}
\label{p IH inclusions}
The ideals introduced in this paper satisfy
\begin{align}
\label{eq:p inclusions 1}
& \Icommp{P} \subset \Iplacp{P}\subset  I_H^P \subset I_\text{\rm pol}^P \, , \\
\label{eq:p inclusions 2}
& \QQ(t) \tsr I_H^P \subset I_T^P.
\end{align}
where \eqref{eq:p inclusions 1} holds for \threeone-free posets $P$ and \eqref{eq:p inclusions 2} for natural unit interval orders $P$.
\end{proposition}
\begin{proof}
The first inclusion of \eqref{eq:p inclusions 1} is Theorem \ref{th:IplacP elem commute}, the second holds since
the generators \eqref{poset rel knuth etc bca} and \eqref{poset rel knuth etc cab}
are multiples of \eqref{poset rel e pos}, and the third is clear.
For \eqref{eq:p inclusions 2}, we need to check $u_c u_a u_b - u_b u_c u_a  \in I_T^P$, which was already done in the proof of Proposition \ref{pr:Itcomm and plac}.
\end{proof}

Just like in \S\ref{ss P Knuth}, we can define graphs corresponding to elements of  $(I_H^P)^\perp$.

\begin{definition}[\emph{$H$-graphs} \cite{Hwang}]
\label{def-Hgraph}
Let  $P$ be a \threeone-free poset.
Define $\Gamma^n$ to be the  graph (undirected, without edge labels)
on the vertex set of words of length~$n$
in the alphabet $P$ and with the following two types of edges:
\begin{align*}
&
\begin{tikzpicture}[xscale = 3,yscale = 1.5]
\tikzstyle{vertex}=[inner sep=0pt, outer sep=4pt]
\tikzstyle{framedvertex}=[inner sep=3pt, outer sep=4pt, draw=gray]
\tikzstyle{aedge} = [draw, thin, ->,black]
\tikzstyle{edge} = [draw, thick, -,black]
\tikzstyle{LabelStyleH} = [text=black, anchor=south]
\node[vertex] (v1) at (0,1){$\, \e{\cdots ca\cdots} \, $};
\node[vertex] (v2) at (1,1){$\, \e{\cdots ac\cdots} \, $};
\node[vertex,anchor=west] (v3) at (1.6,1){with $a \to_P c$ \ (comparable edge)};
\draw[edge] (v1) to node[LabelStyleH]{} (v2);
\end{tikzpicture}
\\[1mm]
&
\begin{tikzpicture}[xscale = 3,yscale = 1.5]
\tikzstyle{vertex}=[inner sep=0pt, outer sep=4pt]
\tikzstyle{framedvertex}=[inner sep=3pt, outer sep=4pt, draw=gray]
\tikzstyle{aedge} = [draw, thin, ->,black]
\tikzstyle{edge} = [draw, thick, -,black]
\tikzstyle{LabelStyleH} = [text=black, anchor=south]
\node[vertex] (v1) at (0,0){$\e{\cdots cab\cdots} $};
\node[vertex] (v2) at (1,0){$\e{\cdots bca\cdots} $};
\node[vertex,anchor=west] (v3) at
(1.6,-.04){\parbox{9cm}{with $a \dote_P b \dote_P c \ \text{ and } \  a \to_P c$ \\[1mm] (a $P$-Knuth transformation
as in Def. \ref{def-switchboard} (3))}};
\draw[edge] (v1) to node[LabelStyleH]{} (v2);
\end{tikzpicture}
\end{align*}
More precisely, words $\e{w}= \e{w}_1\cdots \e{w}_n \in \U_P^*$ and $\e{w'}= \e{w}_1'\cdots
  \e{w}_n' \in \U_P^*$  are the ends of a \emph{comparable edge}
if there is an index  $i$ such that
$\e{w}_j=\e{w}_j'$ for $j\notin\{i,i+1\}$
and $\{\e{w_iw_{i+1}},
\e{w}'_i\e{w}'_{i+1}\} = \{\e{ca}, \e{ac} \}$  with $a \to_P c$.
An \emph{$I_H^P$-graph}, or just \emph{$H$-graph} when  $P$ is clear from context,
is a subgraph of  $\Gamma^n$ which is a disjoint union of connected
components of  $\Gamma^n$.
\end{definition}

\begin{figure}
        \centerfloat
\begin{tikzpicture}[xscale = 1.77,yscale = 1.4]
\tikzstyle{vertex}=[inner sep=0pt, outer sep=2.7pt]
\tikzstyle{framedvertex}=[inner sep=3pt, outer sep=4pt, draw=gray]
\tikzstyle{aedge} = [draw, thin, -,black]
\tikzstyle{thickedge} = [draw,very thick, -,black]
\tikzstyle{thinedge} = [draw, thin, -,black!48]
\tikzstyle{edge} = [draw, thick, -,black]
\tikzstyle{doubleedge} = [draw, thick, double distance=1pt, -,black]
\tikzstyle{hiddenedge} = [draw=none, thick, double distance=1pt, -,black]
\tikzstyle{dashededge} = [draw, very thick, dashed, black]
\tikzstyle{LabelStyleH} = [text=black, anchor=south]
\tikzstyle{LabelStyleHn} = [text=black, anchor=north]
\tikzstyle{LabelStyleV} = [text=black, anchor=east]

\node[vertex] (v1) at (4,4){\scriptsize$\e{26413} $};
\node[framedvertex] (v2) at (5,3){\scriptsize$\e{42136} $};
\node[vertex] (v3) at (2,3){\scriptsize$\e{26314} $};
\node[framedvertex] (v4) at (5,4){\scriptsize$\e{21643} $};
\node[vertex] (v5) at (3,1){\scriptsize$\e{64123} $};
\node[framedvertex] (v6) at (2,1){\scriptsize$\e{61423} $};
\node[framedvertex] (v7) at (1,4){\scriptsize$\e{23164} $};
\node[vertex] (v8) at (4,5){\scriptsize$\e{24361} $};
\node[vertex] (v9) at (4,1){\scriptsize$\e{46213} $};
\node[vertex] (v10) at (3,4){\scriptsize$\e{26341} $};
\node[framedvertex] (v11) at (4,2){\scriptsize$\e{62143} $};
\node[framedvertex] (v12) at (5,5){\scriptsize$\e{24316} $};
\node[framedvertex] (v13) at (3,2){\scriptsize$\e{62413} $};
\node[framedvertex] (v14) at (2,2){\scriptsize$\e{62314} $};
\node[vertex] (v15) at (3,3){\scriptsize$\e{62341} $};
\node[vertex] (v16) at (2,4){\scriptsize$\e{23614} $};
\node[vertex] (v17) at (5,1){\scriptsize$\e{42613} $};
\node[vertex] (v18) at (4,3){\scriptsize$\e{26143} $};
\node[vertex] (v19) at (5,2){\scriptsize$\e{42163} $};
\node[vertex] (v20) at (2,5){\scriptsize$\e{23641} $};
\node[framedvertex] (v21) at (1,1){\scriptsize$\e{16423} $};
\node[vertex] (v22) at (3,5){\scriptsize$\e{24631} $};

\node[vertex] (vv) at (3,.33){\footnotesize $F_{\Gamma_1} = 2h_5+2h_{41}+h_{32}$};

\draw[thinedge] (v1) to (v13);
\draw[thickedge] (v1) to (v18);
\draw[thickedge] (v1) to (v10);
\draw[thinedge] (v2) to (v19);
\draw[thickedge] (v3) to (v14);
\draw[thickedge] (v3) to (v16);
\draw[thinedge] (v3) to (v10);
\draw[thinedge] (v4) to (v18);
\draw[thickedge] (v5) to (v6);
\draw[thickedge] (v5) to (v13);
\draw[thinedge] (v6) to (v21);
\draw[thinedge] (v7) to (v16);
\draw[thickedge] (v8) to (v22);
\draw[thinedge] (v8) to (v12);
\draw[thickedge] (v9) to (v17);
\draw[thickedge] (v9) to (v13);
\draw[thickedge] (v10) to (v22);
\draw[thinedge] (v10) to (v20);
\draw[thickedge] (v10) to (v15);
\draw[thickedge] (v11) to (v18);
\draw[thinedge] (v11) to (v13);
\draw[thickedge] (v13) to (v15);
\draw[thinedge] (v14) to (v15);
\draw[thickedge] (v16) to (v20);
\draw[thickedge] (v17) to (v19);

\end{tikzpicture}
\quad \ \ \ \ \
        \centerfloat
\begin{tikzpicture}[xscale = 1.77,yscale = 1.4]
\tikzstyle{vertex}=[inner sep=0pt, outer sep=2.7pt]
\tikzstyle{framedvertex}=[inner sep=3pt, outer sep=4pt, draw=gray]
\tikzstyle{thickedge} = [draw,very thick, -,black]
\tikzstyle{thinedge} = [draw, thin, -,black!48]

\node[vertex] (v1) at (4,4){\scriptsize$\e{14632} $};
\node[vertex] (v2) at (3,2){\scriptsize$\e{34162} $};
\node[framedvertex] (v3) at (5,3){\scriptsize$\e{63124} $};
\node[vertex] (v4) at (1,3){\scriptsize$\e{32416} $};
\node[vertex] (v5) at (2,2){\scriptsize$\e{34126} $};
\node[vertex] (v6) at (4,3){\scriptsize$\e{14362} $};
\node[framedvertex] (v7) at (3,4){\scriptsize$\e{41632} $};
\node[framedvertex] (v8) at (5,1){\scriptsize$\e{31264} $};
\node[framedvertex] (v9) at (2,3){\scriptsize$\e{41326} $};
\node[vertex] (v10) at (1,2){\scriptsize$\e{32461} $};
\node[vertex] (v11) at (3,3){\scriptsize$\e{41362} $};
\node[vertex] (v12) at (5,4){\scriptsize$\e{16342} $};
\node[vertex] (v13) at (3,1){\scriptsize$\e{31462} $};
\node[vertex] (v14) at (5,2){\scriptsize$\e{36124} $};
\node[framedvertex] (v15) at (5,5){\scriptsize$\e{61342} $};
\node[framedvertex] (v16) at (4,5){\scriptsize$\e{13642} $};
\node[framedvertex] (v17) at (2,1){\scriptsize$\e{31426} $};
\node[vertex] (v18) at (4,1){\scriptsize$\e{31624} $};
\node[framedvertex] (v19) at (1,4){\scriptsize$\e{32146} $};
\node[vertex] (v20) at (4,2){\scriptsize$\e{34612} $};
\node[framedvertex] (v21) at (2,4){\scriptsize$\e{14326} $};
\node[vertex] (v22) at (3,5){\scriptsize$\e{46132} $};

\node[vertex] (vv) at (3,.33){\footnotesize $F_{\Gamma_2} = 2h_5+2h_{41}+h_{32}$};

\draw[thickedge] (v1) to (v12);
\draw[thickedge] (v1) to (v6);
\draw[thinedge] (v1) to (v7);
\draw[thinedge] (v2) to (v5);
\draw[thickedge] (v2) to (v20);
\draw[thickedge] (v2) to (v11);
\draw[thinedge] (v2) to (v13);
\draw[thickedge] (v3) to (v14);
\draw[thinedge] (v4) to (v19);
\draw[thickedge] (v4) to (v10);
\draw[thickedge] (v4) to (v5);
\draw[thickedge] (v5) to (v9);
\draw[thinedge] (v5) to (v17);
\draw[thinedge] (v6) to (v21);
\draw[thinedge] (v6) to (v11);
\draw[thickedge] (v7) to (v22);
\draw[thickedge] (v7) to (v11);
\draw[thinedge] (v8) to (v18);
\draw[thinedge] (v9) to (v21);
\draw[thinedge] (v9) to (v11);
\draw[thinedge] (v12) to (v16);
\draw[thickedge] (v12) to (v15);
\draw[thickedge] (v13) to (v18);
\draw[thinedge] (v13) to (v17);
\draw[thickedge] (v14) to (v18);
\end{tikzpicture}
\caption{\label{fig:H disconnected}
\small
For  $P =\PP_{3}$, $\beta = \{1,2,3,4,6\}$,
and  $d=3$, the $H$-graph  $\Gamma$ on the set of words
$\big\{\e{w} :  \text{\rm content}(\e{w}) = \beta, \, \inv_P(\e{w}) = d \big\}$,
with components $\Gamma^{1}$ and  $\Gamma^{2}$.
This graph is disconnected which shows that Conjecture \ref{cj:m pos} is stronger than the Shareshian-Wachs  $e$-positivity conjecture.
The thick edges are the edges which are also  $P$-Knuth transformations.  The Schur expansion of each component can be read off from the \boxed{\text{vertices}}
which are column reading words of $P$-tableaux, while finding a similar way to read off the  $h$-expansion of any $H$-graph is tantamount to solving Conjecture \ref{cj:m pos}.
}
\end{figure}
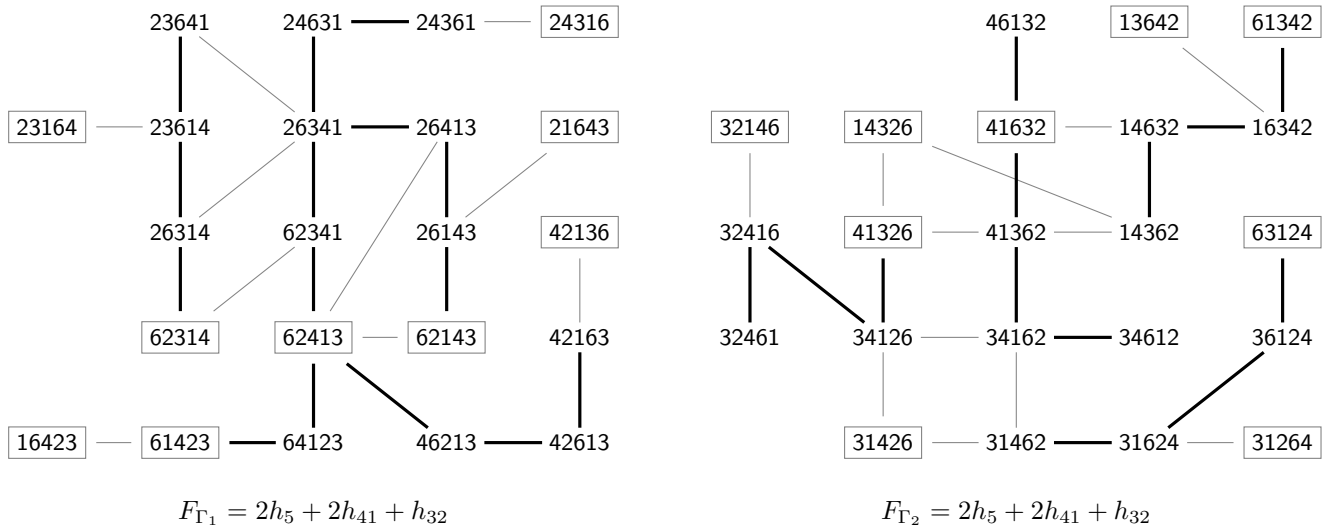

\begin{example}
\label{ex H disconnected}
Let $P =\PP_{3}$,  $\beta = \{1,2,3,4,5\}$, and let $\Gamma^\beta_d$ denote the $H$-graph  on
the words $\{\e{w} :  \text{\rm content}(\e{w}) = \beta, \, \inv_P(\e{w}) = d \}$.
The following table gives the symmetric functions of the components of $\Gamma^\beta_d$ for each possible  $d$.
\[
\begin{array}{c|c|c}
\inv =d& \text{\small \# of components of $\Gamma^\beta_d$} & \parbox{6.5cm}{\small the symmetric functions  $F_{\hat{\Gamma}}(\mathbf{x})$ over \\ the connected components  $\hat{\Gamma}$ of $\Gamma^\beta_d$}\\[4mm]
\hline  \\[-5mm]
& &   \\[-4mm]
0& 1&  h_{5} \\
1& 4 &   h_{5} , \ \, h_{5} , \ \, h_{5} , \ \, h_{5}  \\
2& 5 &  2 h_{5} +  h_{41} , \ \, 2 h_{5} +  h_{41} , \ \, h_{5} , \ \, h_{5} , \ \, h_{5}  \\
3& 4 & 5 h_{5} + 3 h_{41} +  h_{32} , \ \, h_{5} +  h_{41} , \ \, h_{5} , \ \, h_{5}  \\
4 & 4 & 5 h_{5} + 3 h_{41} +  h_{32} , \ \, h_{5} +  h_{41} , \ \, h_{5} , \ \, h_{5} \\
5 & 5 &   2 h_{5} +  h_{41} , \ \, 2 h_{5} +  h_{41} , \ \, h_{5} , \ \, h_{5} , \ \, h_{5}  \\
6 & 4 & h_{5} , \ \, h_{5} , \ \, h_{5} , \ \, h_{5}  \\
 7& 1 & h_{5}
\end{array}\]
\end{example}

Just as for  $P$-Knuth equivalence graphs,
the $0,1$-vectors in  $(I_H^P)^\perp$ which are supported on words of a fixed length are exactly the vectors $\sum_{\e{w} \in  \ver(\Gamma)} \e{w}$ for  $H$-graphs  $\Gamma$.

For an  $H$-graph~$\Gamma$, we set
$F_\Gamma(\mathbf{x})=F_\gamma(\mathbf{x})$ where $\gamma=\sum_{\e{w} \in  \ver(\Gamma)} \e{w}$.

Recall the Stanley-Stembridge and Shareshian-Wachs $e$-positivity conjectures (Conjectures \ref{cj SS intro} and  \ref{cj SW intro}).
We make the following conjecture.

\begin{conjecture}
\label{cj:m pos}
For any \threeone-free poset  $P$,
the elements  $\mathfrak{m}^P_\lambda(\mathbf{u})$ defined in \eqref{eq:P m def} are  $\mathbf{u}$-monomial positive modulo  $I_H^P$.  Hence, by Corollary \ref{c positivity} (ii),
for any  $I_H^P$-graph $\Gamma$, the symmetric function $\omega\spa F_\Gamma(\mathbf{x})$ is $e$-positive.
\end{conjecture}

\begin{remark}
\label{r how generalize}
This is related to \cite[Conjecture 3.10]{Hwang} and
 Conjectures \ref{cj SS intro} and \ref{cj SW intro} as follows:
\begin{itemize}
\item[(a)]
This is the same as \cite[Conjecture 3.10]{Hwang} except that we
state it in the generality of \threeone-free posets rather than natural unit interval orders.
Though but a small addition,
we do offer some additional evidence that \threeone-free is
the right level of generality here---namely, Theorems \ref{th:IplacP-positivity}, \ref{u pos hook}, and \ref{t m two col} all support this.
\item[(b)] For natural unit interval orders  $P$, the elements  $\Wo{\beta}(t)$ lie in
$(I_T^P)^\perp \subset (\QQ(t) \tsr I_H^P)^\perp$
and hence Conjecture \ref{cj:m pos} (and \eqref{eq:F-via-m}) imply the $e$-positivity of $X_{\inc(P)}(\mathbf{x},t)$ and its multicolored versions $X_{\inc(P)}^{\beta}(\mathbf{x},t)$.
\item[(c)] Conjecture \ref{cj:m pos} is stronger than (b) since the  $H$-graph
 $\Gamma^{\beta}_d$ on the words $\big\{\e{w} :  \text{\rm content}(\e{w}) = \beta, \, \inv_P(\e{w}) = d\big\}$   can be  disconnected (see Figure \ref{fig:H disconnected} and Example~\ref{ex H disconnected}),
 in which case  the coefficient of  $t^d$ in $X_{\inc(P)}^{\beta}(\mathbf{x},t)$ is expressed as a positive sum of ``smaller''  $e$-positive things (namely, $\omega \spa F_{\hat{\Gamma}}(\mathbf{x})$ over the components  $\hat{\Gamma}$ of  $\Gamma^\beta_d$).
 \\[-3mm]
\item[(d)]
Conjecture \ref{cj SW intro} implies Conjecture \ref{cj SS intro}
but only indirectly through \cite{GPchromatic}.  This leaves the question of whether there is a way to define a  $t$-statistic on \threeone-free posets with an accompanying  $e$-positivity conjecture both encompassing Conjecture \ref{cj SW intro} and specializing to Conjecture \ref{cj SS intro} at  $t=1$.  There does not seem to be a natural way to do this,
but Conjecture \ref{cj:m pos} remedies this by offering a refinement similar to but finer
than the  $t$-statistic which makes sense for any \threeone-free poset.
\end{itemize}
\end{remark}

In partial progress towards Conjecture \ref{cj:m pos}, we now state theorems showing that $\mathfrak{m}^P_\lambda(\mathbf{u})$ is $\mathbf{u}$-monomial positive modulo $I_H^P$ when $\lambda$ is a hook or two-column  shape.
The proofs of these theorems are given in Section \ref{s two col hook}.
Hwang \cite{Hwang} gives formulas for these same cases, but
our combinatorial objects differ in interesting ways---we use variations of  $P$-tableaux, while Hwang uses objects called heaps.

We introduce the following definitions in order to state the hook shape result.
\begin{definition}
\label{def:power word}
For a  $P$-weakly increasing word $\e{w} = \e{w}_1\e{w}_2\cdots \e{w}_n \in \U_P^*$ (i.e.  $\Des_P(\e{w}) = \varnothing$), we say $\e{w}_i$ is a \emph{right-left $P$-minima} if $\e{w}_i \to_P \e{w}_j$ for each $j = i+1,i+2,\dots,n$.
Each  $P$-weakly increasing word has a trivial right-left $P$-minima at index $n$. We say $\e{w}$ is a \emph{power word} if $\e{w}$ is a  $P$-weakly increasing word with no nontrivial right-left $P$-minima.
\end{definition}

\begin{definition}
Let $T$ be a hook shape $P$-tableau with
entries $\e{w}_1  \dotto \e{w}_2 \dotto \cdots \dotto \e{w}_\ell$ in the first row.
We say $T$ is a \emph{key $P$-tableau} if
$\e{v} \spa \e{w}_2\e{w}_3\cdots \e{w}_\ell$ is a power word for some entry $\e{v}$ in the first column of  $T$.

For a hook shape $\lambda$, let $\KT_P(\lambda)$ denote the set of key  $P$-tableaux of shape $\lambda$, and let $\KT_P^\beta(\lambda)$ denote the set of key  $P$-tableaux of shape $\lambda$ and content  $\beta$.
\end{definition}

\begin{example}
For the poset $P = \PP_{3}$, 
the first $P$-tableau below is a key  $P$-tableau since  $\e{5467}$ is a power word, while the second is not since $\e{1}$ is a right-left  $P$-minima of $\e{1456}$ and the words $\e{7456}$ and $\e{10 \, 4 \, 5\, 6}$ have  $P$-descents.
\[\begin{tabular}{cc}
\\[-1mm]
{\fontsize{8pt}{6pt}\selectfont \tableau{
1 & 4 & 6 & 7\\
5 \\
8 }}
&
\qquad
{\fontsize{8pt}{6pt}\selectfont \tableau{
1 &  4 & 5& 6 \\
7 \\
10 }}
\\[6mm]
 \text{key $P$-tableau} & \qquad \text{not a key $P$-tableau}
\end{tabular}\]
\end{example}

\begin{theorem}
\label{u pos hook}
For a hook shape $\lambda$ and any \threeone-free poset  $P$,
$\mathfrak{m}^P_\lambda(\mathbf{u})$ has the following  monomial positive expression modulo  $I_H^P$\spa:
\begin{align}
\mathfrak{m}^P_\lambda(\mathbf{u}) \equiv \sum_{T \in \KT_P(\lambda)} \mathbf{u}_{\creading(T)}.
\end{align}
\end{theorem}

\begin{corollary}
\label{c epos hook}
Let  $P$ be a \threeone-free poset and  $\lambda$ a hook shape.
For any  $I_H^P$-graph  $\Gamma$, the coefficient of  $e_\lambda(\mathbf{x})$ in the  $e$-expansion of
$\omega F_\Gamma(\mathbf{x})$ is equal to the number of
key $P$-tableaux $T\in \KT_P(\lambda)$ such that
the column reading word $\creading(T)$ appears as a vertex in~$\Gamma$.

In particular, for any  multiset  $\beta$ of elements of  $P$,
\begin{align}
\label{eq: epos hook}
\big( \text{coef of  $e_\lambda(\mathbf{x})$ in the $e$-expansion of
$X_{\inc(P)}^{\beta}(\mathbf{x})$}\big)
=
\big| \KT_P^\beta(\lambda)\big|.
\end{align}
Moreover, if  $P$ is a natural unit interval order,
\begin{align}
\label{eq: epos hook 2}
\big( \text{coef of  $e_\lambda(\mathbf{x})$ in the $e$-expansion of
$X_{\inc(P)}^{\beta}(\mathbf{x},t)$}\big)
=
\sum_{T \in \KT_P^\beta(\lambda)} t^{\inv_P(\creading(T))}.
\end{align}
\end{corollary}
\begin{proof}
Combine Theorem \ref{u pos hook} and
\eqref{eq:F-via-m}
with  $\gamma$ equal to  $\sum_{\e{w} \in \ver(\Gamma)} \e{w}$ or $\Wo{\beta}$ or  $\Wo{\beta}(t)$.  
We also need $\Wo{\beta}(t) \in (\QQ(t) \tsr I_H^P)^\perp$, which follows from Propositions \ref{p IH inclusions} and
\ref{pr:Wo beta t}.
\end{proof}

A formula for the
hook shape case of the Stanley-Stembridge conjecture (i.e. the coefficient of  $e_\lambda$ in  $X_{\inc(P)}$ for a hook shape  $\lambda$) was given by Wolfgang \cite{Wolfgangthesis}.
Hwang \cite{Hwang} gives a result similar to Corollary \ref{c epos hook}, with different combinatorial objects.

\begin{example}
For  $P = \PP_{2}$, $\lambda = (311)$, and  $\beta = \{1,3,4,5,7\}$, the key  $P$-tableaux of shape  $\lambda$
and content  $\beta$ are
\begin{align}
{\fontsize{8pt}{6pt}\selectfont \tableau{
1 & 4 & 5\\
3 \\
7 }}
\quad \quad
{\fontsize{8pt}{6pt}\selectfont \tableau{
1 & 5 & 4\\
3 \\
7 }}
\quad \quad
{\fontsize{8pt}{6pt}\selectfont \tableau{
1 & 4 & 3 \\
5 \\
7 }}.
\end{align}
Hence the coefficient of  $e_\lambda(\mathbf{x})$ in the  $e$-expansion of  $X_{\inc(P)}(\mathbf{x})$ is 3.
\end{example}

We now introduce definitions to state the two-column result.

\begin{definition}
\label{def:ladder}
Let  $(C,D)$ be a pair of (not necessarily disjoint) chains in a poset  $P$. Let  $\inc_P(C,D)$ denote the
bipartite graph with partite sets  $C, D$ and edge set consisting of the
pairs  $c \dote_P d$ of incomparable or equal elements, with $c \in C$,  $d \in D$.
A \emph{ladder $L$ of $(C,D)$} is a connected component of  $\inc_P(C,D)$;  we often identify  $L$ with its vertex set, a subset of the disjoint union $C \sqcup D$.
We say $L$ is \emph{balanced} if $|L\cap C| = |L\cap D|$;
we say  $L$ is a  \emph{$C$-ladder} or is \emph{left unbalanced} if $|L\cap C| > |L\cap D|$,
 and  $L$ is a \emph{$D$-ladder} or \emph{right unbalanced} if $|L\cap C| < |L\cap D|$.

For a  $P$-tableau  $T$ with at most two columns, define a \emph{ladder of  $T$} to be  a ladder of  $(T_1,T_2)$,
where  $T_i$ denotes the  $i$-th column of  $T$.
\end{definition}

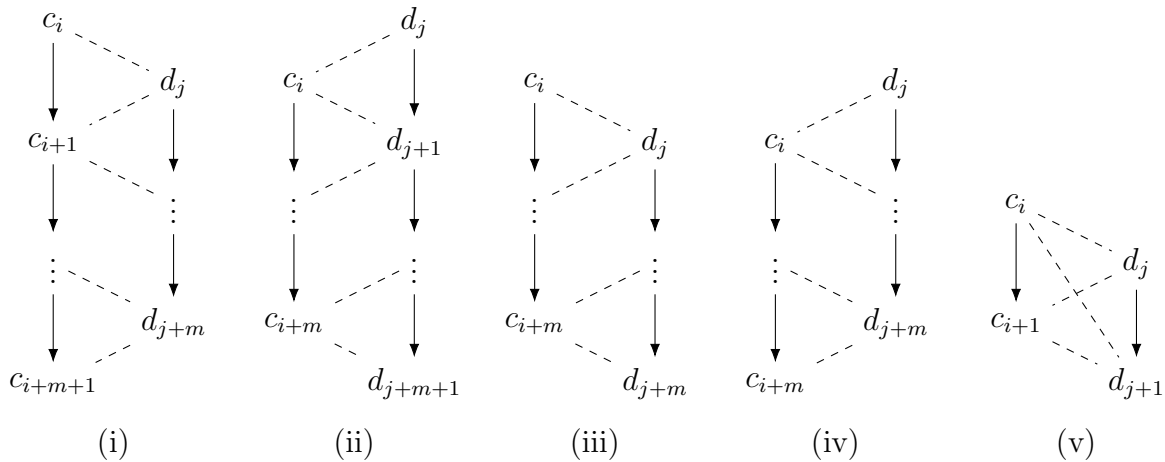
\begin{figure}[h]
\begin{tikzpicture}[scale = 0.8]
\node at (1,-1) {(i)};

\node (11) at (0,6) {$c_{i}$};
\node (12) at (0,4) {$c_{i+1}$};
\node (13) at (0,2) {$\vdots$};
\node (14) at (0,0) {$c_{i+m+1}$};

\node (21) at (2,5) {$d_{j}$};
\node (22) at (2,3) {$\vdots$};
\node (23) at (2,1) {$d_{j+m}$};

\draw[mystealth] (13) to (14);
\draw[mystealth] (11) to (12);
\draw[mystealth] (12) to (13);
\draw[mystealth] (21) to (22);
\draw[mystealth] (22) to (23);

\draw[dashed] (11) to (21);
\draw[dashed] (21) to (12);
\draw[dashed] (12) to (22);
\draw[dashed] (13) to (23);
\draw[dashed] (23) to (14);

\node at (5,-1) {(ii)};

\node (41) at (6,6) {$d_{j}$};
\node (42) at (6,4) {$d_{j+1}$};
\node (43) at (6,2) {$\vdots$};
\node (44) at (6,0) {$d_{j+m+1}$};

\node (31) at (4,5) {$c_{i}$};
\node (32) at (4,3) {$\vdots$};
\node (33) at (4,1) {$c_{i+m}$};

\draw[mystealth] (43) to (44);
\draw[mystealth] (41) to (42);
\draw[mystealth] (42) to (43);
\draw[mystealth] (31) to (32);
\draw[mystealth] (32) to (33);

\draw[dashed] (41) to (31);
\draw[dashed] (31) to (42);
\draw[dashed] (42) to (32);
\draw[dashed] (43) to (33);
\draw[dashed] (33) to (44);

\node at (9,-1) {(iii)};

\node (51) at (8,5) {$c_i$};
\node (52) at (8,3) {$\vdots$};
\node (53) at (8,1) {$c_{i+m}$};

\node (61) at (10,4) {$d_j$};
\node (62) at (10,2) {$\vdots$};
\node (63) at (10,0) {$d_{j+m}$};

\draw[mystealth] (51) to (52);
\draw[mystealth] (52) to (53);
\draw[mystealth] (61) to (62);
\draw[mystealth] (62) to (63);

\draw[dashed] (51) to (61);
\draw[dashed] (61) to (52);

\draw[dashed] (62) to (53);
\draw[dashed] (53) to (63);

\node at (13,-1) {(iv)};

\node (71) at (12,4) {$c_i$};
\node (72) at (12,2) {$\vdots$};
\node (73) at (12,0) {$c_{i+m}$};

\node (81) at (14,5) {$d_j$};
\node (82) at (14,3) {$\vdots$};
\node (83) at (14,1) {$d_{j+m}$};

\draw[mystealth] (81) to (82);
\draw[mystealth] (82) to (83);
\draw[mystealth] (71) to (72);
\draw[mystealth] (72) to (73);

\draw[dashed] (81) to (71);
\draw[dashed] (71) to (82);
\draw[dashed] (72) to (83);
\draw[dashed] (83) to (73);

\node at (17,-1) {(v)};

\node (91) at (16,3) {$c_i$};
\node (92) at (16,1) {$c_{i+1}$};

\node (101) at (18,2) {$d_j$};
\node (102) at (18,0) {$d_{j+1}$};

\draw[mystealth] (91) to (92);
\draw[mystealth] (101) to (102);

\draw[dashed] (91) to (101);
\draw[dashed] (92) to (102);
\draw[dashed] (91) to (102);
\draw[dashed] (101) to (92);
\end{tikzpicture}
\caption{\label{fig:ladder}
For a \threeone-free poset  $P$, the five possible forms of a ladder of a pair of chains $(C,D)$ of  $P$, where
$C = c_k \leftto_P \cdots \leftto_P c_1$,
$D = d_k \leftto_P \cdots \leftto_P d_1$.}
\end{figure}

A ladder in a \threeone-free poset has a very controlled form---see Figure \ref{fig:ladder} and the justification in Lemma \ref{l ladder}.

\begin{definition}
A $P$-tableau $T$ with at most two columns is a \emph{left  $P$-tableau} if every ladder in $T$ is either balanced or left unbalanced. For a partition $\lambda$ with  $\lambda_1 \le 2$, let $\LT_P(\lambda)$ denote the set of left  $P$-tableaux of shape $\lambda$.
\end{definition}

\begin{example}
For the poset $P = \PP_{2}$, 
the first $P$-tableau below is a left  $P$-tableau as its ladders are $\{1,2,3\}$ and  $\{5,6,7\}$ which are both left unbalanced.  The second  $P$-tableau has the same ladders but now $\{5,6,7\}$ is right unbalanced, so it is not a left  $P$-tableau.
\[\begin{tabular}{cc}
\\[-4mm]
{\fontsize{8pt}{6pt}\selectfont \tableau{
1 & 2 \\
3 & 6 \\
5 \\
7}}
&
\qquad
{\fontsize{8pt}{6pt}\selectfont \tableau{
1 & 2 \\
3 & 5 \\
6 & 7}}
\\[8mm]
 \text{a left  $P$-tableau} & \qquad \text{not a left  $P$-tableau}
\end{tabular}\]
\end{example}

\begin{theorem}
\label{t m two col}
For a partition  $\lambda$ with  $\lambda_1 \le 2$ and any \threeone-free poset  $P$,
$\mathfrak{m}^P_\lambda(\mathbf{u})$ has the following  monomial positive expression modulo  $I_H^P$\spa:
\begin{align}
\label{eq: t m two col}
\mathfrak{m}^P_\lambda(\mathbf{u}) \equiv \sum_{T \in \LT_P(\lambda)} \mathbf{u}_{\creading(T)}.
\end{align}
\end{theorem}

Just as Theorem \ref{u pos hook} gave Corollary \ref{c epos hook}, Theorem \ref{t m two col} yields
\begin{corollary}
\label{c epos twocol}
Let  $P$ be a \threeone-free poset and  $\lambda$ be a partition with  $\lambda_1 \le 2$.
For any   $I_H^P$-graph  $\Gamma$, the coefficient of  $e_\lambda(\mathbf{x})$ in the  $e$-expansion of $\omega F_\Gamma(\mathbf{x})$ is equal to the number of
left $P$-tableaux $T\in \LT_P(\lambda)$ such that
 $\creading(T) \in \ver(\Gamma)$.

Moreover, if  $P$ is a natural unit interval order and $\beta$ is a multiset of elements of $P$,
\begin{align}
\Big( \text{coef of  $e_\lambda(\mathbf{x})$ in the $e$-expansion of
$X_{\inc(P)}^{\beta}(\mathbf{x},t)$} \Big)
=
\sum_{\substack{T \in \LT_P(\lambda) \\ \text{$T$ has content $\beta$}}} t^{\inv_P(\creading(T))}.
\end{align}
\end{corollary}

For two column  $\lambda$,
a formula for the coefficient of  $e_\lambda$ in  $X_{\inc(P)}(\mathbf{x})$ was first given
in \cite{Wolfgangthesis}, and
a formula for the coefficient
of  $e_\lambda$ in  $X_{\inc(P)}(\mathbf{x},t)$  was
first given in \cite{CHSS}.
Hwang \cite{Hwang} gives a result similar to Corollary \ref{c epos twocol}, with different combinatorial objects.

\section{The arrow algebra}
\label{sec:poset algebra}


We define and study the arrow algebra  $\A$ and several quotients including
$\A/I_\text{closed}$.  The algebra $\A/I_\text{closed}$ is essentially the set of formal linear combinations of  labeled posets
equipped with a natural concatenation product.  We show that $\A/I_\text{closed}$ injects into the product of  $\U_P$'s over all posets  $P$, via a map given in terms of a natural generalization of  $P$-tableaux called $P$-fillings.
We use this to reformulate our main positivity statement and conjecture from
\S\ref{ss Schur pos P Knuth}, \S\ref{ss e pos} in terms of the arrow algebra (\S\ref{ss strengthening SS}).
This suggests ways to further strengthen this conjecture (Question \ref{q: m diagram}) and leads to variations
on Stanley's $P$-partitions conjecture and related questions of McNamara (\S\ref{ss Stanley variations}).


\subsection{The arrow algebra}


The \textit{arrow space on $d$ vertices} is $\A_d\defeq \bigotimes\limits_{1\leq i<j\leq d}V_{ij}$ where each $V_{i,j}$ is a copy of the free $\ZZ$-module $V=\ZZ\{r,\ell,e,n\}\cong \ZZ^4$, whose elements are called \emph{edge types}.
An \textit{arrow diagram} $M\in \A_d$ is a simple tensor $M=\otimes_{1\leq i<j\leq d} \, M_{i,j}$ where each $M_{i,j} \in V_{i,j}$ (called an \emph{edge} of  $M$) is one of the following eight edge types: $\{r, \ell, e, n, e+n, r+e+n, \ell+e+n, r+\ell+e+n\}$.
We draw arrow diagrams with vertices $1,\ldots,d$ in a row, increasing from left to right, and edge types as shown below. 
The meaning of the edge types and why we have given them an additive structure will become apparent below.
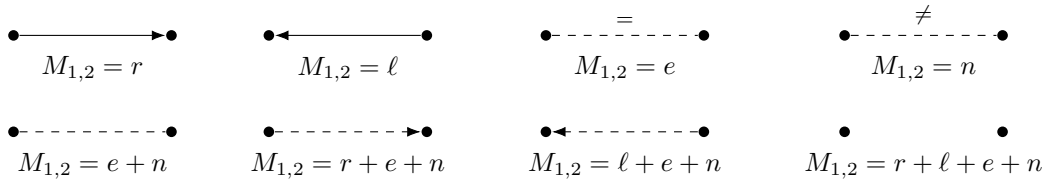
\[\begin{array}{c@{\quad \ \ \ }c@{\quad \ \ \ }c@{\quad \ \ \ }c}
\begin{tikzpicture}[xscale=0.7, yscale=0.9]
\vertex (11) at (0,6) {};
\vertex (12) at (3,6) {};
\draw[mystealth] (11) -- (12);
\node at (1.5,5.5) {\footnotesize $M_{1,2}=r$};
\end{tikzpicture}
&
\begin{tikzpicture}[xscale=0.7, yscale=0.9]
\vertex (11) at (0,6) {};
\vertex (12) at (3,6) {};
\draw[mystealth] (12) -- (11) ;
\node at (1.5,5.5) {\footnotesize $M_{1,2}=\ell$};
\end{tikzpicture}
&\begin{tikzpicture}[xscale=0.7, yscale=0.9]
\vertex (11) at (0,6) {};
\vertex (12) at (3,6) {};
\draw[dashed] (11) -- (12) node [midway, above] {\tiny $=$};
\node at (1.5,5.5) {\footnotesize $M_{1,2}=e$};
\end{tikzpicture}
&
\begin{tikzpicture}[xscale=0.7, yscale=0.9]
\vertex (11) at (0,6) {};
\vertex (12) at (3,6) {};
\draw[dashed] (11) -- (12) node [midway, above] {\tiny $\ne$};
\node at (1.5,5.5) {\footnotesize $M_{1,2}=n$};
\end{tikzpicture}
\\[3mm]
\begin{tikzpicture}[xscale=0.7, yscale=0.9]
\vertex (11) at (0,6) {};
\vertex (12) at (3,6) {};
\draw[dashed] (11) -- (12);
\node at (1.5,5.5) {\footnotesize $M_{1,2}=e+n$};
\end{tikzpicture}
&
\begin{tikzpicture}[xscale=0.7, yscale=0.9]
\vertex (11) at (0,6) {};
\vertex (12) at (3,6) {};
\draw[mystealth][dashed] (11) -- (12);
\node at (1.5,5.5) {\footnotesize $M_{1,2}=r+e+n$};
\end{tikzpicture}
&
\begin{tikzpicture}[xscale=0.7, yscale=0.9]
\vertex (11) at (0,6) {};
\vertex (12) at (3,6) {};
\draw[mystealth][dashed] (12) -- (11);
\node at (1.5,5.5) {\footnotesize $M_{1,2}=\ell+e+n$};
\end{tikzpicture}
&
\begin{tikzpicture}[xscale=0.7, yscale=0.9]
\vertex (11) at (0,6) {};
\vertex (12) at (3,6) {};
\node at (1.5,5.5) {\footnotesize $M_{1,2}=r+\ell+e+n$};
\end{tikzpicture}
\end{array}\]
An arrow diagram $M\in \A_d$ is \emph{primitive} if every $M_{i,j}\in\{r,\ell,e,n\}$. The primitive arrow diagrams on $d$ vertices form a $\ZZ$-basis for $\A_d$.

\begin{example}
On the left is the arrow diagram  $M \in \A_3$ with  $M_{1,2} = e+ n$,  $M_{2,3} = r$, and  $M_{1,3}= e+n$, and on the right is its expansion into primitive arrow diagrams.
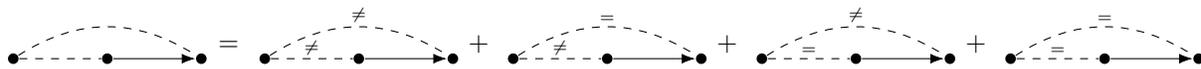
\[{\footnotesize \text{$
\raisebox{-1.7mm}{
\begin{tikzpicture}[xscale=1.25, yscale=1.38]
\tikzstyle{aedgecurve} = [dashed, draw, black, bend left=30]
\vertex (a) at (1,0) { };
\vertex (b) at (2,0) { };
\vertex (c) at (3,0) { };
\draw[dashed] (a) to  (b);
\draw[mystealth]  (b) to (c);
\draw[aedgecurve]  (a) to (c);
\end{tikzpicture}}
\,  = \,
\raisebox{-1.7mm}{
\begin{tikzpicture}[xscale=1.25, yscale=1.38]
\tikzstyle{aedgecurve} = [draw, black, dashed, bend left=30]
\vertex (a) at (1,0) { };
\vertex (b) at (2,0) { };
\vertex (c) at (3,0) { };
\draw[dashed] (a) to node[above,inner sep = 1.5pt] {\Tiny$\ne$} (b);
\draw[mystealth]  (b) to (c);
\draw[aedgecurve]  (a) to node[above, inner sep = 1.5pt] {\Tiny$\ne$} (c);
\end{tikzpicture}}
\,  + \,
\raisebox{-1.7mm}{
\begin{tikzpicture}[xscale=1.25, yscale=1.38]
\tikzstyle{aedgecurve} = [draw, black, dashed, bend left=30]
\vertex (a) at (1,0) { };
\vertex (b) at (2,0) { };
\vertex (c) at (3,0) { };
\draw[dashed] (a) to node[above,inner sep = 1.5pt] {\Tiny$\ne$} (b);
\draw[mystealth]  (b) to (c);
\draw[aedgecurve]  (a) to node[above,inner sep = 1.8pt] {\Tiny$=$} (c);
\end{tikzpicture}}
\,  + \,
\raisebox{-1.7mm}{
\begin{tikzpicture}[xscale=1.25, yscale=1.38]
\tikzstyle{aedgecurve} = [draw, black, dashed, bend left=30]
\vertex (a) at (1,0) { };
\vertex (b) at (2,0) { };
\vertex (c) at (3,0) { };
\draw[dashed] (a) to node[above,inner sep = 1.8pt] {\Tiny$=$} (b);
\draw[mystealth]  (b) to (c);
\draw[aedgecurve]  (a) to node[above, inner sep = 1.5pt] {\Tiny$\ne$} (c);
\end{tikzpicture}}
\,  + \,
\raisebox{-1.7mm}{
\begin{tikzpicture}[xscale=1.25, yscale=1.38]
\tikzstyle{aedgecurve} = [draw, black, dashed, bend left=30]
\vertex (a) at (1,0) { };
\vertex (b) at (2,0) { };
\vertex (c) at (3,0) { };
\draw[dashed] (a) to node[above,inner sep = 1.8pt] {\Tiny$=$} (b);
\draw[mystealth]  (b) to (c);
\draw[aedgecurve]  (a) to node[above,inner sep = 1.8pt] {\Tiny$=$} (c);
\end{tikzpicture}}
$} }
\]
\end{example}

The \textit{arrow algebra} $\A$ is the graded ring $\bigoplus_{d\geq0}\A_d$ (by convention $\A_1=\A_0=\ZZ$)
equipped with a kind of concatenation product, defined as follows.  Its homogeneous component $*: \A_d\otimes \A_{d'}\rightarrow \A_{d+d'}$ is the  $\ZZ$-linear map such that 
for simple tensors $L=\otimes_{1\leq i<j\leq d} \, L_{i,j}$ and $M=\otimes_{1\leq i<j\leq d'} \, M_{i,j}$, the product  $L*M$ is the simple tensor $\otimes_{1\leq i<j\leq d+d'} \, (L*M)_{i,j}$ with
\begin{itemize}
\item $(L*M)_{i,j}=L_{i,j}$ \ when $i,j\in[d]$,
\item $(L*M)_{i,j}=M_{i-d,j-d}$ \ when $i,j\in\{d+1,\ldots,d+d'\}$,
\item $(L*M)_{i,j}=r+\ell+e+n$ \ when $i\in[d]$ and $j\in\{d+1,\ldots,d+d'\}.$
\end{itemize}
The drawing of  $L*M$ is obtained by placing  $M$ to the right of  $L$. For example,
%
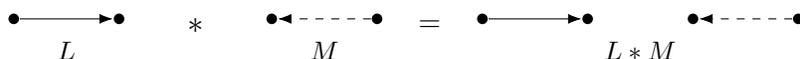
\begin{center}
\begin{tikzpicture}[xscale=0.7, yscale=0.85]

\vertex (11) at (0-0.9,0) {};
\vertex (12) at (2-0.9,0) {};
\node at (3-0.9*0.5,-0.1) {$*$};
\vertex (21) at (4,0) {};
\vertex (22) at (6,0) {};
\node at (7,-0.1) {\small $=$};
\vertex (31) at (8,0) {};
\vertex (32) at (10,0) {};
\vertex (33) at (12,0) {};
\vertex (34) at (14,0) {};
\node at (1-0.9,-0.5) {\footnotesize $L$};
\node at (5,-0.5) {\footnotesize $M$};
\node at (11,-0.5){\footnotesize $L*M$};
\draw[mystealth][dashed] (22) -- (21);
\draw[mystealth] (11) -- (12);
\draw[mystealth][dashed] (34) -- (33);
\draw[mystealth] (31) -- (32);
\end{tikzpicture}
\end{center}
This element  $L*M \in \A_4$ is equal to the sum of  $3\cdot4^4$ primitive arrow diagrams.

The notion of reversing an edge type is formalized as follows: define $\text{\rm rev}:V\rightarrow V$
to be the automorphism which fixes $e$ and $n$ and satisfies $\rev(r)=\ell$ and $\rev(\ell)=r$. For a simple tensor $M=\otimes_{1\leq i<j\leq d} \, M_{i,j}$, we define $M_{j,i}=\rev(M_{i,j})$.

\subsection{$P$-fillings of an arrow diagram}
Given a  poset  $P$ and an arrow diagram $M\in \A_d$, a \textit{$P$-filling of $M$} is a word  $\e{w}= \e{w}_1 \cdots \e{w}_d \in \U_P^*$ satisfying (for all distinct $i,j\in[d]$)
\begin{itemize}
\item $\e{w}_i\to_P\e{w}_j$ when $M_{i,j}=r$,
\item $\e{w}_i=\e{w}_j$ when $M_{i,j}=e$,
\item $\e{w}_i \dote_P \e{w}_j$ and $\e{w}_i\not=\e{w}_j$ when $M_{i,j}=n$,
\item $\e{w}_i \dote_P \e{w}_j$ when $M_{i,j}=e+n$,
\item $\e{w}_i \dotto_P \e{w}_j$ when $M_{i,j}=r+e+n$,
\item No additional condition on $\e{w}$ when $M_{i,j}=r+\ell+e+n$.
\end{itemize}
Let $\text{\rm Fill}_P(M)$ denote the set of all $P$-fillings of $M$.

Let $\text{\rm Eval}_P:\A\rightarrow \U_P$ denote the $\ZZ$-linear map determined by
\begin{align}
\label{eq:def Eval}
\text{\rm Eval}_P(M)=\sum\limits_{\e{w} \in \text{\rm Fill}_P(M)} \mathbf{u}_{\e{w}}
\end{align}
whenever $M\in \A_d$ is a primitive arrow diagram.

\begin{example}
\label{ex P fillings}
For $P$ the poset  \twotwo, with incomparable chains  $a \to_P b$ and  $c \to_P d$,
\begin{align}
\label{eq:ex P fillings}
\text{\rm Eval}_P \Big(
\raisebox{-3mm}{
\begin{tikzpicture}[scale=1.7]
\tikzstyle{aedgecurve} = [draw, black, dashed, bend left=30]
\vertex (a) at (1,0) { };
\vertex (b) at (2,0) { };
\vertex (c) at (3,0) { };
\vertex (d) at (4,0) { };
\draw[mystealth, bend left = 47]  (a) to (d);
\draw[mystealth]  (c) to (b);
\draw[dashed] (a) to node[above,inner sep = 1.5pt] {\Tiny$\ne$} (b);
\draw[dashed]  (c) to node[above,inner sep = 1.5pt] {\Tiny$\ne$} (d);
\draw[aedgecurve]  (a) to node[above,inner sep = 1.5pt] {\Tiny$\ne$} (c);
\draw[aedgecurve]  (b) to node[above,inner sep = 1.5pt] {\Tiny$\ne$} (d);
\end{tikzpicture}} \,
\Big)  = u_au_du_c u_b + u_cu_bu_a u_d.
\end{align}
\end{example}
In \S\ref{ss strengthening SS}, we will see that $P$-fillings and  $P$-tableaux are related in a natural way.

\begin{proposition}
(i) Formula \eqref{eq:def Eval} holds for any arrow diagram $M\in \A_d$.

(ii) $\text{\rm Eval}_P(L*M)=\text{\rm Eval}_P(L)\text{\rm Eval}_P(M)$ for any primitive arrow diagrams $L$ and $M$.

(iii) $\text{\rm Eval}_P$ is a ring homomorphism.
\end{proposition}
\begin{proof}
Statement (i) is immediate from the definition of $P$-filling, comparing $P$-fillings of $M$ to the $P$-fillings of every primitive arrow diagram in the expansion of $M$ into primitive arrow diagrams. Then, (ii) is a consequence of (i), noting that $L*M$ is an arrow diagram whose fillings are in straightforward bijection with $\text{\rm Fill}_P(L)\times \text{\rm Fill}_P(M)$.  Statement (iii) then follows from (ii) by linearity.
\end{proof}

\begin{remark}
\label{r P partition}
Stanley's  $(P,\omega)$-partitions \cite{StanleyPpartition} and a generalization studied by McNamara can be viewed as  $P$-fillings.
McNamara \cite[\S2]{McNamara} defines an \emph{oriented poset}  $(D,O)$ to be a poset  $D$ together with
an assignment  $O$ of the covering relations of  $D$ as being either \emph{strict} or \emph{weak}.
He then defines a $(D,O)$-\emph{partition} to be a map  $\kappa \colon D \to [N]$ satisfying
\begin{itemize}
\item If  $a \to_D b$, then  $\kappa(a) \le \kappa(b)$,
\item If  $a \to_D b$ is a strict covering relation, then   $\kappa(a) < \kappa(b)$.
\end{itemize}
(We have slightly modified McNamara's definition to $(D,O)$-partitions taking values in  $[N]$ rather than  $\ZZ_{\ge 1}$.)

An oriented poset  $(D,O)$ is naturally encoded by an arrow diagram.
Assume that  $D$ is finite of size  $d$.
Label the elements of $D$ arbitrarily by  $[d]$. Let $M \in \A_d$ be the arrow diagram
with an edge $M_{i,j} = r$ for each strict covering relation $i \to_D j$, and $M_{i,j} = r+e+n$ for each weak covering relation  $i \to_D j$.
Let  $\PP$ be the total order  $[N]$.
Then $\kappa$ is a $(D,O)$-partition if and only if $\kappa(1)\cdots \kappa(d)$ is a  $\PP$-filling of  $M$.
Hence also, the image of $\text{\rm Eval}_{\PP}(M)$ in  $\U_{\PP}/\Ipol^{\PP} \cong \ZZ[y_1,\dots,y_N]$ is the same as the \emph{generating function} $K_{D,O}(y_1,\dots, y_N)$ of  $(D,O)$ in the sense of \cite[\S2]{McNamara}.
(Note that the $\PP$-fillings of  $M$ are words where the order of the letters depends on our arbitrary choice of labeling of $D$ with  $[d]$ but the image of $\text{\rm Eval}_{\PP}(M)$ in  $\U_{\PP}/\Ipol^{\PP}$ does not.)
\end{remark}

\subsection{The plactic arrow algebra}

We wish to only consider primitive arrow diagrams which naturally correspond to posets, so we quotient by those that do not.
Accordingly, let $I_\text{closed}$ be the ideal of $\A$ consisting of the linear span of the primitive arrow diagrams $M$ which have three distinct indices $i,j,k\in[d]$ satisfying either:
\begin{enumerate}[label=(\Roman*)]
\item $M_{i,j}=M_{j,k}=r$ but $M_{i,k}\not=r$,
\item $M_{i,j}=e$, but $M_{i,k}\not=M_{j,k}$.
\end{enumerate}
These conditions are the obvious ones which ensure a primitive arrow diagram has no  $P$-fillings for any  $P$
(see the proof of Proposition \ref{p I closed} below).
%
%

\begin{definition}
\label{def:clique poset}
For a primitive arrow diagram $M\in \A_d$ which is not in $I_\text{closed}$, condition (II) implies that the vertex set $[d]$ is necessarily partitioned into cliques by the edges with edge type $e$. Then, the \textit{clique poset of $M$} is the (unlabeled) poset $Q_M$ whose elements are in bijection with these cliques
 (let $\text{\rm cliq}_M(a)$ be the clique corresponding to $a\in Q_M$), with the poset relations satisfying $a\to_{Q_M}b$ if any only if some $i\in\text{\rm cliq}_M(a)$ and $j\in\text{\rm cliq}_M(b)$ satisfy $M_{i,j}=r$. In other words, $a\to_{Q_M}b$ whenever there is a solid arrow from $i$ to $j$ in $M$. Condition (I) ensures that the poset relations in $Q_M$ are transitive.
\end{definition}

By ``unlabeled'' in Definition \ref{def:clique poset}, we mean that we only want the isomorphism class of the poset $Q_M$ and we forget the cliques used to define it.  So, for example,
\begin{align*}
\text{clique poset} \Big(
\begin{tikzpicture}[scale=1.4]
\tikzstyle{aedgecurve} = [draw, black, dashed, bend left=30]
\vertex (a) at (1,0) { };
\vertex (b) at (2,0) { };
\vertex (c) at (3,0) { };
\draw[dashed] (a) to node[above,inner sep = 1.5pt] {\Tiny$\ne$} (b);
\draw[mystealth]  (b) to (c);
\draw[aedgecurve]  (a) to node[above,inner sep = 1.5pt] {\Tiny$\ne$} (c);
\end{tikzpicture}
\Big) =
\text{clique poset} \Big(
\begin{tikzpicture}[scale=1.4]
\tikzstyle{aedgecurve} = [draw, black, dashed, bend right=30]
\vertex (a) at (3,0) {};
\vertex (b) at (1,0) {};
\vertex (c) at (2,0) {};
\draw[aedgecurve] (a) to node[above,inner sep = 1.5pt] {\Tiny$\ne$} (b);
\draw[mystealth]  (b) to (c);
\draw[dashed]  (a) to node[above,inner sep = 1.5pt] {\Tiny$\ne$} (c);
\end{tikzpicture}
\Big)
=\twoone,
\end{align*}
where $\twoone$ denotes the poset which is the disjoint union of a 2 element chain and a 1 element chain.

\begin{proposition}
\label{p I closed}
Let $M$ be a primitive arrow diagram. The following are equivalent:
\begin{enumerate}[label=(\roman*)]
\item $M\in I_\text{closed}$,
\item $\text{\rm Fill}_P(M)=\varnothing$ for all (finite) posets $P$.
\end{enumerate}
\end{proposition}
\begin{proof}
Let $M\in I_\text{closed}$ be a primitive arrow diagram. Let $P$ be a poset. Any $P$-filling of $M$ would imply a failure of transitivity in $P$ or a single element of $P$ which has two distinct relationships with another single element of $P$. Therefore $\text{\rm Fill}_P(M)=\varnothing$.

Now, let $M\in \A_d$ be a primitive arrow diagram which is not in $I_\text{closed}$. Let $Q$ be the clique poset of
$M$. There exists a $Q$-filling  $\e{w}$ of $M$ with  $\e{w}_i = a$ for each  $a \in Q$ and $i \in \text{\rm cliq}_M(a)$.
\end{proof}

Let $M\not\in I_\text{closed}$ be a primitive arrow diagram. We say that $M$ is \threeone-free if the clique poset of $M$ is \threeone-free. Let $I_\text{\rm (3+1)}$ denote the linear span of all primitive arrow diagrams which are either in $I_\text{closed}$ or whose clique poset is not \threeone-free, so that $I_\text{closed}\subset I_\text{\rm (3+1)}\subset \A$, and $I_\text{\rm (3+1)}$ is an ideal in $\A$.
\begin{proposition}
\label{p 31free primitive}
Let $M$ be a primitive arrow diagram. The following are equivalent:
\begin{enumerate}[label=(\roman*)]
\item $M\in I_\text{\rm (3+1)}$,
\item $\text{\rm Fill}_P(M)=\varnothing$ for all \threeone-free posets $P$.
\end{enumerate}
\end{proposition}

We can apply local manipulations to arrow diagrams to describe relations analogous to those in
$I_{\text{\rm plac}}^P$, $I_H^P$, etc.
For an arrow diagram  $M \in \A_d$ and $i\in[d-1]$, define the arrow diagram $\text{\rm Swap}_i(M) \in \A_d$ by
\begin{itemize}
\item $\text{\rm Swap}_i(M)_{i,i+1}=\text{\rm rev}(M_{i,i+1})$,
\item $\text{\rm Swap}_i(M)_{i,k}=M_{i+1,k}$ \ \ for all $k\not\in\{i,i+1\}$,
\item $\text{\rm Swap}_i(M)_{i+1,k}=M_{i,k}$ \ \ for all $k\not\in\{i,i+1\}$.
\end{itemize}

\begin{center}
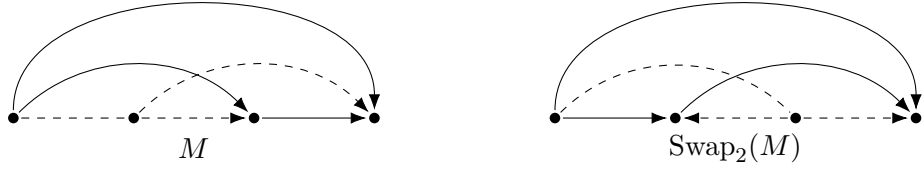
\begin{figure}[h]
\begin{tikzpicture}[scale=0.8]
\tikzstyle{vertex}=[circle, fill, inner sep=0pt, minimum size=4pt, outer sep = 1pt]

\node[vertex] (11) at (0,0) {};
\node[vertex] (12) at (2,0) {};
\node[vertex] (13) at (4,0) {};
\node[vertex] (14) at (6,0) {};

\node[vertex] (21) at (9,0) {};
\node[vertex] (22) at (11,0) {};
\node[vertex] (23) at (13,0) {};
\node[vertex] (24) at (15,0) {};

\draw[mystealthbig] (11) to [out=45, in=135] (13);
\draw[mystealthbig][dashed] (12) to (13);
\draw[dashed] (11) to (12);
\draw[mystealthbig] (13) to (14);
\draw[mystealthbig][dashed] (12) to [out=45, in=135] (14);
\draw[mystealthbig] (11) to [out=90, in=90] (14);

\draw[dashed] (21) to [out=45, in=135] (23);
\draw[mystealthbig][dashed] (23) to (22);
\draw[mystealthbig] (21) to (22);
\draw[mystealthbig] (21) to [out=90, in=90] (24);
\draw[mystealthbig] (22) to [out=45, in=135] (24);
\draw[mystealthbig][dashed] (23) to (24);

\node at (3,-0.5) {\small $M$};
\node at (12,-0.5) {\small $\text{\rm Swap}_2(M)$};
\end{tikzpicture}
\vspace{-1.5mm}
\caption{\label{fig:swap} $M$ and $\text{\rm Swap}_2(M)$. Edges incident to vertex $2$ have moved to instead be incident to vertex $3$, and vice versa. The edge between vertices $2$ and $3$ has been reversed.}
\end{figure}
\end{center}

\begin{definition}
\label{def:arrow Iplac}
The \emph{plactic arrow algebra} is the quotient $\A/\Iplac$, where
$\Iplac$ is the ideal of $\A$ consisting of the linear span of all arrow diagrams in $I_\text{\rm (3+1)}$ as well as the following elements:
\begin{enumerate}[leftmargin=1.54cm]
\item $M-\text{\rm Swap}_{i}(M)$ \qquad \quad \  \ \ \ \ \, for $M_{i-1,i+1}=M_{i,i+1}=\ell$ and $M_{i-1,i}\in\{e,n,r\}$,
\item $M-\text{\rm Swap}_{i-1}(M)$ \qquad  \quad  \ \ \spa \, for $M_{i-1,i}=M_{i-1,i+1}=\ell$ and $M_{i,i+1}\in\{e,n,r\}$,
\item $M-\text{\rm Swap}_{i}(\text{\rm Swap}_{i-1}(M))$ \ \, for $M_{i-1,i}=M_{i-1,i+1}=n$ and $M_{i,i+1}=\ell$,
\end{enumerate}
where  $M$ is any arrow diagram satisfying the stated conditions and  $1 < i < d$ for  $d$ the number of vertices of  $M$; see Figure \ref{fig:diagram Iplac}.
\end{definition}
It is straightforward to verify that the elements (1)--(3) map via $\text{\rm Eval}_P$ into the ideal $\Iplacp{P}$ of $\U_P$ for a \threeone-free poset $P$, i.e.,
\begin{align}
\label{eq: Eval Iplac}
\text{\rm Eval}_P(\Iplac) \subset \Iplacp{P}.
\end{align}

\newcommand{\mysize}{*0.84cm}
\begin{figure}[h]
\begin{center}
\begin{tikzpicture}[xscale = .9, yscale=.7]
\node (00) at (-0.75,-0.1) {$\cdots$};
\vertex (01) at (0,0) {};
\vertex (03) at (2,0) {};
\vertex (05) at (4,0) {};
\node (07) at (4.75, -0.1) {$\cdots$};

\node (08) at (6,0.25) {$\equiv$};
\node (09) at (6,-0.25) {(1a)};

\node (10) at (7.25,-0.1) {$\cdots$};
\vertex (11) at (8,0) {};
\vertex (13) at (10,0) {};
\vertex (15) at (12,0) {};
\node (17) at (12.75,-0.1) {$\cdots$};

\draw[mystealth] (05) to [out=135, in=45] (01);
\draw[dashed] (01) -- (03) node [midway, above = -0.2mm] {\Tiny$=$};
\draw[mystealth] (05) to (03);

\draw[dashed] (11) to [out=45, in=135] (15);
\node at (13) [above=0.7\mysize] {\Tiny$=$};
\draw[mystealth] (13) to (11);
\draw[mystealth] (13) to (15);

\node (100) at (00) [below=2\mysize-0.2cm]{$\cdots$};
\vertex (101) at (01) [below=2\mysize] {};
\vertex (103) at (03) [below=2\mysize] {};
\vertex (105) at (05) [below=2\mysize] {};
\node (107) at (07) [below=2\mysize-0.2cm] {$\cdots$};

\node (108) at (08) [below=2\mysize] {$\equiv$};
\node (109) at (09) [below=2\mysize] {(1b)};

\node (110) at (10) [below=2\mysize-0.2cm] {$\cdots$};
\vertex (111) at (11) [below=2\mysize] {};
\vertex (113) at (13) [below=2\mysize] {};
\vertex (115) at (15) [below=2\mysize] {};
\node (117) at (17) [below=2\mysize-0.2cm] {$\cdots$};

\draw[mystealth] (105) to [out=135, in=45] (101);
\draw[dashed] (101) -- (103) node [midway, above = -0.2mm] {\Tiny$\ne$};
\draw[mystealth] (105) to (103);

\draw[dashed] (111) to [out=45, in=135] (115);
\node at (113) [above=0.7\mysize] {\Tiny$\ne$};
\draw[mystealth] (113) to (111);
\draw[mystealth] (113) to (115);

\node (200) at (00) [below=4\mysize-0.2cm] {$\cdots$};
\vertex (201) at (01) [below=4\mysize] {};
\vertex (203) at (03) [below=4\mysize] {};
\vertex (205) at (05) [below=4\mysize] {};
\node (207) at (07) [below=4\mysize-0.2cm] {$\cdots$};

\node (208) at (08) [below=4\mysize] {$\equiv$};
\node (209) at (09) [below=4\mysize] {(1c)};

\node (210) at (10) [below=4\mysize-0.2cm] {$\cdots$};
\vertex (211) at (11) [below=4\mysize] {};
\vertex (213) at (13) [below=4\mysize] {};
\vertex (215) at (15) [below=4\mysize] {};
\node (217) at (17) [below=4\mysize-0.2cm] {$\cdots$};

\draw[mystealth] (205) to [out=135, in=45] (201);
\draw[mystealth] (201) to (203);
\draw[mystealth] (205) to (203);

\draw[mystealth] (211) to [out=45, in=135] (215);
\draw[mystealth] (213) to (211);
\draw[mystealth] (213) to (215);

\node (300) at (00) [below=6\mysize-0.2cm] {$\cdots$};
\vertex (301) at (01) [below=6\mysize] {};
\vertex (303) at (03) [below=6\mysize] {};
\vertex (305) at (05) [below=6\mysize] {};
\node (307) at (07) [below=6\mysize-0.2cm] {$\cdots$};

\node (308) at (08) [below=6\mysize] {$\equiv$};
\node (309) at (09) [below=6\mysize] {(2a)};

\node (310) at (10) [below=6\mysize-0.2cm] {$\cdots$};
\vertex (311) at (11) [below=6\mysize] {};
\vertex (313) at (13) [below=6\mysize] {};
\vertex (315) at (15) [below=6\mysize] {};
\node (317) at (17) [below=6\mysize-0.2cm] {$\cdots$};

\draw[mystealth] (305) to [out=135, in=45] (301);
\draw[mystealth] (303) to (301);
\draw[dashed] (303) -- (305) node [midway, above = -0.2mm] {\Tiny$=$};

\draw[dashed] (311) to [out=45, in=135] (315);
\draw[mystealth] (311) to (313);
\draw[mystealth] (315) to (313);
\node at (313) [above=0.7\mysize] {\Tiny$=$};

\node (400) at (00) [below=8\mysize-0.2cm] {$\cdots$};
\vertex (401) at (01) [below=8\mysize] {};
\vertex (403) at (03) [below=8\mysize] {};
\vertex (405) at (05) [below=8\mysize] {};
\node (407) at (07) [below=8\mysize-0.2cm] {$\cdots$};

\node (408) at (08) [below=8\mysize] {$\equiv$};
\node (409) at (09) [below=8\mysize] {(2b)};

\node (410) at (10) [below=8\mysize-0.2cm] {$\cdots$};
\vertex (411) at (11) [below=8\mysize] {};
\vertex (413) at (13) [below=8\mysize] {};
\vertex (415) at (15) [below=8\mysize] {};
\node (417) at (17) [below=8\mysize-0.2cm] {$\cdots$};

\draw[mystealth] (405) to [out=135, in=45] (401);
\draw[mystealth] (403) to (401);
\draw[dashed] (403) -- (405) node [midway, above = -0.2mm] {\Tiny$\ne$};

\draw[dashed] (411) to [out=45, in=135] (415);
\draw[mystealth] (411) to (413);
\draw[mystealth] (415) to (413);
\node at (413) [above=0.7\mysize] {\Tiny$\ne$};

\node (500) at (00) [below=10\mysize-0.2cm] {$\cdots$};
\vertex (501) at (01) [below=10\mysize] {};
\vertex (503) at (03) [below=10\mysize] {};
\vertex (505) at (05) [below=10\mysize] {};
\node (507) at (07) [below=10\mysize-0.2cm] {$\cdots$};

\node (508) at (08) [below=10\mysize] {$\equiv$};
\node (509) at (09) [below=10\mysize] {(2c)};

\node (510) at (10) [below=10\mysize-0.2cm] {$\cdots$};
\vertex (511) at (11) [below=10\mysize] {};
\vertex (513) at (13) [below=10\mysize] {};
\vertex (515) at (15) [below=10\mysize] {};
\node (517) at (17) [below=10\mysize-0.2cm] {$\cdots$};

\draw[mystealth] (505) to [out=135, in=45] (501);
\draw[mystealth] (503) to (501);
\draw[mystealth] (503) to (505);

\draw[mystealth] (511) to [out=45, in=135] (515);
\draw[mystealth] (511) to (513);
\draw[mystealth] (515) to (513);

\node (600) at (00) [below=12\mysize-0.2cm] {$\cdots$};
\vertex (601) at (01) [below=12\mysize] {};
\vertex (603) at (03) [below=12\mysize] {};
\vertex (605) at (05) [below=12\mysize] {};
\node (607) at (07) [below=12\mysize-0.2cm] {$\cdots$};

\node (608) at (08) [below=12\mysize] {$\equiv$};
\node (609) at (09) [below=12\mysize] {(3)};

\node (610) at (10) [below=12\mysize-0.2cm] {$\cdots$};
\vertex (611) at (11) [below=12\mysize] {};
\vertex (613) at (13) [below=12\mysize] {};
\vertex (615) at (15) [below=12\mysize] {};
\node (617) at (17) [below=12\mysize-0.2cm] {$\cdots$};
\draw[mystealth] (605) to (603);
\draw[dashed] (601) -- (603) node [midway, above = -0.2mm] {\Tiny$\ne$};
\draw[dashed] (601) to [out=45, in=135] (605);
\node at (603) [above=0.7\mysize] {\Tiny$\ne$};

\draw[mystealth] (613) to (611);
\draw[dashed] (613) -- (615) node [midway, above = -0.2mm] {\Tiny$\ne$};
\draw[dashed] (611) to [out=45, in=135] (615);
\node at (613) [above=0.7\mysize] {\Tiny$\ne$};
\end{tikzpicture}
\end{center}
\caption{\small \label{fig:diagram Iplac} The congruences  $M\equiv \text{\rm Swap}_{i}(M)$ (1a)--(1c),
$M \equiv \text{\rm Swap}_{i-1}(M)$ (2a)--(2c), and $M\equiv\text{\rm Swap}_{i}(\text{\rm Swap}_{i-1}(M))$ (3)
corresponding to the elements defining $\Iplac$. Note that edges $M_{j,k}$ with exactly one of $j,k$ in $\{i-1,i,i+1\}$ are not drawn but must be adjusted as in Figure~\ref{fig:swap}.}
\end{figure}
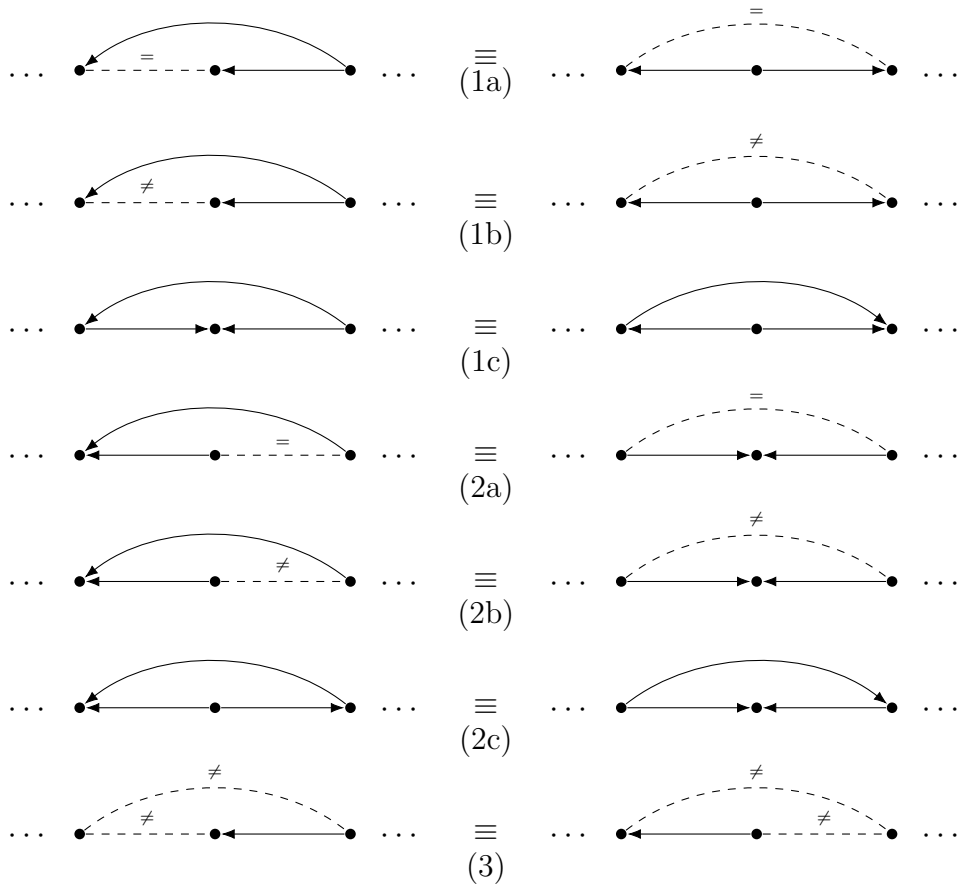

Let $I_H$ be the ideal of $\A$ consisting of the linear span of all arrow diagrams in $\Iplac$ as well as the following elements (for any arrow diagram $M$):
\begin{align}
\text{$M-\text{\rm Swap}_{i}(M)$ \ \ for $1 \le i < |\ver(M)|$ such that $M_{i,i+1}=\ell$}.
\end{align}
Elements in $I_H$ map via $\text{\rm Eval}_P$ into the ideal $I_H^P$ of $\U_P$ for a \threeone-free poset $P$.

%
%

Let $\Ipol$ be the ideal of $\A$ consisting of the linear span of all arrow diagrams in $I_\text{\rm (3+1)}$ as well as the following elements (for any arrow diagram $M$):
\begin{align}
\text{$M-\text{\rm Swap}_{i}(M)$ \ \ \ for all $1\leq i < |\ver(M)|$.}
\end{align}

\subsection{Injection into the product of  $\U_P$'s}

Let $P$ be a  poset. For a word $\e{w}=\e{w}_1\e{w}_2\cdots \e{w}_d\in\U_P^*$, the \textit{origin} of $\e{w}$,
 denoted $\text{\rm orig}_P(\e{w})$, is the primitive arrow diagram $M\in \A_d$ satisfying (for all distinct $i,j\in[d]$)
\begin{itemize}
\item $M_{i,j}=r$ \ whenever $\e{w}_i \to_P \e{w}_j$,
\item $M_{i,j}=e$ \ whenever $\e{w}_i= \e{w}_j$,
\item $M_{i,j}=n$ \ whenever $\e{w}_i \dote_P \e{w}_j$ and $\e{w}_i \ne \e{w}_j$.
\end{itemize}
For example, with  $P$ as in \eqref{eq:ex P fillings},
$\text{\rm orig}_P(\e{adcb}) = \text{\rm orig}_P(\e{cbad})$ is equal to
the primitive arrow diagram  in \eqref{eq:ex P fillings}.
As this example suggests, in general we have
\begin{align}
\label{eq: fibers fillings}
\parbox{10cm}{
the map $\e{w} \mapsto \text{\rm orig}_P(\e{w})$
from words in  $\U_P^*$ to primitive \\[1mm] arrow diagrams has fibers given by $\text{\rm Fill}_P(M)$.}
\end{align}

\begin{remark}
The difference between primitive arrow diagrams and words in  $\U_P$ can be explained informally as follows:
for simplicity, let's compare primitive arrow diagrams  $M \in \A_d \setminus I_\text{closed}$ with no $e$ edges vs.
words $\e{w} \in \U_P^*$ of content $P$ for a poset  $P$ of size  $d$;  the data encoded
by the former is the clique poset $Q_M$ of $M$ 
together with position labels in  $[d]$, while the data encoded by the latter is the
poset $P = Q_{\text{\rm orig}_P(\e{w})}$ together with position labels in  $[d]$ and value labels in  $P$.
The map  $\e{w} \mapsto \text{\rm orig}_P(\e{w})$ removes value labels, and the map  $M \mapsto Q_M$ removes position labels. There are  $|\aut(P)|$ many $\e{w}$ with a fixed origin and  $d!/|\aut(P)|$ many  $M$ whose clique poset
$Q_M$ is $P$.
For example,  with the poset  $P$ and  the arrow diagram  $M$ from \eqref{eq:ex P fillings},  $|\aut(P)| = 2$ is the size of $\{\e{w} : \text{\rm orig}_P(\e{w}) = M\}=\text{\rm Fill}_P(M)$.
\end{remark}

\begin{lemma}
\label{lem:origins}
\
\begin{enumerate}[label=(\alph*)]
\item For $Z\in \A$, the coefficient of a word $\e{w}$ appearing in $\text{\rm Eval}_P(Z)$ is equal to the coefficient of $\text{\rm orig}_P(\e{w})$ in the expansion of $Z$ into primitive arrow diagrams.
\item For all primitive arrow diagrams $M\not\in I_\text{closed}$, if $Q$ is the clique poset of $M$, then there exists a word $\e{w}\in\U_{Q}$ such that $M=\text{\rm orig}_{Q}(\e{w})$.
\item If two words $\e{w}^1,\e{w}^2$ are related by a $P$-Knuth transformation,
then their origins $\text{\rm orig}_P(\e{w}^1)$ and $\text{\rm orig}_P(\e{w}^2)$ are related
by the analogous transformation corresponding to a defining element of $\Iplac$
(cf. Definitions \ref{def-switchboard}/\ref{def:arrow Iplac} (1)--(3) and Figure \ref{fig:diagram Iplac}).
%
\end{enumerate}
\end{lemma}
\begin{proof}
Statement (a) follows from  \eqref{eq: fibers fillings}, (b) from \eqref{eq: fibers fillings} and
Proposition \ref{p I closed},
and (c) follows directly from the definitions of  $\text{\rm Swap}_i$, $\Iplac$, and  $\Iplacp{P}$.
\end{proof}

Combining the homomorphisms $\text{\rm Eval}_P:\A\rightarrow \U_P$ over all (finite) \threeone-free posets $P$, we obtain
a ring homomorphism $\phi: \A\rightarrow \prod\limits_{\text{\threeone-free }P} \U_P$.
By Proposition \ref{p 31free primitive},  $I_\text{\rm (3+1)} \subset \ker \varphi$ and hence we obtain a map
$\varphi: \A/I_\text{\rm (3+1)} \rightarrow \prod\limits_{\text{\threeone-free }P} \U_P$.
It follows from \eqref{eq: fibers fillings} that this map is injective
and its image consists of the elements  with constant
coefficient across the set of monomials  $\bigsqcup_P \{\mathbf{u}_\e{w} \in \U_P : \text{\rm orig}_P(\e{w}) = M\}$,
for each primitive arrow diagram $M$.

We now prove similar statements for the quotients  $\A/\Iplac$,  $\A/I_H$, and  $\A/\Ipol$.
Since $\varphi(\Iplac) \subset \prod\limits_{\text{\threeone-free }P} \Iplacp{P}$ by \eqref{eq: Eval Iplac},
we obtain an induced map $\varphi_\text{\rm plac}: \A/\Iplac \rightarrow \prod\limits_{\text{\threeone-free }P} \U_P/\Iplacp{P}$ between quotient spaces.  We define the maps $\varphi_H$ and $\varphi_\text{\rm pol}$ similarly.

\begin{proposition}
\label{p inject}
Each of the following maps is an injective ring homomorphism:
\begin{align}
\label{eq:p inject1}
\varphi_\text{\rm plac}: \A/\Iplac &\hookrightarrow \prod\limits_{\text{\threeone-free }P} \U_P/\Iplacp{P}\\
\label{eq:p inject2}
\varphi_H: \A/I_H \  &\hookrightarrow \prod\limits_{\text{\threeone-free }P} \U_P/I_H^P\\
\label{eq:p inject3}
\varphi_\text{\rm pol}: \A/\Ipol \spa &\hookrightarrow \prod\limits_{\text{\threeone-free }P} \U_P/\Ipol^P
\end{align}
\end{proposition}

\begin{proof}
We prove only \eqref{eq:p inject1} as the proofs of \eqref{eq:p inject2} and \eqref{eq:p inject3} are similar.

Let $Z\in \A$ be in the kernel of $\varphi_\text{\rm plac}$. It suffices to show that $Z\in \Iplac$. Let $M\not\in I_\text{\rm (3+1)}$ be a primitive arrow diagram with a positive coefficient in $Z$ (if none exist, then we have either $Z\in I_\text{\rm (3+1)}$ or $Z$ has all negative coefficients among $M\not\in I_\text{\rm (3+1)}$, the latter contradicting the fact that $Z$ is in the kernel of $\varphi_\text{\rm plac}$). Let $Q = Q_M$ be the clique poset of $M$. By Lemma \ref{lem:origins} (b) there exists $\e{w}^1\in\U_Q^*$ satisfying $M=\text{\rm orig}_Q(\e{w}^1)$, and by Lemma \ref{lem:origins} (a) $\e{w}$ has a positive coefficient in $\text{\rm Eval}_Q(Z)$. Since $Z$ is in the kernel of $\varphi_\text{\rm plac}$, we have that $\text{\rm Eval}_Q(Z)\in \Iplacp{Q}$, and so there exists $\e{w}^2\in\U_Q^*$ which is congruent to $\e{w}^1$ modulo $\Iplacp{Q}$ and has a negative coefficient in $\text{\rm Eval}_Q(Z)$. Let $L=\text{\rm orig}_Q(\e{w}^2)$. By Lemma \ref{lem:origins} (c),
$M- L \in \Iplac$, and by Lemma \ref{lem:origins} (a) $L$ has a negative coefficient in $Z$.
By \eqref{eq: Eval Iplac}, $M-L$ is in the kernel of $\varphi_\text{\rm plac}$, so let $Z'=Z-(M-L)$ so that $Z'$ is in the kernel of $\varphi_\text{\rm plac}$.
It now suffices to prove $Z'\in \Iplac$ since this implies $Z\in \Iplac$.
We repeat this process on $Z'$ and proceed by induction on the sum of the positive coefficients of $Z$ among primitive arrow diagrams $M\not\in I_\text{\rm (3+1)}$.
%
\end{proof}

\subsection{Further strengthenings of the Stanley-Stembridge conjecture}
\label{ss strengthening SS}


We reformulate our main positivity statement and conjecture from
\S\ref{ss Schur pos P Knuth}, \S\ref{ss e pos} in terms of the arrow algebra.
This naturally suggests ways to further strengthen this conjecture, in the spirit that
a stronger positivity conjecture is easier to solve since it makes it easier to single out a
canonical and elegant positive formula from the multitude.
We begin by defining analogs  $E_\lambda$,  $\mathfrak{J}_\lambda$, $\mathfrak{m}_\lambda$ in the arrow algebra of the elements $\mathfrak{e}_\lambda^P(\mathbf{u})$, $\mathfrak{J}_\lambda^P(\mathbf{u})$, $\mathfrak{m}_\lambda^P(\mathbf{u})$ of $\U_P$.

Define $E_d \in \A_d$ to be the arrow diagram with edge types $(E_d)_{i,j}$ for
$1 \le i < j \le d$ given by
\begin{align}
(E_d)_{i,j}=
\begin{cases}
\ell & \text{ if  $j = i+1$,} \\
r+\ell+e+n & \text{ otherwise.}
\end{cases}
\end{align}
For any sequence  $\alpha \in \ZZ_{\ge 0}^\ell$, define the arrow diagram  $E_\alpha$ by
\begin{align}
E_\alpha = E_{\alpha_1} * E_{\alpha_2} * \cdots * E_{\alpha_\ell}.
\end{align}
For example,
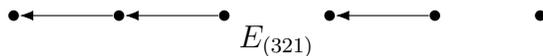
\begin{center}
\begin{tikzpicture}[scale = 0.7]
\vertex (1) at (0,0) {};
\vertex (2) at (2,0) {};
\vertex (3) at (4,0) {};
\vertex (4) at (6,0) {};
\vertex (5) at (8,0) {};
\vertex (6) at (10,0) {};
\draw[mystealth] (2) to (1);
\draw[mystealth] (3) to (2);

\draw[mystealth] (5) to (4);
\node at (5,-0.5) {$E_{(321)}$};
\end{tikzpicture}
\end{center}
Immediate from the definitions, we have
\begin{align}
\label{eq: Eval E}
\text{\rm Eval}_P(E_\alpha)=e_{\alpha_1}^P(\mathbf{u}) \cdots e_{\alpha_\ell}^P(\mathbf{u}).
\end{align}
Also note that by Theorem \ref{th:IplacP elem commute} and \eqref{eq:p inject1},
\begin{align}
\label{eq: Elambdas commute}
\quad \quad E_k E_\ell \equiv E_\ell E_k \ \bmod \spa \Iplac   \ \ \ \text{(and mod $I_H$ and mod  $\Ipol$)}.
\end{align}

Let $\lambda$ be a partition of  $N$,
and let $\bx{\lambda}{1}, \dots,
\bx{\lambda}{N}$ denote the boxes of the Young diagram of $\lambda$ read from bottom to top in each column and with columns read left to right.
Define $D_\lambda \in \A_N$ to be the arrow diagram  whose edge types  $(D_\lambda)_{i,j}$ for
$1 \le i < j \le N$ are given by
\begin{align}
\label{eq: Dlambda}
(D_{\lambda})_{i,j} =
\begin{cases}
\ell & \text{ if $\bx{\lambda}{i}$ is immediately south of $\bx{\lambda}{j}$,} \\
r+e+n & \text{ if $\bx{\lambda}{i}$ is immediately west of $\bx{\lambda}{j}$,} \\
r+\ell+e+n  & \text{ otherwise}.
\end{cases}
\end{align}
These edge types encode the same filling conditions as a $P$-tableau:
$\e{w}$ is a $P$-filling of $D_\lambda$ if and only if it is the column reading word of $P$-tableau of shape  $\lambda$.  Hence also
\begin{align}
\label{eq: Eval D}
\text{\rm Eval}_P(D_\lambda)=  \sum_{T \in \SSYT_P(\lambda)} \mathbf{u}_{\creading(T)}.
\end{align}
The arrow diagram $D_{(433)}$ is drawn below, but with
vertices numbered instead of with our usual convention of displaying vertices left to right:
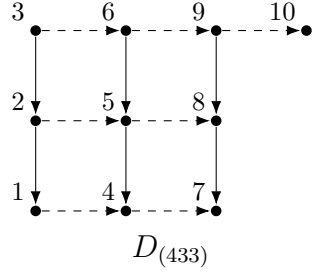
\begin{center}
\begin{tikzpicture}[scale = 0.6]
\vertex (1) at (0,0) {};
\vertex (2) at (0,2) {};
\vertex (3) at (0,4) {};
\vertex (4) at (2,0) {};
\vertex (5) at (2,2) {};
\vertex (6) at (2,4) {};
\vertex (7) at (4,0) {};
\vertex (8) at (4,2) {};
\vertex (9) at (4,4) {};
\vertex (10) at (6,4) {};

\draw[mystealth] (3) to (2);
\draw[mystealth] (2) to (1);

\draw[mystealth] (6) to (5);
\draw[mystealth] (5) to (4);

\draw[mystealth] (9) to (8);
\draw[mystealth] (8) to (7);

\draw[mystealth][dashed] (3) to (6);
\draw[mystealth][dashed] (6) to (9);
\draw[mystealth][dashed] (9) to (10);

\draw[mystealth][dashed] (2) to (5);
\draw[mystealth][dashed] (5) to (8);

\draw[mystealth][dashed] (4) to (7);
\draw[mystealth][dashed] (1) to (4);

\node at (1) [above left] {\footnotesize $1$};
\node at (2) [above left] {\footnotesize $2$};
\node at (3) [above left] {\footnotesize $3$};
\node at (4) [above left] {\footnotesize $4$};
\node at (5) [above left] {\footnotesize $5$};
\node at (6) [above left] {\footnotesize $6$};
\node at (7) [above left] {\footnotesize $7$};
\node at (8) [above left] {\footnotesize $8$};
\node at (9) [above left] {\footnotesize $9$};
\node at (10) [above left] {\footnotesize $10$};
\node at (10) [above right] {\hphantom{\footnotesize $3$}};

\node at (3,-0.9) {$D_{(433)}$};
\end{tikzpicture}
\end{center}


For a partition  $\lambda$, define the \emph{noncommutative arrow Schur function}  $\mathfrak{J}_\lambda \in \A$ by
\begin{equation}
\label{eq:ncschur diagram-via-e}
\mathfrak{J}_\lambda = \sum_{\pi\in \SS_{k}}
\sgn(\pi) \, E_{\lambda'_1+\pi(1)-1} * E_{\lambda'_2+\pi(2)-2} * \cdots *
E_{\lambda'_{k}+\pi(k)-k},
\end{equation}
where  $k = \lambda_1$.
Define the \emph{noncommutative arrow monomial function}  $\mathfrak{m}_\lambda \in \A$ by
\begin{equation}
\label{eq:ncmon diagram-via-e}
\mathfrak{m}_\lambda = \sum_{\mu} (K^{-1})_{\lambda \mu} \,
\mathfrak{J}_\mu.
\end{equation}
By \eqref{eq: Eval E} and the definitions \eqref{eq:P J def}, \eqref{eq:P m def}  of  $\mathfrak{J}_\lambda^P(\mathbf{u})$ and
$\mathfrak{m}_\lambda^P(\mathbf{u})$,
for every  poset  $P$,
\begin{align}
\label{e Eval J}
\text{\rm Eval}_P(\mathfrak{J}_\lambda) = \mathfrak{J}_\lambda^P(\mathbf{u}), 
\qquad \quad
\text{\rm Eval}_P(\mathfrak{m}_\lambda) = \mathfrak{m}_\lambda^P(\mathbf{u}).
\end{align}

An element $Z\in \A$ is \emph{diagram positive modulo an ideal~$I$}
if $Z$ can be written as an element of $I$ plus a nonnegative integer combination of
primitive arrow diagrams.

\begin{proposition}
For any $Z \in \A$, the following are equivalent:
\begin{enumerate}[label=(\roman*)]
\item $Z$ is diagram positive modulo $\Iplac$.
\item $\text{\rm Eval}_P(Z)$ is  $\mathbf{u}$-monomial positive modulo $\Iplacp{P}$ for all \threeone-free posets $P$.
\end{enumerate}
Moreover, the analogous statement holds with  $I_H$ and $I_H^P$ in place of $\Iplac$ and $\Iplacp{P}$, respectively, or with $I_\text{\rm pol}$ and $I_\text{\rm pol}^P$ in place of $\Iplac$ and  $\Iplacp{P}$.
In particular,
\begin{align}
\label{eq: IH arrow alg pos}
 \parbox{13.3cm}{$\mathfrak{m}_\lambda$ is diagram positive modulo $I_H$ (resp.  $I_\text{\rm pol}$) if and only if
$\mathfrak{m}_\lambda^P(\mathbf{u})$ is  $\mathbf{u}$-monomial positive modulo $I_H^P$ (resp.  $I_\text{\rm pol}^P$) for all \threeone-free posets $P$.}
\end{align}
\end{proposition}
\begin{proof}
Let $Z\in \A$ and suppose condition (i) holds. Let $P$ be a \threeone-free poset. Since $\text{\rm Eval}_P(\Iplac) \subset \Iplacp{P}$, we have that $\text{\rm Eval}_P(Z)$ is  $\mathbf{u}$-monomial positive modulo $\Iplacp{P}$.

Now, let $Z\in \A$ and suppose condition (ii) holds. Let $M\not\in I_\text{\rm (3+1)}$ be a primitive arrow diagram with a negative coefficient in $Z$ (if none exist, $Z$ is diagram positive modulo $\Iplac$).
Let $Q = Q_M$ be the clique poset of $M$.
Let  $\e{w}^1\in \text{\rm Fill}_{Q}(M)$, so that  $\e{w}^1$ has a negative coefficient in $\text{\rm Eval}_{Q}(Z)$. Since $\text{\rm Eval}_{Q}(Z)$ is $\mathbf{u}$-monomial positive modulo $\Iplacp{Q}$, there exists a word $\e{w}^2\in\U_{Q}^*$ which satisfies $\e{w}^1\equiv \e{w}^2 \bmod \Iplacp{Q}$ and has a positive coefficient in $\text{\rm Eval}_{Q}(Z)$. Let $L=\text{\rm orig}_{Q}(\e{w}^2)$. Then by Lemma \ref{lem:origins} (c), we have $M- L \in \Iplac$, and by Lemma \ref{lem:origins} (a), $L$ has a positive coefficient in $Z$. Now, let $Z'=Z+(M-L)$. Since $\text{\rm Eval}_P(M-L) \in \Iplacp{P}$ for all \threeone-free posets $P$, we have that $\text{\rm Eval}_P(Z')$ is $\mathbf{u}$-monomial positive modulo $\Iplacp{P}$ for all \threeone-free
posets $P$. Further, since $(M-L)\in \Iplac$, the condition that $Z'$ is diagram positive modulo $\Iplac$ implies that $Z$ is diagram positive modulo $\Iplac$. We repeat this process on $Z'$ and proceed by induction on the sum of the negative coefficients of $Z$ among primitive arrow diagrams $M\not\in I_\text{\rm (3+1)}$.
This shows that condition (i) holds.

This proof also works with $I_H$ and $I_H^P$ (or $\Ipol$ and $\Ipol^P$) in place of $\Iplac$ and $\Iplacp{P}$.
\end{proof}


Proposition \ref{p inject}
 allows us to restate Theorem \ref{th:IplacP-positivity} in an elegant way in the arrow algebra.
Surprisingly, there is now no mention of a specific \threeone-free poset in the statement.

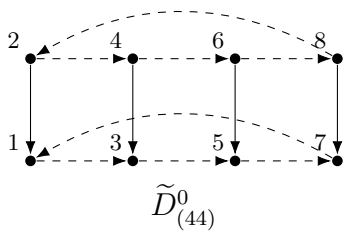
\begin{theorem}
\label{t plactic positivity arrow alg}
For any partition  $\lambda$,
\begin{align}
\mathfrak{J}_\lambda \equiv D_\lambda \ \bmod \, \Iplac
\end{align}
and hence $\mathfrak{J}_\lambda$ is diagram positive modulo  $\Iplac$.
\end{theorem}
\begin{proof}
By \eqref{e Eval J}, \eqref{eq: Eval D}, and Theorem \ref{th:IplacP-positivity},  $\varphi_\text{\rm plac}(\mathfrak{J}_\lambda - D_\lambda) = 0$, hence $\mathfrak{J}_\lambda - D_\lambda = 0$ in  $\A/\Iplac$ by  \eqref{eq:p inject1}.
\end{proof}

For example, for  $\lambda = (22)$, this says that $E_{22}-E_{31}\equiv D_{22}\mod \Iplac$, or in pictures,
\begin{center}
\begin{tikzpicture}[scale = 0.7]
\vertex (1) at (0,0) {};
\vertex (2) at (2,0) {};
\vertex (3) at (4,0) {};
\vertex (4) at (6,0) {};
\node at (7,-0.25) {$-$};
\vertex (5) at (8,0) {};
\vertex (6) at (10,0) {};
\vertex (7) at (12,0) {};
\vertex (8) at (14,0) {};
\node at (3,-0.7) {$E_{22}$};
\node at (11,-0.7) {$E_{31}$};

\draw[mystealth] (4) to (3);
\draw[mystealth] (2) to (1);

\draw[mystealth] (7) to (6);
\draw[mystealth] (6) to (5);

\node at (3,-2.25) {$\equiv$};

\vertex (9) at (4,-2-0.2) {};
\vertex (10) at (6,-2-0.2) {};
\vertex (11) at (8,-2-0.2) {};
\vertex (12) at (10,-2-0.2) {};

\node at (7,-2.7-0.2) {$D_{22}$};

\draw[mystealth] (10) to (9);
\draw[mystealth] (12) to (11);

\draw[mystealth][dashed] (9) to [out=45, in=135] (11);
\draw[mystealth][dashed] (10) to [out=45, in=135] (12);
\end{tikzpicture}
\end{center}

By \eqref{eq: IH arrow alg pos}, Conjecture \ref{cj:m pos} is equivalent to the following conjecture.
\begin{conjecture}
\label{cj: m pos arrow alg}
The elements  $\mathfrak{m}_\lambda$ are
diagram positive modulo  $I_H$.
\end{conjecture}

Since an arrow diagram is a positive sum of primitive arrow diagrams,
a simple way for an element of $\A$
to be diagram positive modulo  $I$ is to be congruent to
a single arrow diagram modulo $I$.  As we just saw in Theorem \ref{t plactic positivity arrow alg}, this happens for  the element $\mathfrak{J}_\lambda$ with the ideal $\Iplac$.
It also happens for certain  $\mathfrak{m}_\lambda$ and  $I_\text{\rm pol}$, as we'll now see.

Let  $\lambda$ be a partition, set $k = \lambda_1$, and let  $m \in  \{\lambda'_1- \lambda'_k, \dots, \lambda'_1\}$.
Define an arrow diagram $\widetilde{D}^m_\lambda$ similar to $D_\lambda$ by
\begin{align*}
(\widetilde{D}^m_{\lambda})_{i,j} =
\begin{cases}
\ell & \text{if $\bx{\lambda}{i}$ is immediately south of $\bx{\lambda}{j}$,} \\
r+e+n & \text{if $\bx{\lambda}{i}$ is immediately west of $\bx{\lambda}{j}$,} \\
\ell+e+n & \text{if $\bx{\lambda}{i} = (r+m,1)$ and $\bx{\lambda}{j} = (r,k)$ for some  $r \in [\lambda'_1-m]$,} \\
r+\ell+e+n  & \text{otherwise},
\end{cases}
\end{align*}
with notation as in \eqref{eq: Dlambda} so that, for instance, $\bx{\lambda}{j} = (r,k)$ means the  $j$-th box in column reading order is the easternmost box in the  $r$-th row of  $\lambda$. Also,  $k=2,$
$m=0$ is a degenerate case in which the second and third cases just above should be replaced by $(\widetilde{D}^m_{\lambda})_{i,j} = e+n$ if $\bx{\lambda}{i}$ is immediately west of $\bx{\lambda}{j}$.

Here is the example $\widetilde{D}^0_{(44)}$:
\begin{center}
\begin{tikzpicture}[scale = 0.68]
\vertex (1) at (0,0) {};
\vertex (2) at (0,2) {};
\vertex (3) at (2,0) {};
\vertex (4) at (2,2) {};
\vertex (5) at (4,0) {};
\vertex (6) at (4,2) {};
\vertex (7) at (6,0) {};
\vertex (8) at (6,2) {};

\draw[mystealth] (2) to (1);

\draw[mystealth] (4) to (3);

\draw[mystealth] (6) to (5);

\draw[mystealth] (8) to (7);

\draw[mystealth][dashed] (2) to (4);
\draw[mystealth][dashed] (4) to (6);
\draw[mystealth][dashed] (6) to (8);

\draw[mystealth][dashed] (1) to (3);
\draw[mystealth][dashed] (3) to (5);
\draw[mystealth][dashed] (5) to (7);

\draw[mystealth][dashed] (8) to [out=150, in=30] (2);
\draw[mystealth][dashed] (7) to [out=150, in=30] (1);

\node at (1) [above left] {\scriptsize $1$};
\node at (2) [above left] {\scriptsize $2$};
\node at (3) [above left] {\scriptsize $3$};
\node at (4) [above left] {\scriptsize $4$};
\node at (5) [above left] {\scriptsize $5$};
\node at (6) [above left] {\scriptsize $6$};
\node at (7) [above left] {\scriptsize $7$};
\node at (8) [above left] {\scriptsize $8$};
\node at (8) [above right] {\hphantom{\scriptsize $3$}};

\node at (3,-0.9) {$\widetilde{D}^0_{(44)}$};
\end{tikzpicture}
\end{center}

The following result is a consequence of \cite[Theorem 4.7 (v-b)]{CHSS}, which is based off
\cite[Theorem 2.8]{StembridgeImmanant} (see also \cite[\S4]{Isaiah}).
\begin{theorem}
\label{t rectangle m}
For a rectangle shape $\lambda = (k^\ell)$ and any \threeone-free poset  $P$,
\begin{align}
\label{eq:t rectangle m}
\mathfrak{m}^P_\lambda(\mathbf{u}) \equiv
 \text{\rm Eval}_P(\widetilde{D}^0_\lambda)
= \sum_{\e{w} \in \text{\rm Fill}_P(\widetilde{D}^0_\lambda)}\mathbf{u}_{\e{w}}  \ \ \bmod \, \Ipol^P.
\end{align}
Hence the coefficient of  $e_\lambda(\mathbf{x})$ in the multicolored chromatic symmetric function  $X_{\inc(P)}^{\beta}(\mathbf{x})$
is the number of  $P$-fillings of  $\widetilde{D}^0_\lambda$ of content  $\beta$.
\end{theorem}

By \eqref{eq:p inject2}, Theorem \ref{t rectangle m} has the following equivalent restatement in the arrow algebra.
\begin{theorem}
\label{t m rectangle positivity arrow alg}
For a rectangle shape $\lambda$,
\begin{align}
\label{eq:t rectangle m}
\mathfrak{m}_\lambda \equiv  \widetilde{D}^0_\lambda  \ \bmod \, \Ipol.
\end{align}
\end{theorem}

Theorems \ref{t plactic positivity arrow alg} and \ref{t m rectangle positivity arrow alg} suggest the following lines of attack for Conjecture \ref{cj: m pos arrow alg}, which further build on the idea
that strengthening positivity conjectures as much as possible is a fruitful way to search for an elegant positive formula.


\begin{question}
\label{q: m diagram}
\
\begin{itemize}
\item[(i)] Is $\mathfrak{m}_\lambda$ congruent to a single arrow diagram modulo  $I_\text{\rm pol}$?
\item[(ii)] We suspect that (i) is too good to hope for in general, but is $\mathfrak{m}_\lambda$ congruent modulo  $I_\text{\rm pol}$ to a sum of arrow diagrams using only solid arrows, dashed arrows, and empty edges?
\item[(iii)] Does the strengthening of (ii) hold in which  $I_\text{\rm pol}$ is replaced with  $I_H$?
\end{itemize}
\end{question}

These questions make precise the idea that in order to rewrite
$\mathfrak{m}_\lambda = \sum_{\mu} (K^{-1})_{\lambda \mu} \,
\mathfrak{J}_\mu \equiv \sum_{\mu} (K^{-1})_{\lambda \mu} \,
D_\mu$ as a positive sum of arrow diagrams mod $I_H$ or $I_\text{\rm pol}$, we should try to cancel positive and negative terms while expanding dashed arrows and empty edges
as little as possible into the edges types $r,\ell,e,$ and  $n$.

\subsection{Variations on Stanley's  $P$-partitions conjecture}
\label{ss Stanley variations}
The following line of inquiry may help with Question \ref{q: m diagram} and also seems interesting in its own right.
Let $\PP$ be the total order $[N]$,
and let $\overline{\text{\rm Eval}}_\PP(M)$ denote
the image of $\text{\rm Eval}_\PP(M)$ in  $\U_\PP/I_\text{\rm pol}^\PP \cong \ZZ[y_1,\dots,y_N]$.
Consider the following properties of an arrow diagram  $M$:
\begin{itemize}
\item[(a)]  $M$ lies in $\text{span}_\ZZ\{E_\lambda : \text{ partitions } \lambda\}$ in  $\A/I_H$.
\item[(b)]  $M$ lies in $\text{span}_\ZZ\{E_\lambda : \text{ partitions } \lambda\}$ in  $\A/\Ipol$.
\item[(c)]  $\overline{\text{\rm Eval}}_\PP(M)$ is an ordinary symmetric function.
\end{itemize}
We have $(\text{a}) \! \implies \! (\text{b}) \! \implies \! (\text{c})$ since
$\overline{\text{\rm Eval}}_\PP(E_\lambda)$ is the ordinary elementary symmetric function  $e_\lambda(\mathbf{y})$.

\begin{question}
\label{q: which symmetric}
\

\begin{itemize}
\item[(i)] Which arrow diagrams  $M$ satisfy (a), (b), or (c)?
\item[(ii)] What symmetric functions arise as $\overline{\text{\rm Eval}}_\PP(M)$ for such  $M$?
\end{itemize}
\end{question}
As we have seen, $D_\lambda$ and  $\widetilde{D}^0_{(k^\ell)}$ satisfy (a)--(c) and they yield the
Schur functions and rectangular monomial symmetric functions, respectively.

Recall from Remark \ref{r P partition} that Stanley's $(P, \omega)$-partitions and McNamara's
$(P,O)$-partitions can be viewed as  $\PP$-fillings of special arrow diagrams  $M$ and their generating functions
are the same as $\overline{\text{\rm Eval}}_\PP(M)$.
Stanley \cite{StanleyPpartition} conjectured that the generating function of a  $(P, \omega)$-partition
is symmetric if and only if the set of $(P, \omega)$-partitions are the semistandard Young tableaux of a skew shape.

McNamara \cite{McNamara} considered the similar question
``which oriented posets have symmetric generating functions?''
He observes that cylindrical Schur functions can be viewed as generating functions of the oriented posets corresponding to the $\widetilde{D}^m_\lambda$.
He also found several oriented posets not of this form
which have symmetric generating functions \cite[Figure 8]{McNamara}.
Question \ref{q: which symmetric} can be viewed as a variation of McNamara's question and Stanley's conjecture in which we allow more general kinds of fillings
but ask for the stronger properties (a) or (b) rather than just (c).

It is now natural to ask whether the arrow diagrams
$\widetilde{D}^m_\lambda$
and those corresponding to \cite[Figure 8]{McNamara} satisfy (a) or (b).
The fourth author \cite{Isaiah} recently showed that the $\widetilde{D}^m_\lambda$ do indeed satisfy (b), as we now explain.
Let  $\lambda$ be a partition, set $k = \lambda_1$, and let  $m \in  \{\lambda'_1- \lambda'_k, \dots, \lambda'_1\}$.
Define  the \emph{noncommutative $P$-cylindrical Schur function} $\mathfrak{J}_{\lambda/m}^P(\mathbf{u}) \in \U_P$ by \begin{align}
\label{eq:J cylinder}
\mathfrak{J}_{\lambda/m}^P(\mathbf{u}) = \sum_{\substack{a_1, \spa \dots, a_k \in \ZZ \\ a_1 + \spa \cdots \spa + a_k = 0}}
\sum_{\pi\in \SS_{k}}
\sgn(\pi) \, e^P_{a_1(k+m)+\lambda'_1+\pi(1)-1}(\mathbf{u})\,
\cdots e^P_{a_k(k+m)+\lambda'_{k}+\pi(k)-k}(\mathbf{u}).
\end{align}
The \emph{cylindrical Schur function} $s_{\lambda/m}(\mathbf{y})$ is
the image of  $\mathfrak{J}_{\lambda/m}^\PP(\mathbf{u})$ in  $\U_\PP/\Ipol^\PP \cong \ZZ[y_1,\dots,y_N]$,
which is just \eqref{eq:J cylinder} with ordinary  $e_d(\mathbf{y})$'s in place of  $e_d^P(\mathbf{u})$'s
and is the well-known determinantal formula for cylindrical Schur functions.  It is known from \cite{GKcyl,PostnikovCylSchur,McNamara}
that this determinantal definition of  $s_{\lambda/m}(\mathbf{y})$ is equal to a sum over cylindrical tableaux,
which can be written in terms of  $\PP$-fillings of $\widetilde{D}^m_\lambda$ as
\begin{align}
s_{\lambda/m}(\mathbf{y}) = \overline{\text{\rm Eval}}_\PP(\widetilde{D}^m_\lambda).
\end{align}
The following strengthening of this result to any \threeone-free poset is proved in \cite{Isaiah}:

\begin{theorem}
For any \threeone-free poset  $P$, the algebraically defined
$\mathfrak{J}^P_{\lambda/m}(\mathbf{u})$
have the following tableau-style formula modulo  $\Ipol^P$:
\begin{align}
\mathfrak{J}^P_{\lambda/m}(\mathbf{u}) \equiv
\text{\rm Eval}_P(\widetilde{D}^m_\lambda)
= \sum_{\e{w} \in \text{\rm Fill}_P(\widetilde{D}^m_\lambda)}\mathbf{u}_{\e{w}}
 \ \ \ \, \bmod \, \Ipol^P.
\end{align}
\end{theorem}

It then follows from \eqref{eq:p inject3} that
\begin{corollary}
The arrow diagram $\widetilde{D}_\lambda^{m}$ satisfies (b) and in fact, letting  $\mathfrak{J}_{\lambda/m}$ denote the element of  $\A$ obtained by replacing $e_d^P(\mathbf{u})$'s with  $E_d$'s in \eqref{eq:J cylinder},
\begin{align}
\mathfrak{J}_{\lambda/m} \equiv \widetilde{D}^m_{\lambda}  \ \, \bmod \, \Ipol.
\end{align}
\end{corollary}

\section{The  $e_k^P(\mathbf{u})$ commute in the  $P$-plactic algebra}
\label{sec:elem commute}
We prove Theorem \ref{th:IplacP elem commute}.
We begin with preparatory material on ladders and relations in  $\U_P/\Iplacp{P}$
which will be used in other proofs as well.

Throughout this section,  $P$ denotes a  \threeone-free poset.

\subsection{Ladders}
\label{ss ladders}
Recall from Definition \ref{def:ladder} the definition of a ladder and related notions.
We need several facts about ladders and some further definitions.

\begin{lemma}
\label{l ladder}
Let  $H$ be a ladder of a pair of chains  $(C,D)$ in a \threeone-free poset $P$.  Then
\begin{itemize}
\item[(i)]  $H$ is either a 4-cycle or a path.
\item[(ii)] In fact,  $H$ must have one of the forms shown in Figure \ref{fig:ladder}.
\item[(iii)]  $|H\cap C|$ and  $|H \cap D|$ differ in size by at most 1.
\end{itemize}
\end{lemma}
\begin{proof}
As $P$ is \threeone-free, each vertex of $\inc_P(C,D)$
can be incomparable to at most two elements in  $C \sqcup D$, so the degree of each vertex is at most 2. Therefore, as $H$ is connected, $H$ must be a path or an even cycle. As 4-cycles are the only even cycles that can appear as an induced subgraph of an incomparability graph of a poset \cite{Gallai} (it is also easy to prove directly that
this is true of $\inc_P(C,D)$),
$H$ must either be a path or a 4-cycle in $\inc_P(C,D)$.

Statement (ii) follows from (i) and transitivity, using a straightforward case by case analysis.  For instance, it follows directly from transitivity that if $H$ has nonempty intersection with $C$ and  $D$, then  $c \dote_P d$ for  $c$ (resp.  $d$) the smallest element (in  $P$) of $H \cap C$ (resp.  $H\cap D$).
Statement (iii) is immediate from (i).
\end{proof}

\begin{definition}
\label{d ladder stuff}
Let  $C,D$ be chains of a \threeone-free poset  $P$.
The ladders of  $(C,D)$ inherit a total order from  $P$, i.e. for distinct ladders $H, H'$ every element of  $H$
is less than every element of  $H'$ in  $P$ or vice-versa
(this is easily checked with the help of Lemma~\ref{l ladder}~(ii)).
The \emph{ladder decomposition of $(C,D)$} is the sequence $H_m,\dots,H_2,H_1$ of ladders of  $(C,D)$
in decreasing order.

For a ladder  $H$ of  $(C,D)$, we can move its elements in  $C$ into  $D$ and vice-versa
to obtain new chains
$C' = (C \setminus H) \cup (H \cap D)$ and $D' = (D \setminus H) \cup (H \cap C)$.  We say that  $(C',D')$ is the result of a \emph{ladder swap on the ladder $H$ of  $(C,D)$}.  It follows from Lemma \ref{l ladder} (ii) that  $C',D'$ are indeed chains and that they have the same ladder decomposition as  $(C,D)$.
\end{definition}

\subsection{Relations in the $P$-plactic algebra}
\label{ss relations}
We establish several identities in the  $P$-plactic algebra  $\U_P/\Iplacp{P}$ which
can be thought of as versions of
the  $P$-plactic relations \eqref{poset rel knuth etc bca}--\eqref{poset rel bca cab} in which the letters  $\e{a}, \e{b}, \e{c}$ are replaced with words.

\begin{lemma}
\label{l plactic inc word 0}
Let $b$ and  $w$ be elements of  $P$ and $\e{v} \!= \! \e{v}_1 \cdots \e{v}_s$ be a word in  $\U^*_P$  satisfying
\begin{itemize}
\item[(i)] $b \leftto_P v_1 \leftto_P \cdots \leftto_P v_s$,
\item[(ii)] $w \leftto_P v_1$,
\item[(iii)]  $b \dotto_P w$.
\end{itemize}
 Then $u_b \mathbf{u}_\e{v} u_w \equiv u_b u_w \mathbf{u}_\e{v}  \bmod \Iplacp{P}$.
\end{lemma}
\begin{proof}
Apply the relation \eqref{poset rel knuth etc bca}  $s$ times to move  $u_w$ to the left through  $\mathbf{u}_{\e{v}}$.
\end{proof}

\begin{example}
With  $P = \PP_{2}$,
$b = 7,$ $w = 8$, $\e{v} = \e{531},$
 three applications of \eqref{poset rel knuth etc bca}
 gives $u_b \mathbf{u}_{\e{v}} u_w = \mathbf{u}_{\e{75318}} \equiv \mathbf{u}_{\e{75381}} \equiv \mathbf{u}_{\e{75831}} \equiv \mathbf{u}_{\e{78531}} = u_b u_{w}\mathbf{u}_{\e{v}} \bmod   \Iplacp{P}$.
\end{example}

\begin{lemma}
\label{l shuffle chains}
Let $\e{a} = \e{a}_n \e{a}_{n-1} \cdots \e{a}_1$ and  $\e{b} = \e{b}_n \e{b}_{n-1} \cdots \e{b}_1$ and $\e{v}$ be words in  $\U^*_P$  satisfying
\begin{itemize}
\item[(i)]  the words $\e{a} \spa \e{v}$ and $\e{b} \spa \e{v}$ are  $P$-strictly decreasing,
\item[(ii)]  $a_{i} \dotto_P b_i$ for all  $i \in [n]$,
\item[(iii)] $a_{i} \to_P b_{i+1}$ for all  $i \in [n-1]$.
\end{itemize}
 Then $\mathbf{u}_\e{a} \mathbf{u}_\e{v} \mathbf{u}_\e{b}\equiv \mathbf{u}_\e{a} \mathbf{u}_\e{b} \mathbf{u}_\e{v}  \bmod \Iplacp{P}$.
\end{lemma}
\begin{proof}
By repeated applications of Lemma \ref{l plactic inc word 0}, we can move each letter of  $\e{b}$ past  $\e{v}$ and shuffle it into the word  $\e{a}$, starting from  $\e{b}_n$:
\begin{multline*}
\mathbf{u}_\e{a} \mathbf{u}_\e{v} \mathbf{u}_\e{b} \equiv
u_{a_n}   u_{b_n} u_{a_{n-1}} \cdots u_{a_1} \mathbf{u}_\e{v} u_{b_{n-1}} \cdots u_{b_1}
\\
\equiv
u_{a_n}   u_{b_n} u_{a_{n-1}}u_{b_{n-1}}u_{a_{n-2}} \cdots u_{a_1} \mathbf{u}_\e{v} u_{b_{n-2}} \cdots u_{b_1}
\equiv \cdots \equiv
u_{a_n}   u_{b_n} u_{a_{n-1}}u_{b_{n-1}} \cdots  u_{a_1}u_{b_1} \mathbf{u}_\e{v}.
\end{multline*}
We can then undo some of these moves 
to unshuffle  $\e{a}$ and $\e{b}$,
obtaining the desired  $\mathbf{u}_\e{a} \mathbf{u}_\e{b} \mathbf{u}_\e{v}$.
\end{proof}

\begin{lemma}
\label{l shuffle chains 2}
Let $\e{a} = \e{a}_n \e{a}_{n-1} \cdots \e{a}_1$ and  $\e{b} = \e{b}_n \e{b}_{n-1} \cdots \e{b}_1$ and $\e{v}$ be words in  $\U^*_P$  satisfying
\begin{itemize}
\item[(i)]  the words $\e{a} \spa \e{v}$ and $\e{b} \spa \e{v}$ are  $P$-strictly decreasing,
\item[(ii)]  $a_{i} \dote_P b_i$ for all  $i \in [n]$,
\item[(iii)]  $a_{i} \dote_P b_{i+1}$ for all  $i \in [n-1]$,
\item[(iv)] $a_{i} \to_P b_{i+2}$ for all  $i \in [n-2]$.
\end{itemize}
 Then $\mathbf{u}_\e{a} \mathbf{u}_\e{v} \mathbf{u}_\e{b}\equiv \mathbf{u}_\e{a} \mathbf{u}_\e{b} \mathbf{u}_\e{v}  \bmod \Iplacp{P}$.
\end{lemma}
\begin{proof}
Similar to the proof of Lemma \ref{l shuffle chains}, we move the letters of  $\e{b}$ past $\e{v}$
and shuffle them into  $\e{a}$ but now each step uses Lemma \ref{l plactic inc word 0} followed by a single application of the relation \eqref{poset rel bca cab},
\begin{equation}
\label{e shuffle chains}
\begin{gathered}
\mathbf{u}_\e{a} \mathbf{u}_\e{v} \mathbf{u}_\e{b} \equiv
u_{a_n}   u_{a_{n-1}}  u_{b_n} \cdots u_{a_1} \mathbf{u}_\e{v} u_{b_{n-1}} \cdots u_{b_1}
\equiv
 u_{b_n} u_{a_n}   u_{a_{n-1}}  \cdots u_{a_1} \mathbf{u}_\e{v} u_{b_{n-1}} \cdots u_{b_1} \equiv
 \\
u_{b_n} u_{a_n}  u_{b_{n-1}} u_{a_{n-1}}  \cdots u_{a_1} \mathbf{u}_\e{v} u_{b_{n-2}} \cdots u_{b_1}
\equiv \cdots \equiv
u_{b_n} u_{a_n}  u_{b_{n-1}} u_{a_{n-1}}  \cdots u_{b_2} u_{a_2} u_{a_1} u_{b_1} \mathbf{u}_\e{v}.
\end{gathered}
\end{equation}
(The last step to move  $b_1$ through is different---we don't apply the relation \eqref{poset rel bca cab}.)
We can then undo some of these moves 
to unshuffle  $\e{a}$ and $\e{b}$, obtaining the desired  $\mathbf{u}_\e{a} \mathbf{u}_\e{b} \mathbf{u}_\e{v}$.
\end{proof}

\begin{lemma}
\label{l swap chains}
Let $\e{a} = \e{a}_{n+1} \e{a}_{n} \cdots \e{a}_1$ and  $\e{b} = \e{b}_n \e{b}_{n-1} \cdots \e{b}_1$
be $P$-strictly decreasing words in  $\U_P^*$ satisfying
\begin{itemize}
\item[(i)]  $a_{i} \dote_P b_i$ for all  $i \in [n]$,
\item[(ii)]  $a_{i+1} \dote_P b_i$ for all  $i \in [n]$,
\item[(iii)] $a_{i} \to_P b_{i+1}$ for all  $i \in [n-1]$,
\item[(iv)] $a_{i+2} \leftto_P b_{i}$ for all  $i \in [n-1]$.
\end{itemize}
 Then $\mathbf{u}_\e{a}  \mathbf{u}_\e{b}\equiv \mathbf{u}_\e{b} \mathbf{u}_\e{a} \bmod \Iplacp{P}$.
\end{lemma}
\begin{proof}
Almost exactly the same as \eqref{e shuffle chains}, except on the last step we do use the relation \eqref{poset rel bca cab}, we have
\begin{align}
\mathbf{u}_\e{a} \mathbf{u}_\e{b} \equiv
u_{b_n} u_{a_{n+1}}  u_{b_{n}} u_{a_{n}}  \cdots u_{b_3} u_{a_4} u_{b_2} u_{a_3} u_{b_1} u_{a_2} u_{a_1}.
\end{align}
By \eqref{poset rel knuth etc bca} applied to  $u_{b_2} u_{a_3} u_{b_1}$, this is congruent mod  $\Iplacp{P}$ to
\begin{align}
u_{b_n} u_{a_{n+1}}  u_{b_{n}} u_{a_{n}}  \cdots u_{b_3} u_{a_4} u_{b_2} u_{b_1} u_{a_3}  u_{a_2} u_{a_1}.
\end{align}
Next, by Lemma \ref{l plactic inc word 0} with  $b = b_3$,  $\e{v} = b_2 b_1$, and
$w = a_4$, we can move  $b_2 b_1$ past  $a_4$ and next to  $b_3$.  Continuing in this way, we can move chunks
 $b_i b_{i-1}\cdots b_1$ to the left past  $a_i$ and next to  $b_{i+1}$, eventually obtaining  $\mathbf{u}_\e{b} \mathbf{u}_\e{a}$.
\end{proof}

\begin{remark}
\label{r flip}
For a \threeone-free poset  $P$, let  $P^{\text{\rm op}}$ be its dual, with
 $a^{\text{\rm op}} \in P^{\text{\rm op}}$ denoting the counterpart to the element  $a \in P$.
 The
algebra anti-isomorphism $\U_P \to \U_{P^{\text{\rm op}}}$ taking $\e{w}= \e{w}_1\cdots \e{w}_n$ to
$\e{w}_n^{\text{\rm op}} \cdots \e{w}_1^{\text{\rm op}}$ takes  $\Iplacp{P}$ to  $\Iplacp{P^{\text{\rm op}}}$.
Thus the four lemmas above have ``flipped'' counterparts obtained by applying this anti-isomorphism to their hypotheses and conclusions.  For example, here is the flipped version of Lemma \ref{l shuffle chains}:
\end{remark}

\begin{lemma}
\label{l shuffle chains flip}
Let $\e{a} = \e{a}_n \e{a}_{n-1} \cdots \e{a}_1$ and  $\e{b} = \e{b}_n \e{b}_{n-1} \cdots \e{b}_1$ and $\e{w}$ be words in  $\U^*_P$  satisfying
\begin{itemize}
\item[(i)]  the words $\e{w} \spa \e{a} $ and $\e{w}\spa \e{b}$ are  $P$-strictly decreasing,
\item[(ii)]  $a_{i} \dotto_P b_i$ for all  $i \in [n]$,
\item[(iii)] $a_{i} \to_P b_{i+1}$ for all  $i \in [n-1]$.
\end{itemize}
 Then $\mathbf{u}_\e{a} \mathbf{u}_\e{w} \mathbf{u}_\e{b}\equiv \mathbf{u}_\e{w} \mathbf{u}_\e{a} \mathbf{u}_\e{b}   \bmod \Iplacp{P}$.
\end{lemma}

\subsection{Proof of Theorem \ref{th:IplacP elem commute}}
\label{ss proof of e commute}
We need to show that $e_k^P(\mathbf{u}_S)$ and $e_\ell^P(\mathbf{u}_S)$ commute modulo  $\Iplacp{P}$.
Since a subposet of a \threeone-free poset is again \threeone-free, we can without loss of generality assume $S = P$.
Denote by $W^k$ the set of chains of $P$ of length $k$, written in decreasing order. One can identify each such chain with one of the monomials in $e_k^P(\mathbf{u})$. To prove the theorem we exhibit an involution $$\eta: W^k \times W^{\ell} \rightarrow W^{\ell} \times W^k.$$  Such
an $\eta$ can be viewed as an analog of a {\it {combinatorial $R$-matrix}} in the theory of crystals.

Throughout this subsection we let  $(\a,\b) \in  W^k \times W^{\ell}$, with elements
 $\a = \big(a_k \leftto_P a_{k-1} \leftto_P \cdots \leftto_P a_1\big)$ and
  $\b = \big(b_\ell \leftto_P b_{\ell-1} \leftto_P \cdots \leftto_P b_1\big)$.

Let us first review how the combinatorial $R$-matrix works in
the usual case when $P$ is the total order on $\{1,2,\dots, N\}$.
If chains $\a$ and $\b$ have the same length, we have $\eta(\a,\b) = (\a,\b)$.
Otherwise, we pair elements of  $\a$ with elements of  $\b$ as follows: choose any $a_i \le b_j$ such that
there is no element of  $\a$ or  $\b$ in the open interval  $(a_i,b_j)$, and pair  $a_i$ with  $b_j$.
Remove this pair and repeat.  When there are no such pairs, pair the largest element of  $\a$ with the smallest element of  $\b$ (call this a \emph{wraparound pair}). Remove the pair and then repeat until  $\a$ or  $\b$ is used up.  Let us call the result {\it {pairing}}.
Once we have a pairing, $\eta(\a,\b)$ is obtained by just transferring all unpaired elements of $\a$ into $\b$  if  $|\a| > |\b|$, or by transferring all unpaired elements of  $\b$ into  $\a$ if $|\a| < |\b|$.

Note that the set of pairs formed this way does not depend on the order they are chosen.

The pairing can also be thought of in terms of parentheses matching:
associate a word  $w$ in the symbols  ``('' and  ``)'' to  $(\a,\b)$ by marking the set
 $\a \cup \{x+1/2 : x \in \b\}$ on the real line and labeling elements of  $\a$ with
``('' and elements of  $\b$ with ``)''.  The pairing is obtained by matching parentheses in the usual way on the infinite word  $\cdots w w w \cdots$.

\begin{example}
For  $P $ the total order on $\{1,\dots, 6\}$, the chains $\a=(6,4,3,2)$ and $\b=(5,4,1)$ produce a pairing of $4$ with $4$, $5$ with $3$, $1$ with $6$. Then after transferring $2$ we get chains $(6,4,3)$ and $(5,4,2,1)$.
The associated word of parentheses is  $)((())($ and after transferring 2 it becomes $))(())($.
\end{example}

Now we generalize the above construction of the $R$-matrix to the case of a general \threeone-free poset $P$.
Let  $H_m,\dots, H_1$ be the ladder decomposition of  $(\a,\b)$ (cf. Definition~\ref{d ladder stuff}).
Let
\begin{align}
  V^a = \{H_i  : |H_i \cap \a| \ge |H_i \cap \b| \},  \quad \quad
  V^b = \{H_i  : |H_i \cap \a| \le |H_i \cap \b| \}.
\end{align}
That is,  $V^a$ consists of  $\a$-ladders and balanced ladders, and  $V^b$ consists of
$\b$-ladders and balanced ladders.  We will often identify  $V^a$ with  the corresponding subset of the indices $\{1,\dots, m\}$, and likewise for  $V^b$.

We define $\eta(\a,\b)$ as follows.
First, we define the pairing between  $V^a$ and  $V^b$ exactly as above (thinking of these sets as subsets of integers). Suppose that $\a$ is longer than $\b$.
Note that any unpaired element of  $V^a$ must be an $\a$-ladder, not a balanced ladder.
Then $\eta(\a,\b)$ is obtained by
switching the unpaired $\a$-ladders into $\b$-ladders,
by moving their parts in $\a$ into $\b$, and vice versa
(c.f. Definition \ref{d ladder stuff}).
When $\b$ is longer than $\a$,
the pairing between $V^a$ and $V^b$
now leaves some unpaired  $\b$-ladders and $\eta(\a,\b)$ is obtained by switching these to $\a$-ladders.

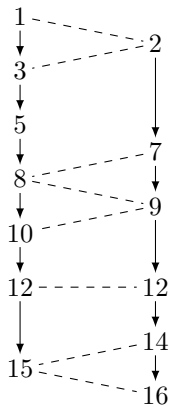
\begin{figure}[h]
\begin{tikzpicture}[xscale = 0.9, yscale = .36]
\tikzstyle{vertex}=[inner sep=0pt, outer sep=2.1pt]
\node[vertex] (v1) at (0,14) {\footnotesize$1$};
\node[vertex] (v2) at (2,13) {\footnotesize$2$};
\node[vertex] (v3) at (0,12) {\footnotesize$3$};
\node[vertex] (v4) at (0,10) {\footnotesize$5$};
\node[vertex] (v5) at (2,9) {\footnotesize$7$};
\node[vertex] (v6) at (0,8) {\footnotesize$8$};
\node[vertex] (v7) at (2,7) {\footnotesize$9$};
\node[vertex] (v8) at (0,6) {\footnotesize$10$};
\node[vertex] (v9) at (0,4) {\footnotesize$12$};
\node[vertex] (v10) at (2,4) {\footnotesize$12$};
\node[vertex] (v11) at (2,2) {\footnotesize$14$};
\node[vertex] (v12) at (0,1) {\footnotesize$15$};
\node[vertex] (v13) at (2,0) {\footnotesize$16$};

\draw[mystealthsmall] (v1) to (v3);
\draw[mystealthsmall] (v3) to (v4);
\draw[mystealthsmall] (v4) to (v6);
\draw[mystealthsmall] (v6) to (v8);
\draw[mystealthsmall] (v8) to (v9);
\draw[mystealthsmall] (v9) to (v12);

\draw[mystealthsmall] (v2) to (v5);
\draw[mystealthsmall] (v5) to (v7);
\draw[mystealthsmall] (v7) to (v10);
\draw[mystealthsmall] (v10) to (v11);
\draw[mystealthsmall] (v11) to (v13);

\draw[dashed] (v1) to (v2);
\draw[dashed] (v2) to (v3);
\draw[dashed] (v5) to (v6);
\draw[dashed] (v6) to (v7);
\draw[dashed] (v7) to (v8);
\draw[dashed] (v9) to (v10);
\draw[dashed] (v11) to (v12);
\draw[dashed] (v12) to (v13);
\end{tikzpicture}
\caption{\label{f ladder 1}
The ladder decomposition of two chains in  $\PP_2$, for the running example in \S\ref{ss proof of e commute}.}
\end{figure}

\begin{example}\label{example e chain ladders drawing}
Figure \ref{f ladder 1} shows two chains $(15,12,10,8,5,3,1)$ and $(16,14,12,9,7,2)$ in the poset  $\PP_{2}$
with ladder decomposition  $H_5 = (16,15,14)$, $H_4 = (12,12)$, $H_3 = (10,9,8,7)$, $H_2 = (5)$, $H_1 = (3,2,1)$.
The ladders $(5)$, $(3,2,1)$ are $\a$-ladders; $(12,12)$, $(10,9,8,7)$ are balanced ladders; $(16,15,14)$ is a $\b$-ladder.
The set
 $V^a = \{1,2,3,4\}$ consists of the  $\a$-ladders $H_1,H_2$ and balanced ladders
$H_3, H_4$, and  $V^b = \{3,4,5\}$ consists of the $\b$-ladder  $H_5$ and balanced ladders  $H_3,H_4$.
The pairing algorithm pairs $H_4$ to $H_4$, $H_3$ to $H_3$, $H_5$ to $H_2$.
\end{example}

\begin{theorem} \label{thm:eta}
 The resulting map $\eta(\a,\b) = (\c,\d)$ is an involution such that
\begin{align}
\label{e thm eta}
u_{a_k} \cdots u_{a_1} u_{b_{\ell}} \cdots u_{b_1} \equiv u_{c_{\ell}} \cdots u_{c_1} u_{d_{k}} \cdots u_{d_1}.
\end{align}
\end{theorem}

Associate a set of  intervals in $\mathbb{R}$ to the $V^a$, $V^b$ pairing as follows:
to each pair  $H_i \in V^a$ with $H_j \in V^b$ with  $i \le j$ (i.e. no wraparound),
we associate the closed interval  $[i,j]$;
to a pair  $H_i \in V^a$ with $H_j \in  V^b$ with  $i > j$ (a wraparound pair) we associate the two intervals
 $[i, m+1]$ and  $[0,j]$ (the numbers 0 and  $m+1$ are just needed to be smaller and larger than the index of any ladder); to each unpaired element  $H_i$ of  $V^a$ or $ V^b$
 we associate the single point $\{i\}$ in  $\mathbb{R}$.
We form {\it {lumps}} of ladders as connected components of the resulting system of intervals. Each unpaired ladder forms a lump by itself. Lumps inherit a total order from ladders.

\begin{example}
In the running example above, $H_5, H_4, H_3, H_2$ form a lump. Ladder $H_1$ by itself forms a lump,
which is smaller than the other lump.
\end{example}

\begin{figure}[h]
\centerfloat
\begin{tikzpicture}[xscale = 0.9, yscale = .36]
\tikzstyle{vertex}=[inner sep=0pt, outer sep=2.1pt]
\node[vertex] (v1) at (2,12) {\footnotesize$1$};
\node[vertex] (v2) at (0,11) {\footnotesize$3$};
\node[vertex] (v3) at (2,10) {\footnotesize$5$};
\node[vertex] (v4) at (0,9) {\footnotesize$6$};
\node[vertex] (v5) at (2,8) {\footnotesize$7$};
\node[vertex] (v6) at (0,7) {\footnotesize$9$};
\node[vertex] (v7) at (2,6) {\footnotesize$10$};
\node[vertex] (v8) at (0,5) {\footnotesize$11$};
\node[vertex] (v9) at (0,3) {\footnotesize$13$};
\node[vertex] (v10) at (0,1) {\footnotesize$15$};
\node[vertex] (v11) at (2,0) {\footnotesize$17$};

\draw[mystealthsmall] (v1) to (v3);
\draw[mystealthsmall] (v3) to (v5);
\draw[mystealthsmall] (v5) to (v7);
\draw[mystealthsmall] (v7) to (v11);

\draw[mystealthsmall] (v2) to (v4);
\draw[mystealthsmall] (v4) to (v6);
\draw[mystealthsmall] (v6) to (v8);
\draw[mystealthsmall] (v8) to (v9);
\draw[mystealthsmall] (v9) to (v10);

\draw[dashed] (v3) to (v4);
\draw[dashed] (v4) to (v5);
\draw[dashed] (v6) to (v7);
\draw[dashed] (v7) to (v8);
\end{tikzpicture}
\caption{\label{f ladder 2}
The ladder decomposition of two chains in  $\PP_2$, for Example \ref{ex ladder 2}.
}
\end{figure}
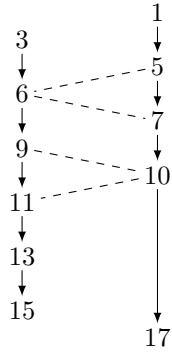

\begin{example}
\label{ex ladder 2}
Consider chains $\a = (15,13,11,9,6,3)$ and $\b=(17,10,7,5,1)$ in the poset $\PP_{2}$ (Figure \ref{f ladder 2}).
Their ladder decomposition is
\[\text{$H_7 = (17),$
$H_6 = (15),$
$H_5 = (13),$
$H_4 = (11,10,9),$
$H_3 = (7,6,5),$
$H_2 = (3),$
$H_1= (1)$}, \]
and $V^a = \{6,5,4,2\}$
and $V^b = \{7,3,1\}$.
As we run the pairing algorithm,
$H_7$ gets paired with $H_6$, $H_3$ with $H_2$,
and $H_1$ with $H_5$ (a wraparound pair).
As a result we get lumps $H_2 \cup H_3 = (7,6,5,3)$,
$H_4 = (11,10,9)$ (a single unpaired ladder lump),  $H_7 \cup H_6 \cup H_5 = (17,15,13)$, and $H_1 = (1)$.
\end{example}

For a lump $L$ let $L^a$ be the part of $L$ which is in $\a$, similarly $L^b$. Let $\mathbf{u}_{L^a}$ be the product of variables corresponding to elements of $L^a$, in order from largest to smallest. Similarly define $\mathbf{u}_{L^b}$. Note that if the lump $L$ consists of only one element, one of $L^a$ or $L^b$ will be empty. The original product can be rewritten as follows:
$$u_{a_k} \cdots u_{a_1} u_{b_{\ell}} \cdots u_{b_1} = \mathbf{u}_{L_r^a} \cdots \mathbf{u}_{L_1^a} \mathbf{u}_{L_r^b} \cdots \mathbf{u}_{L_1^b},$$ where $L_r, \ldots, L_1$ are the lumps in decreasing order.

The map $\eta$ then consists of switching the ladders that were not paired (and which thus form a lump by themselves).

Note that if there are wraparound pairs, then  $L_r^a$ contains all  $\a$-ladders paired via wraparound
and  $L_1^b$ contains all  $\b$-ladders paired via wraparound.

The following lemma is key to the proof of Theorem \ref{thm:eta}.

\begin{lemma} \label{lem:diag}
The $\a$ and the $\b$ parts of the lumps   can be moved next to each other modulo the ideal $\Iplacp{P}$:
\begin{align}
\label{e diag 1}
\mathbf{u}_{L_r^a} \cdots \mathbf{u}_{L_1^a} \mathbf{u}_{L_r^b} \cdots \mathbf{u}_{L_1^b} \equiv \mathbf{u}_{L_r^a} \mathbf{u}_{L_r^b} \cdots \mathbf{u}_{L_1^a} \mathbf{u}_{L_1^b}.
\end{align}
\end{lemma}


\begin{proof}
First assume $|\a| > |\b|$.
Consider the  $i$-th lump consisting of $L_i^a = (\alpha_{n'}, \ldots, \alpha_1)$ and $L_i^b = (\beta_n, \ldots, \beta_1)$.
Assume we have already moved the $\b$ parts of all larger lumps next to their $\a$ parts.
It suffices to show that
we can move $L_i^b$ to the left past $\e{v} = L_{i-1}^a \cdots L_1^a$ and next to  $L_i^a$ using the relations of  $\Iplacp{P}$, i.e.,
\begin{align}
\label{e lump move}
\mathbf{u}_{L_i^a} \mathbf{u}_\e{v} \mathbf{u}_{L_i^b} \equiv \mathbf{u}_{L_i^a} \mathbf{u}_{L_i^b} \mathbf{u}_\e{v}  \ \bmod \, \Iplacp{P}.
\end{align}
There are four cases.
The first case is that the lump  $L_i$ is not a single ladder and does not contain a wraparound pair.
Because of the way ladders are paired and lumped,
as we scan the ladders of  $(V^a \cap L_i) \sqcup (V^b \cap L_i)$ from largest to smallest we have always encountered more elements of  $V^b$ than  $V^a$ (except at the beginning and end where we've seen an equal number from both).
It follows that  $n' =n$, $\alpha_{i} \to_P \beta_{i+1}$, and  $\alpha_i \dotto_P \beta_i$, and so
Lemma~\ref{l shuffle chains} proves \eqref{e lump move}.

The second case is that  $L_i = L_r$ contains a wraparound pair.  This lump now contains more  $\a$-ladders than
$\b$-ladders and  $n' > n$.  However we can just ignore the first  $n'-n$ letters of  $L_r^a$ and do the same thing we did in the previous case; we still have $\alpha_{i} \to_P \beta_{i+1}$ and  $\alpha_i \dotto_P \beta_i$,
and so Lemma \ref{l shuffle chains} gives \eqref{e lump move}.

The third case is that the lump  $L_i$ consists of a single balanced ladder.
Then Lemma~\ref{l shuffle chains 2} (resp. Lemma~\ref{l shuffle chains}) gives
\eqref{e lump move} if  $L_i$ has the form (iv) or (v) (resp. (iii)) in Figure \ref{fig:ladder}.

The fourth case is that  $L_i$ consists of a single  $\a$-ladder.
Then  $n' = n+1$ and we can just ignore the first letter $\alpha_{n+1}$ of  $L_i^a$ and obtain \eqref{e lump move} by
applying Lemma \ref{l shuffle chains} with  $\e{a}= \alpha_{n} \cdots \alpha_1$.

The proof when $|\a| < |\b|$ is similar, keeping in mind Remark \ref{r flip}.
We move $L^a_i$ to the right through  $L^b_r \cdots L^b_{i+1}$ to $L^b_i$, starting with  $i = 1$.
\end{proof}

\begin{example}
This will show why, when  $|\a| < |\b|$, it
is necessary to move $L^a_i$ to the right through $L^b_r \cdots L^b_{i+1}$ to $L^b_i$,
using $L^b_i$ to help with the relations (as in the flip of Lemma \ref{l shuffle chains} or \ref{l shuffle chains 2}), rather than going the other direction.
Consider chains $\a=(5,2)$ and $\b=(6,4,1)$ in the poset  $\PP_{2}$. Their ladder decomposition is
$H_2 = (6,5,4),$ $H_1 =(2,1)$. The pairing process creates only the single pair of  $H_1$ with  $H_1$, and  $H_1$ and $H_2$ are each lumps by themselves.
Since  $|\a| < |\b|$, the proof above says that to obtain
$u_{5}u_2u_6u_4u_1 \equiv u_{5}u_6u_4u_2u_1$,
we move $u_2$ through $u_6u_4$ to $u_1$:
$u_2u_6u_4u_1 \equiv u_6u_2u_4u_1 \equiv u_6u_4u_2u_1.$
On the other hand,
one can check that
$u_{5}u_2u_6u_4 \not \equiv u_{5}u_6u_4u_2$
modulo  $\Iplacp{P}$.
\end{example}

\begin{proof}[Proof of Theorem \ref{thm:eta}]
Since  $(\a,\b)$ and  $(\c,\d)$ play symmetric roles,
we can assume $|\a| > |\b|.$
To go from the left side of \eqref{e thm eta} to the right,
first we shuffle the lumps  using (the $|\a| > |\b|$ case of) Lemma \ref{lem:diag}, switch the single $\a$-ladder lumps
using Lemma \ref{l swap chains}, and then we unshuffle the lumps using
(the $|\a| < |\b|$ case of) Lemma \ref{lem:diag}.
\end{proof}

\begin{example}
In our running example we have
\begin{align}
&\mathbf{u}_{15 \ 12 \ 10 \ 8 \ 5 \ 3 \ 1   \  16 \ 14 \ 12 \ 9 \ 7 \ 2} \equiv \\
&\mathbf{u}_{15 \ 12 \ 10 \ 8 \ 5  \  16 \ 14 \ 12 \ 9 \ 7 \   3 \  1  \  2} \equiv \\
&\mathbf{u}_{ 15 \ 12 \ 10 \ 8 \ 5  \  16 \ 14 \ 12 \ 9 \ 7 \   2 \  3 \  1} \equiv \\
&\mathbf{u}_{ 15 \ 12 \ 10 \ 8 \ 5 \ 2  \  16 \ 14 \ 12 \ 9 \ 7 \      3 \  1}. \phantom{\equiv}
\end{align}
Here there are two lumps. The first step shuffles the lumps, the second step switches the single ladder lump from being an $\a$-lump to being a $\b$-lump, and the third step unshuffles the lumps.
\end{example}

\section{Proof of Theorem \ref{th:IplacP-positivity}}
\label{sec:proof-of-IplacP-positivity}

We will prove Theorem \ref{th:IplacP-positivity} by an inductive
computation of $\mathfrak{J}^P_\lambda(\mathbf{u})$
which  involves peeling off the diagonals of the shape~$\lambda$
and incorporating them into a tableau.
This proof follows a similar format to that of \cite[Theorem 2.14]{BF} and \cite[Theorem 4.8]{BLamLLT}
but the technical details are quite different;
it can be viewed as a noncommutative version of a lattice path canceling argument.
The proof requires the following flagged generalization of the  $\mathfrak{J}^P_{\lambda}(\mathbf{u})$.

\begin{definition}
For a poset  $P$, 
$\alpha\!=\!(\alpha_1,\dots,\alpha_\ell)\!\in\!\ZZ_{\ge 0}^\ell$,
and an $\ell$-tuple
$\mathbf{Z}\!=\!(Z_1,\dots,Z_\ell)$ of subsets of $P$,
we define the \emph{noncommutative column-flagged  $P$-Schur function}
$J^P_{\alpha}(\mathbf{Z})$ by
\begin{align}
J^P_{\alpha}(\mathbf{Z})
=\sum_{\pi\in \SS_{\ell}}
\sgn(\pi) \, e_{\alpha_1+\pi(1)-1}^P(\mathbf{u}_{Z_1})
\,e_{\alpha_2+\pi(2)-2}^P(\mathbf{u}_{Z_2})\cdots
e_{\alpha_{\ell}+\pi(\ell)-\ell}^P(\mathbf{u}_{Z_\ell}). \label{e flag schur}
\end{align}
Here,
$e_k^P(\mathbf{u}_S)$ denotes the noncommutative $P$-elementary function in the variables $\{u_a\}_{a \in S}$ (cf. \eqref{e ek def}); by convention,
$e_0(\mathbf{u}_S)=1$ (even when  $S = \varnothing$) and $e_k^P(\mathbf{u}_S) = 0$ for $k<0$ or $k>|S|$.

Note that the
$\mathfrak{J}^P_{\lambda}(\mathbf{u})$ are recovered as the special case
$\mathfrak{J}^P_{\lambda}(\mathbf{u}) = J^P_{\lambda'}(P,\ldots,P)$.
\end{definition}

\begin{definition}
\label{d diagread}
Let $\mathbf{\alpha}=(\alpha_1,\ldots,\alpha_\ell)\in \ZZ_{\ge 0}^\ell$
be a weak composition satisfying $\alpha_{j+1}\le\alpha_j+1$ for all~$j$.
Write $\alpha'$ for the
diagram whose  $j$th column contains the boxes in rows $1,\dots,\alpha_j$, that is,
$\alpha'=\{(i,j) : j\in [\ell], i\in [\alpha_j]\}\subset
\ZZ_{\ge 1} \times \ZZ_{\ge 1}$, drawn in
English notation as in \S\ref{ss Schur pos P Knuth}.

For an  $\ell$-tuple $\mathbf{Z}=(Z_1,\ldots,Z_\ell)$ of subsets of  $P$,
let $\SSYT_P^{\mathbf{Z}}(\alpha')$ denote the set of
fillings
$T: \ZZ_{\ge 1} \times \ZZ_{\ge 1} \to P$ of the diagram~$\alpha'$ such that
\begin{itemize}
\item Each column of  $T$, read bottom to top, is  $P$-strictly decreasing, i.e.  $T(\alpha_j,j) \leftto_P   \cdots \leftto_P T(2,j) \leftto_P T(1,j)$;
\item  $T(r,c) \dotto_P T(r,c+1)$ whenever the box $(r,c)$ and the box $(r,c+1)$ immediately to its east belong to the diagram $\alpha'$;
\item the entries in column  $j$ lie in $Z_j$;
\item if $\alpha_j < \alpha_{j+1}$, then $T(\alpha_{j+1},j+1) \notin  Z_j$.
\end{itemize}
\end{definition}

\begin{example}
Let $P = \PP_{2}$,
$\alpha = (3,1,2,3)$, and  $\mathbf{Z} = ([6],[3],[6],[7])$.  Then
\begin{align}
\label{eq: PSSYT cols}
{\fontsize{8pt}{6pt}\selectfont \tableau{1&3&2&1 \\3& \bl  & 4 & 3 \\ 5 & \bl &\bl & 7}}
\ \in \, \SSYT_P^\mathbf{Z}(\alpha').
\end{align}
Note that
for any $T \in \SSYT_P^\mathbf{Z}(\alpha')$, the box  $(3,4)$ must be filled with a 7 (the unique element of  $Z_4 \setminus Z_3$) because of the last bullet point above.
\end{example}

The \emph{diagonal reading word} of $T \in \SSYT_P^\mathbf{Z}(\alpha')$, denoted $\diagread(T)$,
is obtained by concatenating the diagonals of~$T$
(reading each diagonal in the southeast  direction), starting with the
southwesternmost diagonal of~$\alpha'$.
For example, the diagonal reading word of the tableau in \eqref{eq: PSSYT cols} is  $\e{5\spa 3\spa 1 \spa 3 47 \spa 23\spa 1}$.

We will establish Theorem~\ref{th:IplacP-positivity} by proving the following
flagged generalization; Theorem ~\ref{th:IplacP-positivity} is obtained
from it by setting  $Z_1 = \cdots = Z_\ell = P$ and $\alpha = \lambda'$ and
then
applying a result about reading words,
Proposition \ref{pr: diag vs col}, which we
 postpone to the end of the section.

\begin{theorem}
\label{th:IplacP-positivity-technical}
Let $P$ be a \threeone-free poset. Let  $\alpha=(\alpha_1,\dots, \alpha_\ell)\in\ZZ_{\ge 0}^\ell$ be a partition,
and let $\mathbf{Z} = (Z_1,\ldots, Z_\ell)$ be an  $\ell$-tuple of lower order ideals of  $P$ satisfying  $Z_1 \subset Z_2 \subset \cdots \subset Z_\ell$.
Then
\begin{align}\label{et:FG'-positivity-technical}
J^P_\alpha(\mathbf{Z})
 \equiv
\sum_{ T \in \SSYT_P^{\mathbf{Z}}(\alpha')} \mathbf{u}_{\diagread(T)} \ \, \bmod \ \Iplacp{P}.
\end{align}
\end{theorem}

We will prove Theorem~\ref{th:IplacP-positivity-technical} by
establishing its more technical variant,
Theorem~\ref{t new statement fgprime} below.
We will need an augmented version of the  $J_\alpha^P(\mathbf{Z})$ and three lemmas.

\begin{definition}
\label{def:Jalpha-augmented}
For a weak composition $\alpha\!=\!(\alpha_1,\dots,\alpha_\ell)\!\in\!\ZZ_{\ge
  0}^\ell$,
an $\ell$-tuple $\mathbf{Z}\!=\!(Z_1,\dots,Z_\ell)$ of subsets of  $P$,
and words
$\e{w}^1, \ldots, \e{w}^{\ell-1} \in \U_P^*$,
define
\begin{align}
&J^P_{\alpha}(Z_1\Jnot{\e{w}^1}Z_2\Jnot{\e{w}^2}\cdots\Jnot{\e{w}^{\ell-1}}Z_\ell)  \label{e
    augmented J1} \\
&=\sum_{\pi\in \SS_{\ell}}
\sgn(\pi) \ e_{\alpha_1+\pi(1)-1}^P(\mathbf{u}_{Z_1}) \, \mathbf{u}_{\e{w}^1} \,
e_{\alpha_2+\pi(2)-2}^P(\mathbf{u}_{Z_2})\, \mathbf{u}_{\e{w}^2}\cdots
\mathbf{u}_{\e{w}^{\ell-1}} \, e_{\alpha_{\ell}+\pi(\ell)-\ell}^P(\mathbf{u}_{Z_\ell}). \label{e
  augmented J2}
\end{align}
We omit  $\e{w}^j$ from the notation in \eqref{e augmented J1} if the word $\e{w}^j$ is empty.
\end{definition}

\begin{lemma}
\label{l plactic inc word}
Let $b$ be an element of  $P$ and $\e{v} \!= \! \e{v}_1 \cdots \e{v}_s$ and $\e{w}  \!= \! \e{w}_1 \cdots \e{w}_t$ be words in  $\U^*_P$  satisfying
\begin{itemize}
\item[(i)] $b \leftto_P v_1 \leftto_P \cdots \leftto_P v_s$,
\item[(ii)]  $b \dotto_P w_1 \dotto_P w_2 \dotto_P \cdots \dotto_P w_t$,
\item[(iii)] $v_i \to_P w_j$ for all  $i$ and  $j$.
\end{itemize}
 Then $u_b \mathbf{u}_\e{v} \mathbf{u}_\e{w} \equiv u_b \mathbf{u}_\e{w} \mathbf{u}_\e{v}  \bmod \Iplacp{P}$.
\end{lemma}

\begin{proof}
Repeatedly apply Lemma \ref{l plactic inc word 0} to move each letter of  $\e{w}$ past  $\e{v}$, starting with $w_1$:
$u_b \mathbf{u}_\e{v} u_{w_1}\cdots u_{w_t} \equiv  u_b  u_{w_1} \mathbf{u}_\e{v}
u_{w_2}\cdots u_{w_t} \equiv $
$u_b u_{w_1} u_{w_2} \mathbf{u}_\e{v} u_{w_3}\cdots
u_{w_t} \equiv \cdots \equiv  u_b \mathbf{u}_\e{w} \mathbf{u}_\e{v}   \bmod   \Iplacp{P}$.
\end{proof}

\begin{lemma}
\label{l elem sym J0}
If   $\alpha_{j+1}=\alpha_j+1 $ and  $Z_{j}=Z_{j+1}$, then
$J^P_\alpha(\mathbf{Z}) \equiv 0 \bmod \Iplacp{P}\,$.
More generally, 
this holds in the augmented case provided $\e{w}^j$ is empty.
\end{lemma}

\begin{proof}
This follows from the fact that the noncommutative $P$-elementary symmetric functions commute modulo~$\Iplacp{P}$.
\end{proof}

For $b \in P$, let $(b)_{<_P} = \{a \in P : a \to_P b \}$, the set of elements less than  $b$ in  $P$.

\begin{lemma}
Let $Z \subset P$,  $z \in Z$, and $k \in \ZZ$. Then
\begin{equation}
\label{e ek induction}
e_k^P(\mathbf{u}_{Z}) =u_{z} \spa e_{k-1}^P(\mathbf{u}_{(z)_{<_P}})+e_k^P(\mathbf{u}_{Z \setminus \{z\}}).
\end{equation}
\end{lemma}

\begin{proof}
This is immediate from the definition of $e_k^P(\mathbf{u}_{Z})$.
\end{proof}

For an index  $j \in [\ell]$ and  $z \in Z_j$, we will apply \eqref{e ek induction} to  $J^P_\alpha(Z_1, \ldots, Z_\ell)$
and its variants by expanding
$e_{\alpha_j+\pi(j)-j}^P(\mathbf{u}_{Z_j})$ in \eqref{e flag schur}
using \eqref{e ek induction} (so that \eqref{e ek induction} is
applied once to each of the  $\ell!$ terms in the sum in \eqref{e flag
  schur}).
We refer to this as a \emph{$(j, z)$-expansion of  $J^P_\alpha(Z_1,\ldots, Z_\ell)$
using~\eqref{e ek induction}}.

For $\alpha\!=\!(\alpha_1,\dots, \alpha_\ell)\in\ZZ_{\ge 0}^\ell$, let
$\Asc(\alpha) = \{j \in [\ell] : \alpha_j < \alpha_{j+1} \}$ and
 $\overline{\Asc}(\alpha) = [\ell] \setminus \Asc(\alpha)$, where by convention
 $\alpha_{\ell+1}=0$ so that  $\ell \in \overline{\Asc}(\alpha)$ always.
The \emph{peeling index} of $\alpha$ is the minimum of
$\overline{\Asc}(\alpha)$.

For a word  $\e{w}$ and $T \in \SSYT_P^\mathbf{Z}(\alpha')$,
define the word $\diagread(\e{w};T)$ as follows:
if  $\alpha_1 = \alpha_2 = 0$, set $\diagread(\e{w};T) = \e{w}\,\diagread(T)$.
 Otherwise, $\diagread(\e{w};T) = \e{d^1 \spa w \spa d^2  \cdots d^t}$, where
$\e{d^1, d^2, \ldots, d^t}$ are the words
 corresponding to the (nonempty) diagonals of~$T$ read in the southeast direction,
numbered starting with the southwesternmost diagonal.
(Hence $\diagread(\e{w};T) = \diagread(T)$ when  $\e{w}$ is empty.)

\begin{theorem}
\label{t new statement fgprime}
Let $P$ be a \threeone-free poset. Let  $\alpha=(\alpha_1,\dots, \alpha_\ell)\in\ZZ_{\ge 0}^\ell$ be such that
$\alpha_{i+1}\le \alpha_i+1$ for all~$i$.
Let $j$ be the peeling index of~$\alpha$.
Let  $\mathbf{Z}= (Z_1, \dots, Z_\ell)$ be an  $\ell$-tuple of subsets  $Z_i \subset P$  and
$\e{w} = \e{w_1\cdots w_t}$ be a word in $\U_P^*$ satisfying
\begin{itemize}
\item[(a)]  each  $Z_i$ is a lower order ideal of  $P$ (which may be the empty set).
\item[(b)] for $i \in \overline{\Asc}(\alpha)$,  $Z_{i} \setminus Z_{i+1}$ contains no comparable pair $a \to b$.
\item[(c)] for  $i \in \Asc(\alpha)$,  $Z_i \subset Z_{i+1}$.
\item[(d)] for  $i \in \Asc(\alpha)$, there is no induced $\twoone$ subposet
\raisebox{-0.7mm}{\begin{tikzpicture}[xscale = 1.84,yscale = 1.4]
\tikzstyle{vertex}=[inner sep=0pt, outer sep= 2.5pt]
\tikzstyle{aedge} = [draw, ->,>=stealth', black]
\tikzstyle{aedgecurve} = [draw, ->,>=stealth', black, bend left=40]
\tikzstyle{dashedcurve} = [draw, dashed, black, bend left=40]
\tikzstyle{edge} = [draw, thick, -,black]
\tikzstyle{dashededge} = [draw, -,  dashed, black]
\node[vertex] (va) at (0,-0.8){${a}$};
\node[vertex] (vb) at (0.5,-0.8){${b}$};
\node[vertex] (vc) at (1,-0.8){${c}$};
\draw[aedge] (va) to (vb);
\draw[dashededge] (vb) to (vc);
\draw[dashedcurve] (va) to (vc);
\end{tikzpicture}}
of $P$ with $a,b \in Z_i$ and $c \in P\setminus Z_i$.
\item[(e)] $a \dotto \e{w}_k$ for all  $a \in Z_j$ and  $k \in [t]$.
\item[(f)] $\e{w_1} \dotto \e{w_2} \dotto \cdots \dotto \e{w_t}$.
\item[(g)] there is no induced $\twoone$ subposet
\raisebox{-1mm}{\begin{tikzpicture}[xscale = 1.84,yscale = 1.4]
\tikzstyle{vertex}=[inner sep=0pt, outer sep= 2.5pt]
\tikzstyle{aedge} = [draw, ->,>=stealth', black]
\tikzstyle{aedgecurve} = [draw, ->,>=stealth', black, bend left=40]
\tikzstyle{dashedcurve} = [draw, dashed, black, bend left=40]
\tikzstyle{edge} = [draw, thick, -,black]
\tikzstyle{dashededge} = [draw, -,  dashed, black]
\node[vertex] (va) at (0,-0.8){${a}$};
\node[vertex] (vb) at (0.5,-0.8){${b}$};
\node[vertex] (vc) at (1,-0.8){${c}$};
\draw[aedge] (va) to (vb);
\draw[dashededge] (vb) to (vc);
\draw[dashedcurve] (va) to (vc);
\end{tikzpicture}}
of  $P$ with $a,b \in Z_j$ and $c~\in~\{w_1,\dots, w_t\}$.
\end{itemize}
Then
\begin{align}
\label{e claim diagread}
 J^P_{\alpha}(Z_1, \dots, Z_{j-1}, Z_j \Jnot{\e{w}}Z_{j+1},\dots) \equiv
 \sum_{ T \in \SSYT_P^\mathbf{Z}(\alpha')} \mathbf{u}_{\diagread(\e{w};T)}
 \bmod \Iplacp{P}.
\end{align}
\end{theorem}

The reader may wish to follow along the proof with Example \ref{ex main proof}. Figures 1 and 2 of \cite{BLamLLT} may also be helpful (they apply to a different technical setup, but the idea is similar).

\begin{proof}
Write $J$ for $J^P_{\alpha}(Z_1,\dots,Z_j\Jnot{\e{w}}Z_{j+1},\dots)$.
The proof is by induction on  $\ell$ and the~$Z_i$.

We first address the base cases $(\alpha_1 = 0, j=1)$ and $(Z_1 = \varnothing, j=1)$.
Suppose $\alpha_1 = 0$ and $j=1$.
Since $\alpha_{i+1} \leq \alpha_{i}+1$ for all  $i$
and $\alpha_1 = \alpha_2 = 0$, $J^P_{\alpha}(Z_1,\dots, Z_\ell)$ is  a noncommutative version of
the determinant of a matrix whose first column, read downwards, is a 1 followed by~0's.
Then by Definition~\ref{def:Jalpha-augmented} and induction,
\begin{align}\label{e j1 base case}
J= \mathbf{u}_{\e{w}}J^P_{\hat{\alpha}}(\mathbf{\hat{Z}}) \equiv
\mathbf{u}_{\e{w}}\sum_{ T \in \SSYT_P^\mathbf{\hat{Z}}(\hat{\alpha}')  } \mathbf{u}_{\diagread(T)}
 = \sum_{ T \in \SSYT_P^\mathbf{Z}(\alpha')} \mathbf{u}_{\diagread(\e{w};T)} \ \bmod \Iplacp{P},
\end{align}
where  $\hat{\alpha} = (\alpha_2,\dots, \alpha_\ell)$ and $\mathbf{\hat{Z}} = (Z_2,\dots, Z_\ell)$.
Now suppose $Z_1 = \varnothing$ and  $j=1$.
By the previous case, we may assume  $\alpha_1 > 0$.
Then  $J^P_{\alpha}(Z_1,\dots, Z_\ell)$ is a noncommutative version of the
determinant of a matrix whose first row is~0, hence $J=0$.
Since  $Z_1 =\varnothing$ implies that $\SSYT_P^\mathbf{Z}(\alpha')$ is empty, the right side of \eqref{e claim diagread} is also 0.

By the base cases above, we may assume that either ($\alpha_1 > 0$ and  $Z_1 \ne \varnothing$ and  $j=1$) or  $j >1$.

\textbf{Case 1: $\alpha_1 > 0$ and  $Z_1 \ne \varnothing$ and  $j=1$.}
Let  $z_1$ be a maximal element of $Z_1$.
Set
\[\text{
$\alpha_- = (\alpha_1-1,\alpha_{2}, \dots)$, \ \,
$\mathbf{Z}_* = ((z_1)_{<_P},Z_{2},Z_3,\ldots)$, \ \,  and \ \
$\mathbf{Z}_- = (Z_1\setminus \{z_1\},Z_{2},Z_3,\ldots)$.}\]
Note that  $(z_1)_{<_P} = (z_1)_{<_P} \cap Z_1$ since  $Z_1$ is a lower order ideal.
A $(1,z_1)$-expansion of $J$ using \eqref{e ek induction} yields
\begin{align}
J&=u_{z_1}J^P_{\alpha_-}((z_1)_{<_P}\Jnot{\e{w}}Z_2,\dots, Z_\ell)
+J^P_{\alpha}(Z_1\setminus{z_1} \Jnot{\e{w}}Z_2,\dots, Z_\ell) \notag 
\\
\label{e un2}
&\equiv
u_{z_1}\mathbf{u}_{\e{w}}J^P_{\alpha_-}((z_1)_{<_P},Z_2,\dots, Z_\ell)
+J^P_{\alpha}(Z_1\setminus{z_1}\Jnot{\e{w}}Z_2,\dots, Z_\ell),
\end{align}
where the congruence is modulo  $\Iplacp{P}$ and holds by Lemma~\ref{l plactic inc word} applied with  $b=z_1$ and  $\e{v}$ any  $P$-strictly decreasing word contributing to the leftmost  $e^P$ in the definition \eqref{e augmented J2} of  $J^P_{\alpha_-}((z_1)_{<_P}\Jnot{\e{w}}Z_2,\dots, Z_\ell)$ (in particular, the letters of  $\e{v}$ belong to $(z_1)_{<_P}$). Assumption (i) of the lemma holds by definition of  $e_k^P$, assumption (ii) holds by (e) and (f),
while (iii) states that $v_i \to w_k$ for any letter  $v_i$ in the word $\e{v}$ and $w_k$ in $\e{w}$.
We know  $v_i \to z_1$, and $v_i \dotto w_k$ and  $z_1 \dotto w_k$ by (e). By (g), we cannot have both
$v_i \dote w_k$ and  $z_1 \dote w_k$, so we must have  $v_i \to w_k$ or  $z_1 \to w_k$ but then  $v_i \to w_k$ always by transitivity.

We will apply induction to both terms in \eqref{e un2} and we need to show that the conditions (a)--(g) hold for both terms.
We now address the first term  $J^P_{\alpha_-}((z_1)_{<_P},Z_2,\dots, Z_\ell) = J^P_{\alpha_-}(Z_*)$.
Since the input word is empty for this term, (e)--(g) hold trivially.
Since  $(Z_*)_i = Z_i$ for  $i > 1$, we know (a)--(d) hold for  $i>1$ and it remains to check them for  $i=1$.
Condition (a) holds because  $z_1$ is a maximal element of  $Z_1$.
Since  $Z_1 \setminus Z_2$ contains no pair of comparable elements, $(z_1)_{<_P} \subset Z_2$, which proves (b) and (c) for  $i=1$
(we don't know whether or not 1 lies in $\Asc(\alpha_-)$ but the stronger condition in (c) holds either way).
For (d), we must show that  $P$ has no induced $\twoone$ subposet
\raisebox{-0.7mm}{\begin{tikzpicture}[xscale = 1.84,yscale = 1.4]
\tikzstyle{vertex}=[inner sep=0pt, outer sep= 2.5pt]
\tikzstyle{aedge} = [draw, ->,>=stealth', black]
\tikzstyle{aedgecurve} = [draw, ->,>=stealth', black, bend left=40]
\tikzstyle{dashedcurve} = [draw, dashed, black, bend left=40]
\tikzstyle{edge} = [draw, thick, -,black]
\tikzstyle{dashededge} = [draw, -,  dashed, black]
\node[vertex] (va) at (0,-0.8){${a}$};
\node[vertex] (vb) at (0.5,-0.8){${b}$};
\node[vertex] (vc) at (1,-0.8){${c}$};
\draw[aedge] (va) to (vb);
\draw[dashededge] (vb) to (vc);
\draw[dashedcurve] (va) to (vc);
\end{tikzpicture}}
with $a, b \in (z_1)_{<_P}$ and  $c \notin (z_1)_{<_P}$.
Suppose it did. Then each possible relationship between  $z_1$ and  $c$ is impossible:  $z_1 \leftarrow c$ is impossible since $c \notin (z_1)_{<_P}$,
$z_1 \to c$ is impossible since then $b \to z_1 \to c$ would contradict  $b \dote c$,
and $c \dote z_1$ is impossible since this would make the chain $\{a,b,z_1\}$ all incomparable to  $c$ contradicting  $P$ being \threeone-free.

The conditions (a)--(g) for the second term $J^P_{\alpha}(Z_1\setminus{z_1}\Jnot{\e{w}}Z_2,\dots, Z_\ell)$
follow easily from (a)--(g) for the original  $J$:  (b), (e), and (g) are weaker statements for this term than for the original  $J$ and the remaining statements are identical (note $1 \in \overline{\Asc}(\alpha)$ so for  $i=1$ we only need check (b) and not  (c) nor (d)).

Hence by induction, \eqref{e un2} is congruent modulo  $\Iplacp{P}$ to
\begin{align}
& ~u_{z_1}\mathbf{u}_{\e{w}}\sum_{ T \in \SSYT_P^{\mathbf{Z}_*}(\alpha_-')}  \mathbf{u}_{\diagread(T)}
+ \sum_{ T \in \SSYT_P^{\mathbf{Z}_-}(\alpha')}  \mathbf{u}_{\diagread(\e{w};T)}.
\label{e diag rec0}
\end{align}
Now \eqref{e diag rec0} is equal to the right side of
\eqref{e claim diagread} because \eqref{e diag rec0} is the
result of partitioning the set of tableaux appearing in
\eqref{e claim diagread} into two, depending on whether or not  $T(\alpha_1,1)$
 is equal to $z_1$.
Note that if  $\alpha_1 = \alpha_2$, then for $T \in \SSYT_P^{\mathbf{Z}}(\alpha')$ with $T(\alpha_1,1) = z_1$,
the condition $z_1=T(\alpha_1,1) \dotto T(\alpha_2,2)$ (which holds since the rows of  $T$ are $P$-weakly increasing)
is equivalent to the condition  $U((\alpha_-)_{2},2) \notin (z_1)_{<_P}$ required of any $U \in \SSYT_P^{\mathbf{Z}_*}(\alpha'_-)$.

\textbf{Case 2: $j > 1$.}
By definition of the peeling index,
 $\alpha_{j}=\alpha_{j-1}+1$ and
$\alpha_j\geq \alpha_{j+1}$.
If $Z_{j-1}=Z_{j}$, then
$J\equiv 0 \bmod \Iplacp{P}$ by Lemma~\ref{l elem sym J0};
in this case $\SSYT_P^\mathbf{Z}(\alpha')$ is empty, so the right side of \eqref{e claim diagread} is also 0.
So we may assume that $Z_{j-1} \subsetneq Z_j$.
Choose a maximal element  $z_j$ of  $Z_j$ which does not belong to  $Z_{j-1}$ (such an element exists because  $Z_{j-1}$ is a lower order ideal).
Set
\begin{center}
\text{$\mathbf{Z}_* = (Z_1,\ldots,Z_{j-1},(z_j)_{<_P},Z_{j+1},\ldots)$,  \ \ \
$\mathbf{Z}_- = (Z_1,\ldots,Z_{j-1},Z_j\setminus \{z_j\},Z_{j+1},\ldots)$,} \\[2mm]
\text{$\alpha_- = (\alpha_1,\dots, \alpha_{j-1}, \alpha_j-1, \alpha_{j+1},\dots)$.}
\end{center}
Now by a  $(j, z_j)$-expansion of  $J$ using \eqref{e ek induction}, we have
\begin{align}
J=&~J^P_{\alpha_-}(Z_1,\dots, Z_{j-1} \Jnot{\e{z_j}} (z_j)_{<_P} \Jnot{\e{w}} Z_{j+1},\dots) \notag\\
&+J^P_{\alpha}(Z_1, \dots, Z_{j-1}, Z_j \setminus \{z_j\} \Jnot{\e{w}}Z_{j+1},\dots)\notag\\
\equiv &~J^P_{\alpha_-}(Z_1,\dots, Z_{j-1}\Jnot{\e{z_jw}}(z_j)_{<_P},Z_{j+1},\dots) \label{e unj3}\\
&~+J^P_{\alpha}(Z_1,\dots, Z_{j-1}, Z_j \setminus \{z_j\}\Jnot{\e{w}}Z_{j+1},\dots),
\label{e unj4}
\end{align}
where the congruence is by Lemma~\ref{l plactic inc word} and the assumptions of the lemma are satisfied by the same argument given in Case 1.

We will apply induction to the terms in \eqref{e unj3}--\eqref{e unj4}.
First consider the term in \eqref{e unj3}.
Note that $j-1$ is the peeling index of $\alpha_-$.
We must show that the flagging $\mathbf{Z}_*$ and the word $\e{z_j w}$ satisfy (a)--(g).
Condition (a) holds since  $(z_j)_{<_P}$  is a lower order ideal.
For (b), we need to show  $Z_{j-1} \setminus (z_j)_{<_P}$ contains no pair  $a \to b$. Suppose that it does.  Then  $a \dote z_j$ and  $b \dote z_j$ which follows from  $a,b \notin (z_j)_{<_P}$ and  $z_j \notin Z_{j-1}$ and  $Z_{j-1}$ being a lower order ideal; but then the induced $\twoone$ subposet
\raisebox{-1mm}{\begin{tikzpicture}[xscale = 1.9,yscale = 1.4]
\tikzstyle{vertex}=[inner sep=0pt, outer sep= 2.5pt]
\tikzstyle{aedge} = [draw, ->,>=stealth', black]
\tikzstyle{aedgecurve} = [draw, ->,>=stealth', black, bend left=40]
\tikzstyle{dashedcurve} = [draw, dashed, black, bend left=40]
\tikzstyle{edge} = [draw, thick, -,black]
\tikzstyle{dashededge} = [draw, -,  dashed, black]
\node[vertex] (va) at (0,-0.8){${a}$};
\node[vertex] (vb) at (0.5,-0.8){${b}$};
\node[vertex] (vc) at (1,-0.8){${z_j}$};
\draw[aedge] (va) to (vb);
\draw[dashededge] (vb) to (vc);
\draw[dashedcurve] (va) to (vc);
\end{tikzpicture}}
contradicts condition (d) (with $i = j-1$) for the original  $J$.
For (c), we have $(z_j)_{<_P} \subset Z_{j+1}$ since  $Z_{j} \setminus Z_{j+1}$
contains no pair of comparable elements (by condition (b) for the original  $J$).
We need to check (d) for $i = j$, which says that there is
 no induced $\twoone$ subposet
 \raisebox{-1mm}{\begin{tikzpicture}[xscale = 1.9,yscale = 1.4]
\tikzstyle{vertex}=[inner sep=0pt, outer sep= 2.5pt]
\tikzstyle{aedge} = [draw, ->,>=stealth', black]
\tikzstyle{aedgecurve} = [draw, ->,>=stealth', black, bend left=40]
\tikzstyle{dashedcurve} = [draw, dashed, black, bend left=40]
\tikzstyle{edge} = [draw, thick, -,black]
\tikzstyle{dashededge} = [draw, -,  dashed, black]
\node[vertex] (va) at (0,-0.8){${a}$};
\node[vertex] (vb) at (0.5,-0.8){${b}$};
\node[vertex] (vc) at (1,-0.8){${c}$};
\draw[aedge] (va) to (vb);
\draw[dashededge] (vb) to (vc);
\draw[dashedcurve] (va) to (vc);
\end{tikzpicture}}
with  $a, b \in (z_j)_{<_P}$ and
$c \notin (z_j)_{<_P}$ (this is only required when  $\alpha_j = \alpha_{j+1}$ but it holds whether or not this is true).  This is proved by the same argument given in Case 1.
For  (e), since  $F_{j-1} \subset F_j$, it suffices to check  $a \dotto z_j$ for every  $a \in Z_{j-1}$; this follows from  $Z_{j-1} \subset Z_j$ and $z_j$ maximal in  $Z_j$.  For  (f) we need  $z_j \dotto w_1$, which holds by (e) for the original  $J$.
Finally, (g)
follows from (d) and (g) for the original  $J$ together with  $Z_{j-1} \subset Z_j$ and $z_j \notin Z_{j-1}$.

Next, let us check that the term in \eqref{e unj4} satisfies (a)--(g).
Condition (a) holds since  $z_j$ is a maximal element of  $Z_j$, (c) holds since  $z_j \notin Z_j$, and the remaining conditions are the same or weaker than those for the original  $J$.

Hence by induction, the expression in  \eqref{e unj3}--\eqref{e unj4} is congruent modulo  $\Iplacp{P}$ to
\begin{align}
&~\sum_{ T \in \SSYT_P^{\mathbf{Z}_*}(\alpha_-')}  \mathbf{u}_{\diagread(\e{z_jw};T)}
+ \sum_{ T \in \SSYT_P^{\mathbf{Z}_-}(\alpha')}  \mathbf{u}_{\diagread(\e{w};T)}.
\label{e diag rec2}
\end{align}
Finally, \eqref{e diag rec2} is equal to the right side of
\eqref{e claim diagread} by the same argument given in Case 1
(with indices  $j$ and  $j+1$ in place of  $1$ and  $2$ and  $z_j$ in place of  $z_1$).
\end{proof}

\begin{example}
\label{ex main proof}
Let $P = \PP_{2}$,  $\mathbf{Z} = (Z_1,Z_2) = ([3],[5]) = (\{1,2,3\}, \{1,2,3,4,  5\})$, and $\alpha = (1,2)$. Then
\[\SSYT_P^\mathbf{Z}(\alpha') = {\fontsize{8pt}{6pt}\selectfont  \tableau{1&1\\ \bl  & 4} \ \ \tableau{1&2\\ \bl  & 4} \ \ \tableau{1&1\\ \bl  & 5} \ \ \tableau{1&2\\ \bl  & 5} \ \ \tableau{1&3\\ \bl  & 5} \ \ \tableau{2&1\\ \bl  & 4} \ \ \tableau{2&2\\ \bl  & 4} \ \ \tableau{2&1\\ \bl  & 5} \ \ \tableau{2&2\\ \bl  & 5} \ \ \tableau{2&3\\ \bl  & 5} \ \ \tableau{3&2\\ \bl  & 4} \ \ \tableau{3&2\\ \bl  & 5} \ \ \tableau{3&3\\ \bl  & 5} }\]
and, according to Theorem \ref{t new statement fgprime}, summing the diagonal reading words of these tableaux gives
\[J_{1,2}^P([3],[5]) \equiv \mathbf{u}_\e{141} + \mathbf{u}_\e{142} + \mathbf{u}_\e{151} + \mathbf{u}_\e{152} + \mathbf{u}_\e{153} + \mathbf{u}_\e{241} + \mathbf{u}_\e{242} + \mathbf{u}_\e{251} + \mathbf{u}_\e{252} + \mathbf{u}_\e{253} + \mathbf{u}_\e{342} + \mathbf{u}_\e{352} + \mathbf{u}_\e{353}.\]

We will explain a step in the proof of Theorem \ref{t new statement fgprime} for the data above. Since $\alpha=(1,2)$, $j=2$. Then Case 2 applies and breaks  $J_{1,2}^P([3],[5])$ into two pieces (as in \eqref{e unj3}--\eqref{e unj4}):
\begin{align}
J_{1,2}^P([3],[5]) \equiv J_{1,1}^P([3]\Jnot{\e{5}}[3])  + J_{1,2}^P([3],[4]).
\end{align}
Here, we have chosen the maximal element $z_j = 5$ of  $Z_2\setminus Z_1$.
To clarify the definition of the $J$'s,
let's also write this line in terms of  $e^P_k$'s:
\begin{multline}
e_1^P(\mathbf{u}_{[3]})e_2^P(\mathbf{u}_{[5]})-e_2^P(\mathbf{u}_{[3]})e_1^P(\mathbf{u}_{[5]}) \equiv \\ e_1^P(\mathbf{u}_{[3]})\spa u_5 \spa  e_2^P(\mathbf{u}_{[3]})-e_2^P(\mathbf{u}_{[3]})\spa u_5 \spa e_1^P(\mathbf{u}_{[3]})
\, + \, e_1^P(\mathbf{u}_{[3]})e_2^P(\mathbf{u}_{[4]})- e_2^P(\mathbf{u}_{[3]})e_1^P(\mathbf{u}_{[4]}).
\end{multline}
Now, the terms $J_{1,1}^P([3]\Jnot{\e{5}}[3])$ and $J_{1,2}^P([3],[4])$ are given as follows:
the former is the sum over  $\mathbf{u}_{\diagread(\e{w};T)}$ with  $\e{w}= \e{5}$ and  $T$ ranging over
 $\SSYT_P^{\mathbf{Z}_*}(\alpha_-')$ where  $\mathbf{Z}_* =([3],[3])$ and  $\alpha_- = (1,1)$:
\begin{align}
\SSYT_P^{\mathbf{Z}_*}(\alpha_-') \ &= \ {\fontsize{8pt}{6pt}\selectfont  \tableau{1&1} \ \ \ \,  \tableau{1&2} \ \ \ \,  \tableau{1&3} \ \ \ \,  \tableau{2&1} \ \ \ \,  \tableau{2&2} \ \ \ \,  \tableau{2&3} \ \ \ \,  \tableau{3&2} \ \ \ \,  \tableau{3&3} }\\
J_{1,1}^P([3]\Jnot{\e{5}}[4]) \ &= \ \mathbf{u}_\e{151} + \mathbf{u}_\e{152} + \mathbf{u}_\e{153} + \mathbf{u}_\e{251} + \mathbf{u}_\e{252} + \mathbf{u}_\e{253} +\mathbf{u}_\e{352} + \mathbf{u}_\e{353}
\end{align}
The term $J_{1,2}^P([3],[4])$ is the sum over diagonal reading words of  $\SSYT_P^{\mathbf{Z_-}}(\alpha')$ where  $\mathbf{Z}_- = ([3],[4])$:
\begin{align}
\SSYT_P^{\mathbf{Z_-}}(\alpha') \  &= \ {\fontsize{8pt}{6pt}\selectfont  \tableau{1&1\\ \bl  & 4} \ \ \ \tableau{1&2\\ \bl  & 4} \ \ \ \tableau{2&1\\ \bl  & 4} \ \ \ \tableau{2&2\\ \bl  & 4} \ \ \ \tableau{3&2\\ \bl  & 4} } \\
J_{1,2}^P([3],[4]) \ &= \ \mathbf{u}_\e{141} + \mathbf{u}_\e{142} + \mathbf{u}_\e{241} + \mathbf{u}_\e{242} + \mathbf{u}_\e{342}.
\end{align}
By inspection, one sees that $\SSYT_P^\mathbf{Z}(\alpha')$ is naturally in bijection with  $\SSYT_P^{\mathbf{Z}_*}(\alpha_-') \sqcup \SSYT_P^{\mathbf{Z_-}}(\alpha')$, as is justified in general by the last sentence of the proof.
\end{example}

We conclude by reconciling the difference between the reading words in
Theorems \ref{th:IplacP-positivity} and~\ref{th:IplacP-positivity-technical}.

\begin{lemma}
\label{lem:poset basics}
\
\begin{itemize}
\item[(i)]  If  $a,b,c$ are (not necessarily distinct) elements of a poset  $P$ and  $a \to b \dotto c$, then  $a \dotto c$.
\item[(ii)] If  $a_1,a_2,a_3,b$ are (not necessarily distinct) elements of a  \threeone-free poset  $P$ with  $a_1 \to a_2 \to a_3$ and
$a_i \dotto b$ for all  $i =1,2,3$, then  $a_1 \to b$.
\end{itemize}
\end{lemma}
\begin{proof}
Statement (i) is an immediate consequence
of transitivity.  To see (ii), if  $a_1 \dote b$, then transitivity implies  $a_2 \dote b $ and $a_3 \dote b $, which would contradict  $P$ being \threeone-free.
\end{proof}

\begin{lemma}
\label{lem:T solid arrow}
Let  $P$ be a \threeone-free poset.
Let  $(r,c)$ and  $(r+d+1,c+d)$ be boxes of a partition diagram  $\lambda$ with  $d \ge 0$.
Then  $T(r,c) \to_P T(r+d+1,c+d)$ for any  $T \in \SSYT_P(\lambda)$.
\end{lemma}
\begin{proof}
The proof is by induction on  $d$.  The base case  $d=0$ holds since the columns of a  $P$-tableau read downwards are  $P$-strictly decreasing. Now let  $d > 0$.
Let  $a_1,a_2,a_3,b \in P$ be  $T(r,c), T(r+d, c+d-1), T(r+d+1, c+d-1), T(r+d+1,c+d)$, respectively.  Then  $a_1 \to a_2$ by induction,
while
$a_2 \to a_3$ 
and $a_3 \dotto b$
by definition of a $P$-tableau. Lemma \ref{lem:poset basics} (i) then gives $a_2 \dotto b$ and  $a_1 \dotto b$.
Hence $a_1 \to b$ by Lemma \ref{lem:poset basics} (ii).
\end{proof}

\begin{proposition}
\label{pr: diag vs col}
For any $P$-tableau $T$, $\mathbf{u}_{\creading(T)} \equiv \mathbf{u}_{\diagread(T)} \bmod \, \Iplacp{P}.$
Hence either the diagonal reading word or the column reading word can
be used in Theorems \ref{th:IplacP-positivity} and~\ref{th:IplacP-positivity-technical}.
\end{proposition}

\begin{proof}
The proof is by induction on the number of boxes of $T$.
Let $\e{d^1, d^2, \ldots, d^t}$ be the words corresponding to the diagonals of  $T$,
so that $\diagread(T) = \e{d^1 d^2 \cdots d^t}$.
Note that each word $\e{d^i}$ is  $P$-weakly increasing
by Lemma \ref{lem:poset basics} (i).
Let $r$ be the number of rows of $T$.
For each $i=1,2,\ldots,r$, set $\e{d^i = c_i\hat{d}^i}$ so that $\e{c_i}$ is the first letter of $\e{d^i}$ and  $\e{\hat{d}^i}$ is the remainder of the word $\e{d^i}$.
The key computation is
\begin{align}
\mathbf{u}_{\e{d^1d^2\cdots d^r}} &= \mathbf{u}_{\e{c_1\hat{d}^1c_2\hat{d}^2\cdots c_r\hat{d}^r}} \equiv \mathbf{u}_{\e{c_1\hat{d}^1c_2\hat{d}^2\cdots c_{r-1}c_r\hat{d}^{r-1}\hat{d}^r }} \equiv \cdots \equiv \mathbf{u}_{\e{c_1c_2\cdots c_r\hat{d}^1\cdots\hat{d}^r}} \ \bmod \Iplacp{P}, \label{e column diag}
\end{align}
where each congruence is by Lemma \ref{l plactic inc word} with  $s=1$; hypothesis (i) of the lemma holds since $c_1 \leftto c_2 \leftto \cdots \leftto c_r$, (ii) holds since each $\e{d^i}$ is  $P$-weakly increasing, and (iii) holds by Lemma \ref{lem:T solid arrow}.
Now let $T'$ be the  $P$-tableau obtained from $T$ by removing the first column.
Note that $\e{c_1\cdots c_r}$ is the first column of $T$ read bottom to top, and $\e{\hat{d}^1\cdots\hat{d}^rd^{r+1}\cdots d^t}$ is the diagonal reading word of $T'$.
Hence by \eqref{e column diag} and induction,
we have
$\mathbf{u}_{\diagread(T)}\equiv \mathbf{u}_{\e{c_1\cdots c_r}}\mathbf{u}_{\diagread(T')} \equiv \mathbf{u}_{\e{c_1\cdots c_r}}\mathbf{u}_{\creading(T')} \equiv \mathbf{u}_{\creading(T)} \,
\bmod \Iplacp{P}$.
\end{proof}

\section{
$m_\lambda^P(\mathbf{u})$ positivity for hook and two-column shapes}
\label{s two col hook}

We give the proofs of Theorems \ref{u pos hook} and \ref{t m two col}.

\subsection{The hook case}

Recall from Definition \ref{def:power word} that a word \emph{power word}
is a $P$-weakly increasing word with no nontrivial right-left $P$-minima.
Let  $\mathcal{N}_\ell$ denote the set of power words of length  $\ell$ and
 $\mathcal{E}_k$ the set of  $P$-strictly decreasing words of length  $k$.
The \emph{arm} of a  $P$-tableau  $T$ of hook shape $(\ell,1^k)$
is the word  $T(1,2)T(1,3) \cdots T(1,\ell)$  consisting
of the entries in the first row and columns  $2,\dots,\ell$ of  $T$.

\begin{proof}[Proof of Theorem \ref{u pos hook}]
The theorem states that for hook shapes  $\lambda$,
\begin{equation}
\label{eq:thm hook again}
\mathfrak{m}_{\lambda}^P(\mathbf{u}) \equiv \sum_{T \in \KT_P(\lambda)} \! \mathbf{u}_\creading(T) \ \ \bmod \, I_H^P.
\end{equation}
We will prove this by induction on  $N = |\lambda|$.  The base case  $N = 1$ is clear.  Now assume  $N > 1$.
We will use the following identities of ordinary symmetric functions:
\begin{align}
\label{eq:e times m 1}
 e_k(\mathbf{y})m_{(1)}(\mathbf{y}) & = m_{(2,1^{k-1})}(\mathbf{y}) + (k+1)m_{(1^{k+1})}(\mathbf{y}), \\
\label{eq:e times m 2}
 e_k(\mathbf{y})m_{(\ell)}(\mathbf{y}) & = m_{(\ell+1,1^{k-1})}(\mathbf{y}) + m_{(\ell,1^k)}(\mathbf{y}) \quad \quad \text{for  $\ell > 1$,  $k > 0$}.
\end{align}
Applying the map  $\psi$ from \eqref{eq:psi def} yields
\begin{align}
\label{eq:e times m 1b}
\!\!\!\!\!\!
& \mathfrak{m}^P_{(2,1^{k-1})}(\mathbf{u}) + (k+1)\mathfrak{m}^P_{(1^{k+1})}(\mathbf{u}) \equiv e^P_k(\mathbf{u})\mathfrak{m}^P_{(1)}(\mathbf{u})
\equiv \sum_{(\e{v},\e{w}) \in \mathcal{E}_k \times \mathcal{N}_1 }   \mathbf{u}_{\e{v}}\mathbf{u}_{\e{w}},
\end{align}
\begin{equation}
\begin{gathered}
\label{eq:e times m 2b}
\phantom{a} \!\!\!\!\!\!\!\!\!\! \!\!\!\!\!\!\!\!\!\!
 \mathfrak{m}^P_{(\ell+1,1^{k-1})}(\mathbf{u}) + \mathfrak{m}^P_{(\ell,1^k)}(\mathbf{u})
\equiv  e^P_k(\mathbf{u})\mathfrak{m}^P_{(\ell)}(\mathbf{u})
\equiv \sum_{(\e{v},\e{w}) \in \mathcal{E}_k \times \mathcal{N}_\ell }   \mathbf{u}_{\e{v}}\mathbf{u}_{\e{w}}
\quad \quad \quad \quad \quad
\\
\phantom{a} \quad \quad \quad \quad \quad \quad \quad \quad \quad \quad \quad \quad \quad \quad \quad \quad \ \ \ \ \ \ \ \ \
\text{for  $\ell > 1, k >0$,  $N = k+\ell$},
\end{gathered}
\end{equation}
where the congruences are modulo  $I_H^P$, and for the second congruence of \eqref{eq:e times m 2b} we
used the identity  $\mathfrak{m}^P_{(\ell)}(\mathbf{u})
\equiv \sum_{\e{w} \in \mathcal{N}_\ell }  \mathbf{u}_{\e{w}}$, which holds for  $\ell < N$ by induction.
We will connect these expressions to key  $P$-tableaux by defining, for each  $k,\ell > 0$, a map
\begin{align}
& \phi_{k,\ell}: \mathcal{E}_k \times \mathcal{N}_\ell \to \KT_P(\ell+1,1^{k-1}) \cup \KT_P(\ell,1^k)
\end{align}
satisfying (A) for $(\e{v},\e{w}) = (\e{v}_k \cdots \e{v}_1, \e{w}_1 \cdots \e{w}_\ell )\in \mathcal{E}_k \times \mathcal{N}_\ell$ we have $\e{v}\spa \e{w} \equiv \creading(\phi_{k,\ell}(\e{v}, \e{w}))$ mod  $I_H^P$, and (B) the preimage $\phi_{k,\ell}^{-1}(T)$ of a tableau $T$  has $k+1$ elements when $\ell = 1$ and $T$ has shape $1^{k+1}$, and 1 element otherwise.

The map $\phi_{k,\ell}(\e{v}, \e{w})$ is defined as follows:
\begin{itemize}
\item[(1)] If there is some $\e{v}_i$ so that $\e{v}_i\e{w}_1\e{w}_2\cdots \e{w}_\ell \in \mathcal{N}_{\ell+1}$, let $\phi_{k,\ell}(\e{v},\e{w})$ be the hook-shape $P$-tableau with first column $\e{v}$ and arm $\e{w}$. As $\e{v}_i\e{w}_1\e{w}_2\cdots \e{w}_\ell \in \mathcal{N}_{\ell+1}$, $\phi_{k,\ell}(\e{v},\e{w}) \in \KT_P(\ell+1,1^{k-1})$. Furthermore we have $\e{v}\spa \e{w} = \creading(\phi_{k,\ell}(\e{v}, \e{w}))$.
\item[(2)] Otherwise, there is an index  $i$ so that $\e{v}_1 \to_P \cdots \to_P \e{v}_i \to_P \e{w}_1 \to_P \e{v}_{i+1} \to_P \cdots \to_P \e{v}_k$ and $\e{v}_j$ is a right-left $P$-minima of $\e{v}_j\e{w}_1\e{w}_2\cdots \e{w}_\ell$ for each $j \le i$.
    Let $\phi_{k,\ell}(\e{v},\e{w})$ be the hook shape tableau with first column $\e{v}_1,\e{v}_2,\dots,\e{v}_i, \e{w}_1, \e{v}_{i+1},\dots,\e{v}_k$ and arm $\e{w}_2 \e{w}_3 \cdots \e{w}_\ell$. As $\e{w} \in \mathcal{N}_\ell$, $\phi_{k,\ell}(\e{v},\e{w})$ is in $\KT_P(\ell,1^{k})$. As comparable elements commute modulo $I_H^P$ and $\e{w}_1$ is comparable to each entry of $\e{v}$, we have $\e{v} \spa \e{w} \equiv \creading(\phi_{k,\ell}(\e{v},\e{w}))$.
\end{itemize}
It remains to show that the preimages of key $P$-tableaux have the desired number of elements. We begin with the case of a $P$-tableau $T$ of shape $1^{k+1}$. We want to show that $\phi_{k,1}^{-1}(T)$ has $k+1$ elements. Now if $\e{v}_1 \to_P \e{v}_2 \to_P \cdots \to_P \e{v}_{k+1}$ are the entries of $T$, then the elements of $\phi_{k,1}^{-1}(T)$ are $(\e{v}_{k+1}\e{v}_{k}\cdots \e{v}_{i+1}\e{v}_{i-1}\cdots \e{v}_{1}, \e{v}_i)$ for each $i \in [k+1].$

Now for $\ell > 1$ and a key $P$-tableau $T$ of shape $(\ell,1^k)$ we want to show that $\phi_{k,\ell}^{-1}(T)$ and $\phi_{k+1,\ell-1}^{-1}(T)$ each have one element.
Write $\e{v}_1 \to_P \e{v}_2 \to_P \cdots \to_P \e{v}_{k+1}$ for the entries of the first column of $T$ and $\e{w}_1 \dotto_P \e{w}_2 \dotto_P \cdots \dotto_P \e{w}_{\ell-1}$ for the entries of the arm of $T$. Then there is some smallest index $i$ so that $\e{v}_{i} \e{w}_1 \e{w}_2 \cdots \e{w}_{\ell-1}$ is a power word. Then $\phi_{k,\ell}^{-1}(T) =
(\e{v}_{k+1}\e{v}_{k}\cdots \e{v}_{i+1}\e{v}_{i-1}\cdots \e{v}_{1}, \e{v}_i \e{w}_1 \e{w}_2 \cdots \e{w}_{\ell-1})$ and $\phi_{k+1, \ell-1}^{-1}(T) =
(\e{v}_{k+1} \e{v}_k \cdots \e{v}_{1}, \e{w}_1 \e{w}_2 \cdots \e{w}_{\ell-1})$, as desired.
To see why the smallest index  $i$ is the correct choice here, note that if $\e{v}_{j} \e{w}_1 \e{w}_2 \cdots \e{w}_{\ell-1}$ is also a power word for some $j > i$, then $\e{v}_i \e{v}_{j}  \e{w}_1 \e{w}_2 \cdots \e{w}_{\ell-1}$ is a power word, so
$\phi_{k,\ell}$ applied to
$(\e{v}_{k+1}\e{v}_{k}\cdots \e{v}_{j+1}\e{v}_{j-1}\cdots \e{v}_{1}, \e{v}_j \e{w}_1 \e{w}_2 \cdots \e{w}_{\ell-1})$ is not  $T$
but rather the tableau with first column $\e{v}_1 \e{v}_2 \cdots \e{v}_{j-1} \e{v}_{j+1} \cdots \e{v}_{k+1}$
and arm $\e{v}_j \e{w}_1 \e{w}_2 \cdots \e{w}_{\ell-1}$.

Facts (A) and (B) about  $\phi_{k,\ell}$ and \eqref{eq:e times m 1b}--\eqref{eq:e times m 2b} now yield the following congruences modulo  $I_H^P$:
\begin{align}
\label{eq:e times m 1c}
& \mathfrak{m}_{(2,1^{N-2})}^P(\mathbf{u}) + N\mathfrak{m}_{(1^{N})}^P(\mathbf{u})
\equiv \sum_{T \in \KT_P(2,1^{N-2})} \!\!\! \mathbf{u}_{\creading(T)} + N \sum_{T\in \KT_P(1^{N})} \!\!\! \mathbf{u}_{\creading(T)}; \\
\label{eq:e times m 2c}
& \mathfrak{m}_{(\ell,1^k)}^P(\mathbf{u}) + \mathfrak{m}_{(\ell+1,1^{k-1})}^P(\mathbf{u})
\equiv \!\! \sum_{T \in \KT_P(\ell,1^k)\cup \KT_P(\ell+1,1^{k-1})} \!\!\!\!\! \mathbf{u}_{\creading(T)} \ \  \text{for  $\ell > 1$,  $k>0$,  $N = k+\ell$}.
\end{align}
Taking an alternating sum of \eqref{eq:e times m 1c} and  \eqref{eq:e times m 2c} over ranging $\ell$ and fixed  $N = k+\ell$,
\begin{align*}
& \mathfrak{m}_{(\ell,1^k)}^P(\mathbf{u}) +
 (-1)^{\ell-2}N\mathfrak{m}_{(1^{N})}^P(\mathbf{u})=\\
&
 (-1)^{\ell-2}\Big(\mathfrak{m}_{(2,1^{N-2})}^P(\mathbf{u}) + N\mathfrak{m}_{(1^{N})}^P(\mathbf{u})\Big) +
 \sum_{j=2}^{\ell-1} (-1)^{\ell-1-j}\Big(\mathfrak{m}_{(j,1^{N-j})}^P(\mathbf{u}) + \mathfrak{m}_{(j+1,1^{N-j-1})}^P(\mathbf{u})\Big)
 \\
& \equiv  \sum_{T \in \KT_P(\ell,1^k)} \mathbf{u}_{\creading(T)} + (-1)^{\ell-2} N
\sum_{T \in \KT_P(1^N)} \mathbf{u}_{\creading(T)}.
\end{align*}
Noting that $\mathfrak{m}_{1^N}^P(\mathbf{u}) = e_{N}^P(\mathbf{u}) = \sum_{T \in \KT_P(1^N) } \mathbf{u}_{\creading(T)}$, we obtain the formula \eqref{eq:thm hook again}.
\end{proof}

\subsection{The two-column case}

We freely use the definitions and results  on ladders from Definition \ref{def:ladder} and \S\ref{ss ladders}.
For a chain  $C = \big( c_k \leftto_P c_{k-1} \leftto_P \cdots \leftto_P c_1 \big)$, set  $\mathbf{u}_C = u_{c_k} \cdots u_{c_1}$.

\begin{lemma}
\label{lem:ladder swap IH}
Let $(A_1,A_2)$ be a pair of chains in  $P$ and let $(A'_1,A'_2)$ be the pair of chains obtained by
performing a ladder swap on an unbalanced ladder of $(A_1,A_2)$. Then
\begin{equation}
\mathbf{u}_{A_1}\mathbf{u}_{A_2} \equiv \mathbf{u}_{A'_1}\mathbf{u}_{A'_2}  \ \bmod \spa I_H^P.
\end{equation}
\end{lemma}
\begin{proof}
Let $L_\ell,\dots, L_2, L_1$ be the ladder decomposition of $(A_1,A_2)$. Suppose $(A'_1,A'_2)$ is obtained by performing a ladder swap on an unbalanced ladder $L_i$.
We have
\begin{align}
\mathbf{u}_{A_1} \mathbf{u}_{A_2}
& = \mathbf{u}_{A_1 \cap L_\ell} \cdots \mathbf{u}_{A_1 \cap L_2} \mathbf{u}_{A_1 \cap L_1}
\spa \mathbf{u}_{A_2 \cap L_\ell} \cdots \mathbf{u}_{A_2 \cap L_2}\mathbf{u}_{A_2 \cap L_1}.
\end{align}
As each element of $L_s$ is comparable to every element of $L_t$ when $s \neq t$, we have
\begin{align}
\label{eq:ladder swap far comm}
\mathbf{u}_{A_p \cap L_s} \mathbf{u}_{A_q \cap L_t} \equiv \mathbf{u}_{A_q \cap L_t}\mathbf{u}_{A_p \cap L_s}  \ \bmod \spa I_H^P
\end{align}
for $p,q \in \{1,2\}$ and $s \neq t$.
Therefore we can shuffle the ladders as follows:
\begin{equation}
\label{eq:ladder swap shuffle}
\mathbf{u}_{A_1} \mathbf{u}_{A_2}  \equiv \mathbf{u}_{A_1 \cap L_\ell} \mathbf{u}_{A_2 \cap L_\ell} \cdots \mathbf{u}_{A_1 \cap L_2} \mathbf{u}_{A_2 \cap L_2} \mathbf{u}_{A_1 \cap L_1} \mathbf{u}_{A_2 \cap L_1}.
\end{equation}
Next, by Lemma \ref{l swap chains}, $\mathbf{u}_{A_1 \cap L_i} \mathbf{u}_{A_2 \cap L_i} \equiv \mathbf{u}_{A_2 \cap L_i} \mathbf{u}_{A_1 \cap L_i}$ mod  $\Iplacp{P}$ and therefore mod  $I_H^P$ as well.  Applying this relation to the right side of \eqref{eq:ladder swap shuffle} and then unshuffling the ladders again using \eqref{eq:ladder swap far comm}, we obtain  $\mathbf{u}_{A'_1} \mathbf{u}_{A'_2}$.
 \end{proof}

We now state an equivalent definition of $P$-tableaux in terms of ladders. This is similar to
\cite[Theorem 5.11]{Ehrhard}.
A \emph{$P$-array} is a filling of a Young diagram with elements of $P$ so that each column, read bottom to top, is   $P$-strictly decreasing.

\begin{theorem}
\label{t two col ladder cond}
A  $P$-array  $T$ is a $P$-tableau if and only if, for each pair of adjacent columns $T_i,T_{i+1}$, the ladder
decomposition $L_\ell,\dots,L_1$ of  $(T_i,T_{i+1})$ satisfies
\begin{equation}
\label{ladder cond}
 \ \ \sum_{j=1}^k |L_j \cap T_i| \ge \sum_{j=1}^k |L_j \cap T_{i+1}|  \quad \!  \text{ for all $k \in [\ell]$}.
\end{equation}
\end{theorem}
\begin{proof}
Let $T$ be a $P$-array so that \eqref{ladder cond} is not satisfied for some columns $i$ and $i+1$. Let $k$ be the smallest positive integer so that
\begin{equation}
\sum_{j=1}^k |L_j \cap T_i| < \sum_{j=1}^k |L_j \cap T_{i+1}|,
\end{equation}
and let $r$ be the last row so that $T(r,i) \in \bigcup_{j=1}^k L_j$.
Then, since each $|L_j \cap T_i|$ and $|L_j \cap T_{i+1}|$ differ by at most 1, $r+1$ is the last row so that $T(r+1,i+1) \in \bigcup_{j=1}^k L_j$.
But then $T(r+1,i)$ must belong to a ladder with index greater than $k$, and since ladders are totally ordered,
$T(r+1,i) \leftto_P T(r+1,i+1)$. Thus
$T$ is not a $P$-tableau.

Now let $T$ be a  $P$-array that is not a $P$-tableau. Then for some columns $i$ and $i+1$ there is a row $r$ so that $T(r,i) \leftto_P T(r,i+1)$. Let $L_\ell,\dots,L_1$ be the ladder decomposition of columns $(T_i,T_{i+1})$.

First consider the case $T(r,i)$ and $T(r,i+1)$ lie in different ladders.
Then since ladders are totally ordered,
there is some $k$ so that $T(r, i+1) \in \bigcup_{j=1}^k L_j$ but $T(i,r) \not\in \bigcup_{j=1}^k L_j$. Therefore
\begin{equation}
\sum_{j=1}^k |L_j \cap T_i| \le r-1 < r \le \sum_{j=1}^k |L_j \cap T_{i+1}|,
\end{equation}
so \eqref{ladder cond} is not satisfied.

Now consider the case $T(r,i)$ and $T(r,i+1)$ lie in the same ladder $L_m$.
Then it follows from $T(r,i) \leftto_P T(r,i+1)$, transitivity, and Lemma \ref{l ladder} (ii) that the smallest row $a$ such that  $T(a,i) \in L_m$ and smallest row  $a'$ such that  $T(a', i+1) \in L_m$ satisfy  $a < a'$.
But this means $\sum_{j=1}^{m-1} |L_j \cap T_i| = a-1 < a'-1 = \sum_{j=1}^{m-1} |L_j \cap T_{i+1}|$,
so \eqref{ladder cond} is not satisfied in this case as well.
\end{proof}

\begin{lemma}
\label{l Kostka}
The Kostka coefficient $K_{2^a1^b, \, 2^{c}1^{d}}$ is the number of Yamanouchi words of length  $d$
consisting of $\frac{d+b}{2}$ 1's and $\frac{d-b}{2}$ 2's.
\end{lemma}
\begin{proof}
Let  $\SYT(\nu)$ denote the set of standard Young tableau of shape  $\nu$.
We have $K_{2^a1^b, \, 2^{c}1^{d}} = K_{2^{a-c}1^b, \, 1^{d}} = |\SYT(2^{a-c}1^b)| = |\SYT(\nu)|$ for 
$\nu = (2^{a-c}1^b)' = (a+b-c,a-c)=  (\frac{d+b}{2},\frac{d-b}{2})$.  This is the same as the number of Yamanouchi words of content  $\nu$ since
the map taking a word to its recording tableau
gives a bijection between Yamanouchi words of content $\nu$ and  $\SYT(\nu)$.
\end{proof}

\begin{proof}[Proof of Theorem \ref{t m two col}]
Let  $P$ be a \threeone-free poset.  For a two-column partition $\mu$, set
\begin{equation}
g_\mu^P(\mathbf{u}) = \sum_{T \in \LT_P(\mu)} \mathbf{u}_{\creading(T)}.
\end{equation}
The theorem states that $g_\mu^P(\mathbf{u}) \equiv \mathfrak{m}_\mu^P(\mathbf{u}) \spa \bmod \spa I_H^P$.
We will prove this by showing that
\begin{equation}
\label{eq: two col main}
\mathfrak{J}_\lambda^P(\mathbf{u}) \equiv \sum_{\mu} K_{\lambda \mu} \ g_\mu^P(\mathbf{u})  \ \ \bmod \,  I_H^P
\end{equation}
for any two-column partition  $\lambda$.

Let  $T$ be a $P$-tableau of shape  $\lambda= (2^a 1^b)$ with ladder decomposition  $L_\ell, \dots, L_1$ and unbalanced ladders  $L_{i_1}, \dots, L_{i_d}$ with  $i_1 > \cdots > i_d$.
Let  $v(T) = v_1 \cdots v_d$ be the word in which $v_j = 1$ (resp.  $v_j =2$) if  $L_{i_j}$ is left (resp. right) unbalanced.
By Theorem \ref{t two col ladder cond},  $v(T)$ is a Yamanouchi word.
Let  $L(T)$ denote the left  $P$-tableau obtained by swapping all right unbalanced ladders of  $T$
($L(T)$ is a  $P$-tableau and not just a  $P$-array, again by Theorem \ref{t two col ladder cond}).
Define a map
\begin{align}
\theta \colon \SSYT_P(\lambda) \to \bigsqcup_{\mu = 2^c 1^d} \LT_P(\mu) \times Y_{d}, \quad T \mapsto (L(T), v(T)), \end{align}
where the disjoint union is over all two-column shapes  $\mu$ which are less than or equal to  $\lambda$ in dominance order, and $Y_{d}$ denotes the set of Yamanouchi words of length  $d$ which consist of  $\frac{d+b}{2}$ 1's and $\frac{d-b}{2}$ 2's.
(Note that  $v(T)\in Y_d$ since  $d$ is the number of unbalanced ladders in a left  $P$-tableau of shape  $2^c 1^d$.)

It follows from Theorem \ref{t two col ladder cond} that  $\theta$ is a bijection.  Hence, since
 $|Y_d| = K_{\lambda \mu}$ by Lemma~\ref{l Kostka},
 Lemma \ref{lem:ladder swap IH} yields
\begin{equation}
\sum_{T\in \SSYT_P(\lambda)} \mathbf{u}_{\creading(T)} \equiv \sum_{\mu} K_{\lambda \mu} \ g_\mu^P(\mathbf{u}) \ \bmod \,  I_H^P.
\end{equation}
The left side is congruent modulo  $I_H^P$ to  $\mathfrak{J}_\lambda^P(\mathbf{u})$ by Theorem \ref{th:IplacP-positivity}, hence \eqref{eq: two col main} is proved.
\end{proof}

\bibliographystyle{plain}
\bibliography{mycitations}

\def\cprime{$'$} \def\cprime{$'$} \def\cprime{$'$}
\begin{thebibliography}{10}

\bibitem{APchromatic}
Per Alexandersson and Greta Panova.
\newblock L{LT} polynomials, chromatic quasisymmetric functions and graphs with
  cycles.
\newblock {\em Discrete Math.}, 341(12):3453--3482, 2018.

\bibitem{SamiOct13}
Sami~H. {Assaf}.
\newblock Dual equivalence graphs and a combinatorial proof of {LLT} and
  {M}acdonald positivity.
\newblock {\em {\tt arXiv:1005.3759v5}}, October 2013.

\bibitem{SamiForum}
Sami~H. Assaf.
\newblock Dual equivalence graphs {I}: {A} new paradigm for {S}chur positivity.
\newblock {\em Forum Math. Sigma}, 3:Paper No. e12, 33, 2015.

\bibitem{BLamLLT}
Jonah Blasiak.
\newblock Haglund's conjecture on 3-column {M}acdonald polynomials.
\newblock {\em Math. Z.}, 283(1-2):601--628, 2016.

\bibitem{BF}
Jonah Blasiak and Sergey Fomin.
\newblock Noncommutative {S}chur functions, switchboards, and {S}chur
  positivity.
\newblock {\em Selecta Math. (N.S.)}, 23(1):727--766, 2017.

\bibitem{BLKronecker}
Jonah Blasiak and Ricky~Ini Liu.
\newblock Kronecker coefficients and noncommutative super {S}chur functions.
\newblock {\em J. Combin. Theory Ser. A}, 158:315--361, 2018.

\bibitem{bc18}
Patrick Brosnan and Timothy~Y. Chow.
\newblock Unit interval orders and the dot action on the cohomology of regular
  semisimple {H}essenberg varieties.
\newblock {\em Adv. Math.}, 329:955--1001, 2018.

\bibitem{ChoHong}
Soojin Cho and Jaehyun Hong.
\newblock Positivity of chromatic symmetric functions associated with
  {H}essenberg functions of bounce number 3.
\newblock {\em Electron. J. Combin.}, 29(2):Paper No. 2.19, 37, 2022.

\bibitem{cho2019positivity}
Soojin Cho and JiSun Huh.
\newblock On e-positivity and e-unimodality of chromatic quasi-symmetric
  functions.
\newblock {\em SIAM Journal on Discrete Mathematics}, 33(4):2286--2315, 2019.

\bibitem{ChowDescents}
Timothy~Y. Chow.
\newblock Descents, quasi-symmetric functions, {R}obinson-{S}chensted for
  posets, and the chromatic symmetric function.
\newblock {\em J. Algebraic Combin.}, 10(3):227--240, 1999.

\bibitem{CHSS}
Samuel Clearman, Matthew Hyatt, Brittany Shelton, and Mark Skandera.
\newblock Evaluations of {H}ecke algebra traces at {K}azhdan-{L}usztig basis
  elements.
\newblock {\em Electron. J. Combin.}, 23(2):Paper 2.7, 56, 2016.

\bibitem{Dahlberg}
Samantha Dahlberg.
\newblock A new formula for {S}tanley's chromatic symmetric function for unit
  interval graphs and {$e$}-positivity for triangular ladder graphs.
\newblock {\em S\'{e}m. Lothar. Combin.}, 82B:Art. 59, 12, 2020.

\bibitem{dahlberg2018lollipop}
Samantha Dahlberg and Stephanie van Willigenburg.
\newblock Lollipop and lariat symmetric functions.
\newblock {\em SIAM Journal on Discrete Mathematics}, 32(2):1029--1039, 2018.

\bibitem{Ehrhard}
Henry {Ehrhard}.
\newblock {A Crystal Analysis of $P$-Arrays}.
\newblock {\em {\tt arXiv:2205.01834}}, May 2022.

\bibitem{FG}
Sergey Fomin and Curtis Greene.
\newblock Noncommutative {S}chur functions and their applications.
\newblock {\em Discrete Math.}, 193(1-3):179--200, 1998.
\newblock Selected papers in honor of Adriano Garsia (Taormina, 1994).

\bibitem{Gallai}
T.~Gallai.
\newblock Transitiv orientierbare {G}raphen.
\newblock {\em Acta Math. Acad. Sci. Hungar.}, 18:25--66, 1967.

\bibitem{Gasharov}
Vesselin Gasharov.
\newblock Incomparability graphs of $(3+1)$-free posets are $s$-positive.
\newblock {\em Discrete Math.}, 157:193--197, 1996.

\bibitem{gebhard2001chromatic}
David~D Gebhard and Bruce~E Sagan.
\newblock A chromatic symmetric function in noncommuting variables.
\newblock {\em Journal of Algebraic Combinatorics}, 13(3):227--255, 2001.

\bibitem{GesselPPartition}
Ira~M. Gessel.
\newblock Multipartite {$P$}-partitions and inner products of skew {S}chur
  functions.
\newblock In {\em Combinatorics and algebra ({B}oulder, {C}olo., 1983)},
  volume~34 of {\em Contemp. Math.}, pages 289--317. Amer. Math. Soc.,
  Providence, RI, 1984.

\bibitem{GKcyl}
Ira~M. Gessel and C.~Krattenthaler.
\newblock Cylindric partitions.
\newblock {\em Trans. Amer. Math. Soc.}, 349(2):429--479, 1997.

\bibitem{GouldenJackson}
I.~P. Goulden and D.~M. Jackson.
\newblock Immanants of combinatorial matrices.
\newblock {\em J. Algebra}, 148(2):305--324, 1992.

\bibitem{GreeneImmanant}
Curtis Greene.
\newblock Proof of a conjecture on immanants of the {J}acobi-{T}rudi matrix.
\newblock {\em Linear Algebra Appl.}, 171:65--79, 1992.

\bibitem{GPchromatic}
Mathieu Guay-Paquet.
\newblock A modular relation for the chromatic symmetric functions of
  (3+1)-free posets.
\newblock {\em { \tt arxiv:1306.2400}}, 2013.

\bibitem{g16}
Mathieu Guay-Paquet.
\newblock A second proof of the {S}hareshian--{W}achs conjecture, by way of a
  new {H}opf algebra.
\newblock {\em {\tt arxiv:1601.05498}}, 2016.

\bibitem{Himmanant}
Mark Haiman.
\newblock Hecke algebra characters and immanant conjectures.
\newblock {\em J. Amer. Math. Soc.}, 6(3):569--595, 1993.

\bibitem{harada2019cohomology}
Megumi Harada and Martha~E. Precup.
\newblock The cohomology of abelian {H}essenberg varieties and the
  {S}tanley-{S}tembridge conjecture.
\newblock {\em Algebr. Comb.}, 2(6):1059--1108, 2019.

\bibitem{Hwang}
Byung-Hak Hwang.
\newblock Chromatic quasisymmetric functions and noncommutative {$P$}-symmetric
  functions.
\newblock {\em {\tt arXiv:2208.09857}}, August 2022.

\bibitem{KPchromatic}
Dongkwan Kim and Pavlo Pylyavskyy.
\newblock Robinson-{S}chensted correspondence for unit interval orders.
\newblock {\em Selecta Math. (N.S.)}, 27(5):Paper No. 97, 66, 2021.

\bibitem{LamRibbon}
Thomas Lam.
\newblock Ribbon {S}chur operators.
\newblock {\em European J. Combin.}, 29(1):343--359, 2008.

\bibitem{LS}
Alain Lascoux and Marcel-P. Sch{\"u}tzenberger.
\newblock Le mono\"\i de plaxique.
\newblock In {\em Noncommutative structures in algebra and geometric
  combinatorics ({N}aples, 1978)}, volume 109 of {\em Quad. ``Ricerca Sci.''},
  pages 129--156. CNR, Rome, 1981.

\bibitem{Macdonald95}
I.~G. Macdonald.
\newblock {\em Symmetric functions and {H}all polynomials}.
\newblock The Clarendon Press, Oxford University Press, New York, second
  edition, 1995.
\newblock With contributions by A.~Zelevinsky, Oxford Science Publications.

\bibitem{McNamara}
Peter McNamara.
\newblock Cylindric skew {S}chur functions.
\newblock {\em Adv. Math.}, 205(1):275--312, 2006.

\bibitem{PostnikovCylSchur}
Alexander Postnikov.
\newblock Affine approach to quantum {S}chubert calculus.
\newblock {\em Duke Math. J.}, 128(3):473--509, 2005.

\bibitem{Sch}
M.-P. Sch{\"u}tzenberger.
\newblock La correspondance de {R}obinson.
\newblock In {\em Combinatoire et repr\'esentation du groupe sym\'etrique
  ({A}ctes {T}able {R}onde {CNRS}, {U}niv. {L}ouis-{P}asteur {S}trasbourg,
  {S}trasbourg, 1976)}, pages 59--113. Lecture Notes in Math., Vol. 579.
  Springer, Berlin, 1977.

\bibitem{SWchromatic}
John Shareshian and Michelle~L. Wachs.
\newblock Chromatic quasisymmetric functions.
\newblock {\em Adv. Math.}, 295:497--551, 2016.

\bibitem{Isaiah}
Isaiah {Siegl}.
\newblock {Cylindric $P$-tableaux for 3+1-free posets}.
\newblock {\em {\tt arXiv:2211.03953}}, November 2022.

\bibitem{StanleyPpartition}
Richard~P. Stanley.
\newblock {\em Ordered structures and partitions}.
\newblock Memoirs of the American Mathematical Society, No. 119. American
  Mathematical Society, Providence, R.I., 1972.

\bibitem{Stanleychromatic}
Richard~P. Stanley.
\newblock A symmetric function generalization of the chromatic polynomial of a
  graph.
\newblock {\em Adv. Math.}, 111(1):166--194, 1995.

\bibitem{StanleyStembridge}
Richard~P. Stanley and John~R. Stembridge.
\newblock On immanants of {J}acobi-{T}rudi matrices and permutations with
  restricted position.
\newblock {\em J. Combin. Theory Ser. A}, 62(2):261--279, 1993.

\bibitem{StembridgeImmanant}
J.~R. Stembridge.
\newblock Some conjectures for immanants.
\newblock {\em Canad. J. Math.}, 44(5):1079--1099, 1992.

\bibitem{SWW}
Thomas~S. Sundquist, David~G. Wagner, and Julian West.
\newblock A {R}obinson-{S}chensted algorithm for a class of partial orders.
\newblock {\em J. Combin. Theory Ser. A}, 79(1):36--52, 1997.

\bibitem{Wolfgangthesis}
Harry~Lewis Wolfgang, III.
\newblock {\em Two interactions between combinatorics and representation
  theory: {M}onomial immanants and {H}ochschild cohomology}.
\newblock ProQuest LLC, Ann Arbor, MI, 1997.
\newblock Thesis (Ph.D.)--Massachusetts Institute of Technology.

\end{thebibliography}

\end{document}